%% file: BFV-complex.tex
\documentclass[11pt]{amsart}
\usepackage{amsmath}
\usepackage{amssymb}
\usepackage{amsfonts}
\usepackage{amsthm}
\usepackage{hyperref}
\usepackage{latexsym}
\usepackage{graphics}
\usepackage{graphicx}
\usepackage[all]{xy}
\usepackage{verbatim}

\newtheorem{Definition}{Definition}
\newtheorem{Lemma}{Lemma}

\newtheorem{Theorem}{Theorem}
\newtheorem{Proposition}{Proposition}
\newtheorem{Corollary}{Corollary}

\theoremstyle{definition}
\newtheorem*{Ack}{Acknowledgement}

\newcommand{\N}{\mathbb N}

\DeclareMathOperator{\sign}{sign}
\DeclareMathOperator{\dec}{dec}

\title{\bf BFV-complex and higher homotopy structures}

\author[F. Sch\"atz]{Florian Sch\"atz}
\address{Institut f\"ur Mathematik, Universit\"at Z\"urich--Irchel, Winterthurerstrasse 190, CH-8057 Z\"urich, Switzerland}
\email{florian.schaetz@math.uzh.ch}

\thanks{The author acknowledges partial support by the joint graduate school of mathematics of the ETH and the University of Z\"urich, by SNF-grant Nr.20-113439,
by the European Union through the FP6 Marie Curie RTN ENIGMA (contract number MRTN-CT-2004-5652), and by the European Science Foundation through the MISGAM program. Moreover he thanks the ESI for Mathematical Physics
for support during the author's visit in July and August 2007.}

\begin{document}

\maketitle

\begin{abstract}

We present a connection between the BFV-complex (appreviation for Batalin-Fradkin-Vilkovisky complex) and the so-called strong homotopy Lie algebroid associated
to a coisotropic submanifold of a Poisson manifold. We prove that the latter structure can be derived from the BFV-complex by means of homotopy transfer along contractions. Consequently the BFV-complex and the strong homotopy Lie algebroid structure
are $L_{\infty}$ quasi-isomorphic and control the same formal deformation problem.

However there is a gap between the non-formal
information encoded in the BFV-complex and in the strong homotopy Lie algebroid respectively. We prove that there is a one-to-one correspondence between coisotropic submanifolds given by graphs of sections and equivalence classes of normalized Maurer-Cartan elemens of the BFV-complex. This does not hold if one uses the strong homotopy Lie algebroid instead.

\end{abstract}

\tableofcontents

\begin{Ack} I thank A. Cattaneo for many encouraging and inspiring discussions and his general support.
I also thank D. Fiorenza, D. Indelicato,
B. Keller, P. Mn\"ev, T. Preu,  C. Rossi, S. Shadrin, J. Stasheff and M. Zambon for clarifying
discussions and helpful remarks on a draft of this paper. Moreover I thank
M. Bordemann and H.-C. Herbig
for pointing me to their globalisation of the BFV-complex. The referee contributed a lot to the
form of this work with his/her insightful suggestions.
\end{Ack}

\section{Introduction}

The geometry of coisotropic submanifolds inside Poisson manifolds is a very rich subject with connections to topics such as folitation theory, momentum maps, constrained systems and symplectic groupoids -- see \cite{Weinstein} for instance. Recently a new algebraic structure called the ``strong homotopy Lie algebroid'' associated to such submanifolds has been investigated, e.g. \cite{OhPark} in the symplectic setting or \cite{CattaneoFelder} in the Poisson case. This structure is related to the deformation problem of a given coisotropic submanifold (\cite{OhPark}) on the one hand and to the quantization of constrained systems (\cite{CattaneoFelder}) on the other. Moreover it captures subtle properties of the foliation associated to a coisotropic submanifold (\cite{Kieserman}).

The first main result of this paper is to reveal that the strong homotopy Lie algebroid is in some sense equivalent to a construction known as the BFV-complex -- for a precise formulation see Theorem \ref{thm:big} in subsection \ref{s:main}. The BFV-complex originated from physical considerations concerning the quantization of field theories with so-called open gauge symmetries (\cite{BatalinFradkin}, \cite{BatalinVilkovisky}). It was given an interpretation in terms of homological algebra in \cite{Stasheff} and globalized to coisotropic submanifolds of arbitrary finite dimensional Poisson manifolds in \cite{Bordemann} and \cite{Herbig}. 

Theorem \ref{thm:big} provides a connection between the BFV-complex and the strong homotopy Lie algebroid. In fact, we show that the two structures are isomorphic up to homotopy. In particular this implies (Corollary \ref{cor:formaldef}) that the formal deformation problem associated to both structures is equivalent. In \cite{OhPark} this formal deformation problem was investigated in the setting of the strong homotopy Lie algebroid (in the symplectic case).

Remarkably there is a gap between the strong homotopy Lie algebroid and the BFV-complex in the non-formal regime: we present a simple example of a coisotropic submanifold inside a Poisson manifold where the strong homotopy Lie algebroid does not capture obstructions to deformations. However the BFV-complex always does, see Theorem \ref{thm:BFV_MC} in subsection \ref{s:(n)MC} for the precise statement. Hence the BFV-complex is able to capture non-formal aspects of the geometry of coisotropic submanifolds. This is also supported by the example considered in subsection \ref{s:example_Marco} where the treatment using the BFV-complex reproduces a criterion for finding coisotropic submanifolds which was derived in \cite{Zambon}.

The paper is organized as follows: Section \ref{s:preliminaries} collects known facts concerning algebraic and geometric structures that are used in the main body of the paper. In section \ref{s:BFV} we present the global construction of the BFV-complex. We mainly follow \cite{Stasheff}, \cite{Bordemann} and \cite{Herbig} there. The only original part is the conceptual construction of the global BFV-bracket (see subsection \ref{s:lifting}). Section \ref{s:Liealgebroid} introduces the strong homotopy Lie algebroid and connects it to the BFV-complex (Theorem \ref{thm:big}). In section \ref{s:deformation} we establish a link between the BFV-complex and the geometry of coisotropic sections (Theorem \ref{thm:BFV_MC}) and give an example to demonstrate that this link does not exist if one considers the strong homotopy Lie algebroid instead. In the Appendix we give details on the homotopy transfer along contraction data which is one of our main tools. The material there is well known to the experts.

\section{Preliminaries}\label{s:preliminaries}
\input{pre.tex}

\section{The BFV-complex}\label{s:BFV}
\input{BFV.tex}

\section{Connection to the strong homotopy Lie Algebroid}\label{s:Liealgebroid}
\input{Liealgebroid.tex}

\section{The Deformation Problem}\label{s:deformation}
\input{Deformation.tex}

\begin{appendix}
\input{Appendix.tex}
\end{appendix}

\end{document}

%% file: pre.tex
For the convenience of the reader and in order to fix conventions we recall some basic definitions and facts concerning $L_{\infty}$-algebras (subsection \ref{s:infty}), the derived bracket formalism (subsection \ref{s:derivedbrackets}), homotopy transfer of $L_{\infty}$-algebras along contraction data (subsection \ref{s:homotopytransfer}), smooth graded manifolds (subsection \ref{s:sgmfs}) and Poisson geometry (subsection \ref{s:Poisson}). Readers familiar to these topics might skip this section.

\subsection{$L_{\infty}$-algebras}\label{s:infty}

Let $V$ be a $\mathbb{Z}$-graded vector space over $\mathbb{R}$ (or any other field of characteristic $0$); i.e.,
$V$ is a collection $\{V_{i}\}_{i \in \mathbb{Z}}$ of vector spaces $V_{i}$ over $\mathbb{R}$.
Homogeneous elements of $V$ of degree $i\in \mathbb{Z}$ are the elements of $V_{i}$. We denote the degree
of a homogeneous element $x \in V$ by $|x|$.  A morphism $f: V \to W$ of graded vector
spaces is a collection $\{f_{i}: V_{i} \to W_{i}\}_{i \in \mathbb{Z}}$ of
linear maps.
The $n$th suspension functor $[n]$ from the category of graded vector spaces to itself is defined as follows: given a graded vector space $V$,
$V[n]$ denotes the graded vector space given by the collection $V[n]_{i}:=V_{n+i}$.
The $n$th supsension of a morphism $f: V \to W$ of graded vector spaces is given
by the collection $\{f[n]_{i}:=f_{n+i}:V_{n+i} \to W_{n+i}\}_{i \in \mathbb{Z}}$.

One can consider the tensor algebra $T(V)$ associated to a graded vector space $V$ which is a graded vector space
with components
\begin{align*}
T(V)_{m}:=\bigoplus_{k\ge 0}\bigoplus_{j_{1}+\cdots +j_{k}=m}V_{j_{1}}\otimes \cdots \otimes V_{j_{k}}.
\end{align*}
$T(V)$ naturally carries the structure of a cofree coconnected coassociative coalgebra given by the 
deconcatenation coproduct:
\[
\Delta(x_{1}\otimes \dots \otimes x_{n}):=\sum_{i=0}^{n}(x_{1}\otimes \dots \otimes x_{i})\otimes(x_{i+1}\otimes \cdots \otimes x_{n}).
\]
There are two natural representations of the symmetric group $\Sigma_{n}$ on $V^{\otimes n}$: the even one which
is defined by multiplication with the sign $(-1)^{|a||b|}$ for the transposition interchanging $a$ and $b$ in $V$ 
and the odd one by multiplication with the sign $-(-1)^{|a||b|}$ respectively. These two actions naturally 
extend to $T(V)$.
The fix point set of the first action on $T(V)$ is 
denoted by $S(V)$
and called the graded symmetric algebra of $V$ while the fix point set of the latter action is denoted by $\Lambda(V)$
and called the graded skew--symmetric algebra of $V$. The graded symmetric algebra
$S(V)$ inherits a coalgebra structure from $T(V)$ which is cofree
coconnected coassociative and graded cocommutative.

Let $V$ be a graded vector space together with a family of linear maps
\begin{align*}
(m^{n}\colon S^{n}(V) \to V[1])_{n\in \N}.
\end{align*}
Given such a family one defines the associated family of Jacobiators
\begin{align*}
(J^{n}\colon S^{n}(V) \to V[2])_{n \ge 1}
\end{align*}
by
\begin{multline*}
\label{Jacobiators}
J^{n}(x_{1} \cdots x_{n}):=\\=
\sum_{r+s=n} \sum_{\sigma \in (r,s)-\text{shuffles}} \mspace{-36mu} \sign(\sigma)\, m^{s+1}(m^{r}(x_{\sigma(1)} \otimes \cdots \otimes x_{\sigma(r)})
\otimes x_{\sigma(r+1)} \otimes \cdots \otimes x_{\sigma(n)}) 
\end{multline*}
where $\sign(\cdot)$ is the Koszul sign, i.e., the one induced from the natural even representation of $\Sigma_{n}$ on $T^{n}(V)$, and
$(r,s)$-shuffles are permutations $\sigma$ of $\{1,\dots,n\}$ such that $\sigma(1) < \cdots < \sigma(r)$ and
$\sigma(r+1) < \cdots < \sigma(n)$.

\begin{Definition}
\label{L-infty}
A family of maps $(m^{n}\colon S^{n}(V) \to V[1])_{n\in \N}$ defines the structure of an \textsf{$L_\infty[1]$-algebra} 
on the graded vector space $V$ whenever the associated family of Jacobiators vanishes identically.
\end{Definition}

This definition is the one given in \cite{Voronov}.
We remark that this definition deviates from
the classical notion of $L_\infty$-algebras (see \cite{LadaStasheff} for instance) in two points. First it makes use of the graded symmetric algebra over $V$ instead of the the graded
skew--symmetric one. The transition between
these two settings uses the so called d\'ecalage-isomorphism
\[
\dec^{n}\colon \begin{array}[t]{ccc}
S^{n}(V) &\to& \Lambda^{n}(V[-1])[n]\\
x_{1} \cdots x_{n} &\mapsto& (-1)^{\sum_{i=1}^{n}(n-i)(|x_{i}|)}x_{1} \wedge \cdots \wedge x_{n}.
\end{array}
\]
The connection between $L_{\infty}[1]$-algebras and $L_{\infty}$-algebras is easy:

\begin{Lemma}
Let $W$ be a graded vector space.
There is a one-to-one correspondence between $L_{\infty}[1]$-algebra structures on $W[1]$ and
$L_{\infty}$-algebra structures on $W$.
\end{Lemma}
%\label{decalage}
More important is the fact that we also allow a map $m_{0}: \mathbb{R} \to V[1]$ as part of the structure given by an
$L_\infty[1]$-algebra. This piece can be interpreted as an element of $V_{1}$. In the traditional definition
$m_{0}$ is assumed to vanish. Relying on a widespread terminology,
we call structures with $m_{0}=0$ ``flat''. Observe that in the case of a flat $L_\infty[1]$-algebra $m_1$ is a differential. Moreover $L_{\infty}[1]$-algebras with $m_{k}=0$ for all $k\neq 1,2$
correspond exactly to differential graded Lie algebras under the d\'ecalage-isomorphism:

\begin{Definition}
\label{DGLA}
A \textsf{graded Lie algebra} $(\mathfrak{h},[-,-])$ is
a graded vector space $\mathfrak{h}$ equipped with a linear map $[-,-]\colon \mathfrak{h}\otimes \mathfrak{h} \to \mathfrak{h}$ 
satisfying the following conditions:
\begin{itemize}
\item graded skew-symmetry: $[x,y]=-(-1)^{|x||y|}[y,x]$ and
\item graded Jacobi identity: $[x,[y,z]]=[[x,y],z]+(-1)^{|x||y|}[y,[x,z]]$,
\end{itemize}
for all homogeneous $x \in \mathfrak{h}_{|x|}$, $y \in \mathfrak{h}_{|y|}$ and $z \in \mathfrak{h}$.

A \textsf{differential graded Lie algebra} is a triple $(\mathfrak{h},d,[-,-])$ where $(\mathfrak{h},[-,-])$ is a graded Lie algebra and $d$ is a linear map of degree $+1$ such that
$d \circ d =0$ and $d[x,y]=[dx,y]+(-1)^{|x|}[x,dy]$ holds for all $x\in \mathfrak{h}_{|x|}$ and $y \in \mathfrak{h}$.
\end{Definition}

If one goes from the category of graded vector spaces to the category of graded commutative associative algebras,
the reasonable replacement of the notion of (differential) graded Lie algebra is that of a
(differential) graded Poisson algebra: 

\begin{Definition}
\label{DGPA}
A \textsf{graded Poisson algebra} is a triple $(A,\cdot,[-,-])$ where $(A,\cdot)$ is a graded
commutative associative algebra and $(A,[-,-])$ is a graded Lie algebra such that $[x,y\cdot z]=
[x,y]\cdot z+(-1)^{|x||y|}y\cdot[x,z]$ holds for $x\in A_{|x|}$, $y\in A_{|y|}$ and $z\in A$.

A \textsf{differential graded Poisson algebra} is a quadruple $(A,d,\cdot,[-,-])$ where
$(A,\cdot,[-,-])$ is a graded Poisson algebra, $(A,d,[-,-])$ is a differential graded Lie algebra and
$d(x\cdot y)=dx \cdot y + (-1)^{|x|} x\cdot dy$ holds for all $x\in A_{|x|}$ and $y\in A$.
\end{Definition}

We briefly review a description of $L_\infty[1]$-algebras, equivalent to the one given in Definition \ref{L-infty}, which
goes back to Stasheff \cite{Stasheff1}. We remarked before that the graded commutative algebra $S(V)$
associated to a graded vector space $V$
is a cofree coconnected graded cocommutative coassociative coalgebra with respect to the coproduct $\Delta$
inherited from $T(V)$.
A linear map $Q\colon S(V) \to S(V)$ that satisfies $\Delta \circ Q= (Q \otimes \text{id} + \text{id} \otimes Q)\circ \Delta$ is called a \textsf{coderivation}
of $S(V)$. By cofreeness of the coproduct $\Delta$ it follows that every linear map from $S(V)$ to $V$
can be extended to a coderivation of $S(V)$ and that every coderivation $Q$ is uniquely determined
by $\text{pr}\circ Q$ where $\text{pr}\colon S(V)\to V$ is the natural projection. So there is a one-to-one
correspondence between families of linear maps $(m^{n}\colon S^{n}(V) \to V[1])_{n\in \mathbb{N}}$ 
and coderivations of $S(V)$ of degree $1$. Moreover, the graded commutator equips
$\oplus_{k\in \mathbb{Z}}\text{Hom}(S(V),S(V)[k])[-k]$ with the structure of a graded Lie algebra
and this Lie bracket restricts to the subspace of coderivations of $S(V)$. Odd coderivations $Q$ that satisfy $[Q,Q]=0$
are in one-to-one correspondence with families of maps whose associated Jacobiators
vanish identically. Consequently, Maurer-Cartan elements of the space of coderivations of $S(V)$ correspond
exactly to $L_\infty[1]$-algebra structures on $V$. Since
$Q \circ Q =\frac{1}{2}[Q,Q]=0$, Maurer--Cartan elements of the space of coderivations are exactly the
codifferentials of $S(V)$.

We remark that the approach to $L_\infty[1]$-algebras outlined above makes the notion of $L_\infty[1]$-morphisms
especially transparent: these are just coalgebra morphisms that are
chain maps between the graded symmetric algebras equipped with the codifferentials
that define the $L_\infty[1]$-algebra structures. There are two special kinds of
$L_{\infty}[1]$-morphisms. As usual $L_{\infty}[1]$-isomorphisms are
$L_{\infty}[1]$-morphism with an inverse. Moreover there is the notion of $L_{\infty}[1]$
quasi-isomorphisms, i.e. those $L_{\infty}[1]$-morphisms which admit ``inverses up to homotopy'': consider an
$L_{\infty}[1]$-morphism between flat $L_{\infty}[1]$-algebras, hence the unary structure
maps are coboundary operators. The given $L_{\infty}[1]$-morphism also has a unary component which is a chain map for these coboundary operators. Consequently this map induces
a map between the cohomologies. An $L_{\infty}[1]$ quasi-isomorphism is an
$L_{\infty}[1]$-morphism between flat $L_{\infty}$-algebras such that this induced map between cohomologies is an isomorphism.
The notions of $L_{\infty}$-morphisms, isomorphisms and quasi-isomorphisms are obtained from the corresponding
notions in the category of $L_{\infty}[1]$-algebras using the identification under the d\'ecalage-isomorphism.

Associated to every $L_{\infty}$-algebra structure $(m^{n}:\bigwedge^{n}(V) \to V[2-n])_{n\in \mathbb{N}}$ on a graded vector space $V$ is a subset of $V_{1}$ given by the zero set of the so called \textsf{MC-equation} (MC stands for Maurer-Cartan from now on) which reads
\begin{align*}
\sum_{n \ge 0}\frac{1}{n!}m^{n}(\mu\otimes \cdots \otimes \mu)=0.
\end{align*}
Elements of $V_{1}$ satisfying this equation are called \textsf{MC-elements}. We denote the set of all these elements by $MC(V)$. It is well-known
that there is a natural action of $V_{0}$ on $V$ by inner derivations. Integrating these one obtains a subgroup $Inn(V)$
of the automorphism group $Aut(V)$ of the $L_{\infty}$-algebra $V$. There is an induced action on
$MC(V)$.
We will give a complete definition of the action of $V_{0}$ on $MC(V)$ for $V$ being the BFV-complex in subsection \ref{s:(n)MC}.

\subsection{Derived Brackets Formalism}\label{s:derivedbrackets}

We describe the derived brackets formalism essentially following \cite{Voronov}.

\begin{Definition}
\label{V-algebra}
We call the triple $(\mathfrak{h},\mathfrak{a},\Pi_{\mathfrak{a}})$ a \textsf{V-algebra} (V for Voronov)
if $(\mathfrak{h},[\cdot,\cdot])$ is a graded Lie algebra, $\mathfrak{a}$ is an abelian Lie subalgebra of $\mathfrak{h}$ -- i.e.
$\mathfrak{a}$ is a graded vector subspace of $\mathfrak{h}$ and $[\mathfrak{a},\mathfrak{a}]=0$ -- and
$\Pi_{\mathfrak{a}}\colon \mathfrak{h} \to \mathfrak{a}$ is a projection such that
\begin{align}
\Pi_{\mathfrak{a}}[x,y]=\Pi_{\mathfrak{a}}[\Pi_{\mathfrak{a}}x,y]+\Pi_{\mathfrak{a}}[x,\Pi_{\mathfrak{a}}y]
\label{projection}
\end{align}
holds for every $x, y \in \mathfrak{h}$.
\end{Definition}

Instead of condition \eqref{projection} 
one can require that $\mathfrak{h}$ splits
into $\mathfrak{a}\oplus \mathfrak{p}$ as a graded vector space where $\mathfrak{p}$ is also a graded Lie subalgebra of $\mathfrak{h}$.
In terms of the projection, $\mathfrak{p}$ is given by the kernel of $\Pi_{\mathfrak{a}}$.

Let $(\mathfrak{h},\mathfrak{a},\Pi_{\mathfrak{a}})$ be a V-algebra and pick an element
$P \in \mathfrak{h}$ of degree $+1$. One can define the multilinear maps on $\mathfrak{a}$
\begin{equation}\label{derivedbrackets}
D_{P}^{n}\colon\begin{array}[t]{ccc} \mathfrak{a}^{\otimes n}&\to& \mathfrak{a}[1]\\
x_{1}\otimes \dots \otimes x_{n} &\mapsto& \Pi_{\mathfrak{a}}[[\dots[[P,x_{1}],x_{2}],\dots ],x_{n}]
\end{array}
\end{equation}
for every $n \ge 0$. These maps are called the \textsf{higher derived brackets} associated to $P$.
It is easy to check that all these maps are graded commutative, namely:
\begin{multline*}
D_{P}^{n}(x_{1}\otimes \dots \otimes x_{i} \otimes x_{i+1} \otimes \dots \otimes x_{n})=\\
=(-1)^{|x_{i}||x_{i+1}|}D_{P}^{n}(x_{1}\otimes \dots \otimes x_{i+1} \otimes x_{i} \otimes \dots \otimes x_{n})
\end{multline*}
for every $1 \le i \le n-1$.
We restrict the higher derived brackets constructed from $P$
to the symmetric algebra $S(\mathfrak{a})$ and obtain a family of maps 
$(D_{P}^{n}\colon S^{n}(\mathfrak{a}) \to \mathfrak{a}[1])_{n \in \mathbb{N}}$.

In \cite{Voronov} it is proven that the Jacobiators of the higher derived brackets
$(D_{P}^{n}\colon S^{n}(\mathfrak{a})\to \mathfrak{a}[1])_{n\in \mathbb{N}}$ associated
to $P$
are given by the higher derived brackets associated to $\frac{1}{2}[P,P]$:
\begin{align*}
J_{D_{P}}^{n}=D_{\frac{1}{2}[P,P]}^{n}.
\end{align*}
{}It follows that all Jacobiators vanish identically
if we assume that $[P,P]=0$ holds. Elements $P$ of degree $1$ that satisfy $[P,P]=0$ are exactly the MC-elements
of the graded Lie algebra $\mathfrak{h}$.
Hence one obtains:

\begin{Theorem}\label{Voronov}
Let $(\mathfrak{h},\mathfrak{a},\Pi_{\mathfrak{a}})$ be a V-algebra and $P$ 
a MC-element of $(\mathfrak{h},[-,-])$. 
Then
the family of higher derived brackets associated to $P$
\begin{align*}
(D_{P}^{n}\colon S^{n}(\mathfrak{a})\to\mathfrak{a}[1])_{n\in \mathbb{N}},
\end{align*}
equips $\mathfrak{a}$ with the structure of an $L_{\infty}[1]$-algebra (see Definition~\ref{L-infty}).
\end{Theorem}

\subsection{Homotopy Transfer}\label{s:homotopytransfer}
We describe a way to transfer $L_{\infty}$-algebras along retractions. Since we are not primiliarly interested in this transfer-procedure for its own sake but rather as a tool, we will not state the results of this subsection in the largest possible generality.

The two most serious restrictions are that we will assume 1. that the $L_{\infty}$-algebra we desire to transfer is a differential graded Lie algebra and 2. that the target of the transfer is the cohomology. We remark that a straightforward generalization of the procedure we are going to present works 1. for arbitrary $L_{\infty}$-algebras and 2. more general subcomplexes than the cohomology can be treated.

The situation is as follows:
Let $X$ be a graded vector space and $d$ a coboundary operator on $X$ (i.e. $d: X \to X[1]$ and $d \circ d=0$).
We denote the cohomology $H(X,d)$ by $H$. Assume that there are linear maps
\begin{itemize}
\item $h: X \to X[-1]$,
\item $pr: X \to H$ surjective and
\item $i: H \to X$ injective
\end{itemize}
such that the following conditions hold:
\begin{itemize}
\item $i$ and $pr$ are chain maps (i.e. $d\circ i = 0$ and $pr \circ d =0$),
\item $pr \circ i = id_{H}$,
\item $id_{X} - i \circ pr = d \circ h + h \circ d$ and
\item $h \circ h =0$, $h \circ i =0$ and $pr \circ h =0$ (sideconditions).
\end{itemize}
The tupel $(X,d,h,i,pr)$ is called \textsf{contraction data} and can be encoded in the following diagram:
\begin{align*}
\xymatrix{(H,0) \ar@<2pt>[r]^<<<<<{i}& \ar@<2pt>[l]^>>>>>{pr} (X,d)},h.
\end{align*}

\begin{Theorem}\label{thm:transfer}
Let $(X,d,h,i,pr)$ be a graded vector space equipped with contraction data and a finite compatible filtration, 
i.e. a collection of graded vector subspaces
\begin{align*}
X=\mathcal{F}_{0}X \supseteq \mathcal{F}_{1}X \supseteq \cdots \supseteq \mathcal{F}_{n}X \supseteq \mathcal{F}_{(n+1)}X \supseteq \cdots
\end{align*}
such that $\mathcal{F}_{N}X = \{0\}$ for $N$ large enough, satisfying
\begin{itemize}
\item $d(\mathcal{F}_{k}X) \subset \mathcal{F}_{k}X$ for all $k \ge 0$ and
\item $h(\mathcal{F}_{k}X) \subset \mathcal{F}_{k}X$ for all $k \ge 0$.
\end{itemize}
Futhermore suppose $X$ is equipped with the structure of a differential graded Lie algebra $(X,D,[-,-])$ such that
\begin{itemize}
\item $(D-d)(\mathcal{F}_{k}X) \subset \mathcal{F}_{(k+1)}X$.
\end{itemize}
Then the cohomology $H$ of $(X,d)$ is naturally equipped with the structure
of a (flat) $L_{\infty}$-algebra and there is a well-defined $L_{\infty}$-morphism
$\hat{i}:H \leadsto X$.
\end{Theorem}

In all the cases where we apply Theorem \ref{thm:transfer} it will be straightforward to check that the $L_{\infty}$-morphism described in Lemma \ref{l:inducedembedding} is in fact an $L_{\infty}$ quasi-isomorphsim.

The conceptual proof of Theorem \ref{thm:transfer} is straightforward and can be found in \cite{GugenheimLambe} for instance. One makes use of the interpretation of the $L_{\infty}$-algebra structure on $X$ as a codifferential $Q$ on $S(X[1])$ and uses transfer formulae for $Q$ to obtain a codifferential
$\mathcal{Q}$ on $S(H[1])$, i.e. a $L_{\infty}$-algebra structure on $H$. Moreover there are well-known formuals for $\hat{i}$.  

Although Theorem \ref{thm:transfer} establishes the existence of a transfer-procedure along contraction data, we need a more concrete description of the induced $L_{\infty}$-algebra and of the $L_{\infty}$ quasi-isomorphism between
$H$ and $X$. Such a description was first given in the setting of $A_{\infty}$-algebras: in \cite{Merkulov} inductive formulae where presented for the structure maps of the induced structure and in \cite{KontsevichSoibelman} an interpretation in terms of Feynman diagrams was provided.
Similar descriptions are known to hold for the transfer of $L_{\infty}$-algebras as well, although we need
a slight generalization of the setting presented in \cite{Merkulov} and \cite{KontsevichSoibelman} since we allow the coboundary operator $D$ to deviate from $d$.

We present the description of the transfer along contraction data using diagrams. Since we do not claim any
originality on the material which is well-known to the experts, we only state the results. The interested reader can find the proofs in the Appendix.

An \textsf{oriented decorated tree} $T$ is a finite
connected graph without loops of any kind that only consists of directed edges and trivalent interior vertices with two incoming edges
and one outcoming one. There are two kind of exterior vertices: ones with an outgoing edge -- we call
these leaves -- and exactly one with an incoming edge that we call the root. The orientation
is given by an association of two numbers to any pair of edges with the same vertex as their target
that tells us which of the two edges is the
``right'' and which is the ``left'' one. The decoration is an assignment of a non-negative integer to each edge.

\vspace{0.2cm}
\begin{minipage}[c]{\textwidth}
\centering \includegraphics[width=5cm]{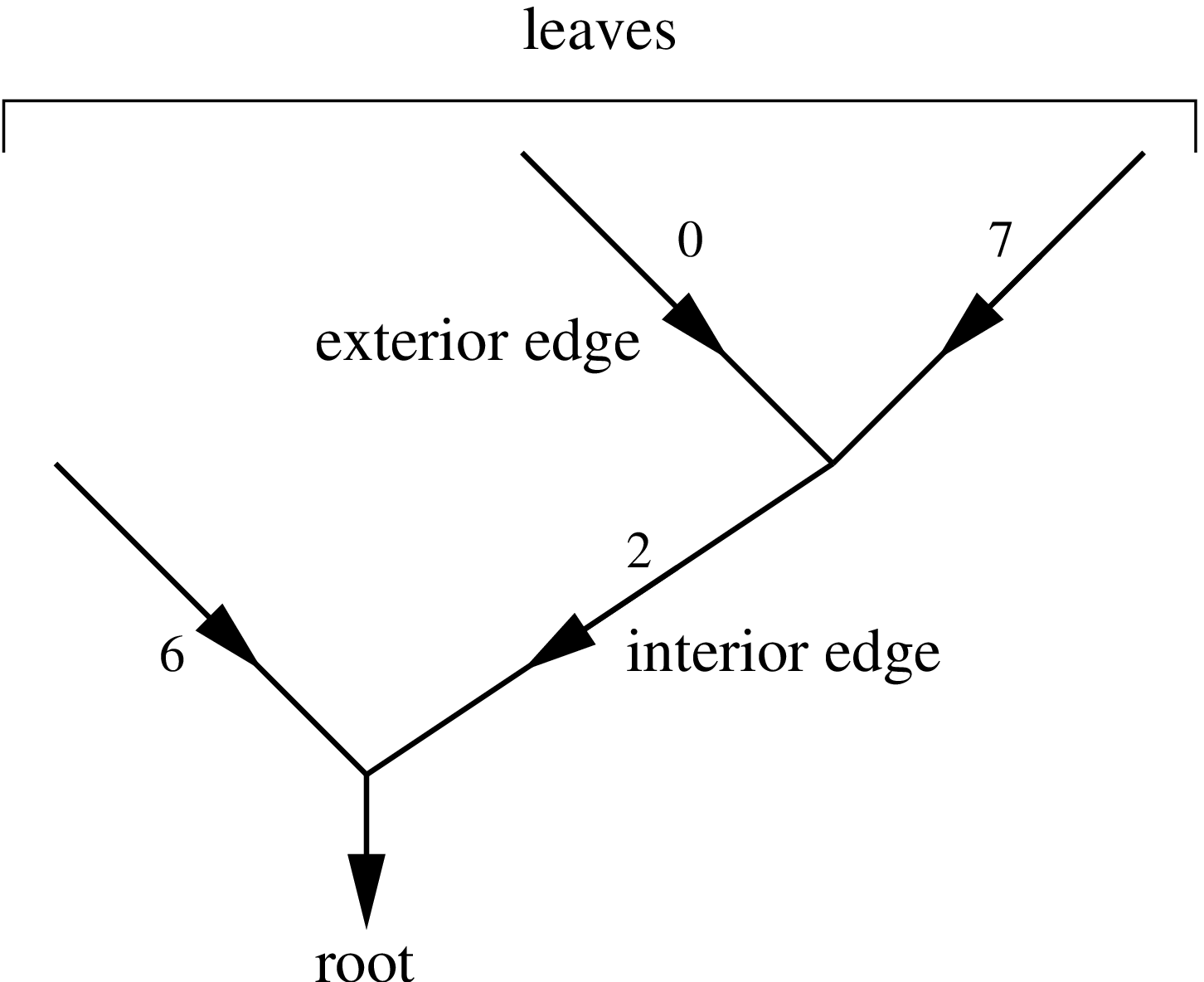}
\vspace{0.2cm}
\end{minipage}

The edge of the diagram with consists of only one leaf which is connected to the root must be decorated by a positive integer.
Clearly we have a decomposition
\begin{align*}
\mathbb{T} = \bigsqcup_{n \ge 1} \mathbb{T}(n)
\end{align*}
where $\mathbb{T}(n)$ denotes the set of trees with exactly $n$ leaves. 
We will denote the set of unoriented decorated trees by $[\mathbb{T}]$. 
There
is a natural projection
\begin{align*}
[\cdot]: \mathbb{T} \to [\mathbb{T}]
\end{align*}
that respects the decomposition of $\mathbb{T}$ and that of $[\mathbb{T}]$:
\begin{align*}
[\mathbb{T}]= \bigsqcup_{n \ge 1} [\mathbb{T}](n) = \bigsqcup_{n \ge 1}[\mathbb{T}(n)].
\end{align*}
We define $|Aut(T)|$ for $T$ an oriented decorated graph to be the cardinality of the group of automorphisms of the  underlying unoriented
decorated graph.

%
%An \textsf{oriented decorated tree} $T$ is a non-empty connected finite graph without cycles and double edges containing vertices that are unary or trivalent only. Unary vertices are called \textsf{exterior vertices} while the trivalent ones are called \textsf{interior vertices}. Edges connecting interior vertices are called \textsf{interior edges} and edges connecting an exterior vertex with an interior one are called \textsf{exterior edges}. We single out one of the exterior vertices to which we refer to as the \textsf{root}, the remaining exterior vertices are the \textsf{leaves}. It is possible to orient the edges such that all edges are oriented towards the root, in particular all the exterior edges connected to leaves point away from the leaves. The \textsf{decoration} is an assignment of a non-negative integer to each edge.
%We denote the cardinality of the group of automorphisms of $T$ by $|Aut(T)|$.
%
%
%The set of such trees $\mathbb{T}$ carries a decomposition given by the number of leaves, i.e. $\mathbb{T}(k)$ is the set of oriented decorated trees with $k$ leaves. Moreover there is a well-defined notion of isomorphisms of such trees given by an isomorphism of trees that respects the orientation and the decoration. We denote the set of
%isomorphism classes of oriented decorated trees by $[\mathbb{T}]$. It decomposes into components given by the set of isomorphism classes of oriented decorated trees with $k$ leaves which we denote by $[\mathbb{T}](k)$.
%
Consider $X$ equipped with contraction data $(X,d,h,i,pr)$ and the structure of a differential graded Lie algebra $(X,D,[-,-])$ satisfying all conditions stated in Theorem \ref{thm:transfer}. Then one assignes to any tree $T\in \mathbb{T}(k)$ a map
\begin{align*}
m_{T}: (H[1])^{\otimes k} \to H[2]
\end{align*}
as follows: Using the d\'ecalage-isomorphism we equip $X[1]$ with the structure of an $L_{\infty}[1]$-algebra with structure maps $\mu^{1}$ and $\mu^{2}$ (corresponding to $D$ and $[-,-]$ respectively). 
We write $\mu^{1}= d + \mu_{\Delta}^{1}$. Next we put a $\mu^{2}$ at each interior vertex of $T \in \mathbb{T}(k)$ and a number of $\mu_{\Delta}^{1}$s at every edge -- the number of $\mu^{1}_{\Delta}$s is given by the number
decorating the edge under consideration. Between any two consecutive operations one puts $-h$.
Finally one places $i$ at the leaves and $pr$ at the root. The orientation of the tree induces a numbering of the leaves of $T$ and applying all these maps in the order given by the orientation of the tree 
yields the map $m_{T}$.

It is easy to check that the ``symmetrization'' 
\begin{align*}
\sum_{\sigma \in \Sigma_k}\frac{1}{|Aut(T)|}\sigma^{*}(m_T)
\end{align*}
does not depend on the specific choice of the orientation of $T$.

%It is easy to check that $m_{T}$ evaluated on elements of $S(H[1])$ does not change if $T$ is replaced by an isomorphic tree. 
Hence we get a map
\begin{align*}
\hat{m}:[\mathbb{T}] \to Hom(S(H[1]),H[2])
\end{align*}
and consquently
\begin{align*}
\nu^{k}:=\sum_{[T]\in [\mathbb{T}](k)}\hat{m}([T])
\end{align*}
is well-defined.

\begin{Lemma}\label{l:inducedstructure}
The sequence of maps $(\nu^{k}:S^{k}(H[1]) \to H[2])_{k \ge 1}$ defines the structure of an
$L_{\infty}[1]$-algebra on $H[1]$.
\end{Lemma}

See the Appendix for a proof of this statement.

The $L_{\infty}[1]$-morphism $\hat{i}: H[1] \leadsto X[1]$ is also given in terms of oriented decorated trees. This time we associate the following map
\begin{align*}
n_{T}: H[1]^{\otimes k} \to X[1]
\end{align*}
to a tree $T$ in $\mathbb{T}(k)$: again place $\mu^{2}$ at all interior vertices, $l$ copies of $\mu^{1}_{\Delta}$ at edges decorated by $l$ and between two consecutive operations of this kind place $-h$. As before put $i$ at the leaves. The only difference is that we put a $-h$ at the root instead of $pr$.
Again it is straightforward to check that the ``symmetrization'' 
\begin{align*}
\sum_{\sigma \in \Sigma_k}\frac{1}{|Aut(T)|}\sigma^{*}(n_T)
\end{align*}
does not depend on the choice of orientation of $T$ and we obtain a map
%Again it is straightforward to check that $n_{T}$ restricted to $S^{k}(H[1])$ depends only on the isomorphism class of $T$ in $\mathbb{T}$ and consequently we obtain a map
\begin{align*}
\hat{n}: [\mathbb{T}] \to Hom(S(H[1]),X[1]).
\end{align*}
One defines a family of maps
\begin{align*}
\lambda^{k}:=\sum_{[T]\in [\mathbb{T}](k)}\hat{n}([T])
\end{align*}
that satisfies

\begin{Lemma}\label{l:inducedembedding}
The sequence of maps $(\lambda^{k}:S^{k}(H[1]) \to X[1])_{k \ge 1}$ defines an $L_{\infty}[1]$-morphism between $(H[1],\nu^{1},\nu^{2},\dots)$ and $(X[1],\mu^{1},\mu^{2})$.
\end{Lemma}

The interested reader can find a proof of this statement in the Appendix.

\subsection{Smooth graded Manifolds}\label{s:sgmfs}

\begin{Definition}\label{def:smgfs}
Let $M$ be a smooth finite dimensional manifold.

A \textsf{(bounded) graded vector bundle} over $M$ is a collection $E_{\bullet}=\{E_{i}\}_{i\in \mathbb{Z}}$ of finite rank vector bundles over $M$ such that $E_{k}=\{0\}$ for $k$ smaller than some $k_{min}$ or larger than some $k_{max}$. Since we only consider bounded graded vector bundles we will drop the adjective bounded from now on. 

The \textsf{algebra of smooth functions} on a graded vector bundle $E_{\bullet}$ is the graded commutative associative algebra
\begin{align*}
\mathcal{C}^{\infty}(E_{\bullet}):=\Gamma(\otimes_{k\in \mathbb{Z}}\mathcal{T}^{-k\bullet}(E^{*}_{k}))
\end{align*}
where $\mathcal{T}^{-k\bullet}(E^{*}_{k})$ is $\bigwedge^{-k\bullet}(E^{*}_{k})$ for $k$ odd and
$S^{-k\bullet}(E^{*}_k)$ for $k$ even.
The symbol $\otimes$ refers to the completed tensor product over $\mathcal{C}^{\infty}(M)$. Moreover the algebraic structure on the tensor product of two graded associative algebras is declared to be $(a\otimes x) \cdot (b \otimes y):=(-1)^{|x||b|} (a\cdot b)\otimes (x\cdot y)$.

A \textsf{morphism} between two graded vector bundles $E_{\bullet}$ and $F_{\bullet}$ is a morphism
of unital graded commutative associative algebras between $\mathcal{C}^{\infty}(F_{\bullet})$ and $\mathcal{C}^{\infty}(E_{\bullet})$. 
\end{Definition}

We define the shift operator $[n]$ on smooth graded vector bundles by $E_{\bullet}[n]:=\{E_{i+n}\}_{i\in \mathbb{Z}}$. 

\begin{Definition}
A \textsf{smooth graded manifold} $\mathcal{M}$ is a unital graded commutative asssociative algebra
$\mathcal{A}_{M}$ that is isomorphic to $\mathcal{C}^{\infty}(E_{\bullet})$ for some graded vector bundle
$E_{\bullet}$. We define $\mathcal{C}^{\infty}(\mathcal{M}):=\mathcal{A}_{M}$.

A \textsf{morphism} between two smooth graded manifolds $\mathcal{M}$ and $\mathcal{N}$ is a morphism
of unital graded commutative algebras from $\mathcal{C}^{\infty}(\mathcal{N})$ to $\mathcal{C}^{\infty}(\mathcal{M})$.
\end{Definition}

We remark that a specific isomorphism between $\mathcal{C}^{\infty}(\mathcal{M})$ and $\mathcal{C}^{\infty}(E_{\bullet})$
is not part of the data that define the smooth graded manifold $\mathcal{M}$.

Let $\mathcal{M}$ be a smooth graded manifold and let $\mathcal{X}(\mathcal{M})$ be
the vector space of graded derivations of $\mathcal{C}^{\infty}(\mathcal{M})$, i.e. $\phi \in \mathcal{X}_{k}(\mathcal{M})$ iff
\begin{align*}
\phi:  \mathcal{C}^{\infty}(\mathcal{M}) \to \mathcal{C}^{\infty}(\mathcal{M})[k]
\end{align*}
satisfies $\phi(a\cdot b)=\phi(a) \cdot b + (-1)^{k|a|} a \cdot \phi(b)$ for homogeneous $a$ and $b$ in $\mathcal{C}^{\infty}(\mathcal{M})$.

\begin{Definition}
Let $\mathcal{M}$ be a smooth graded manifold. The algebra of \textsf{multivector fields} on $\mathcal{M}$
is the graded commutative associative algebra
\begin{align*}
\mathcal{V}(\mathcal{M}):= S_{\mathcal{C}^{\infty}(\mathcal{M})}(\mathcal{X}(\mathcal{M})[-1]),
\end{align*}
i.e. the graded symmetric algebra generated by $\mathcal{X}(\mathcal{M})[-1]$ as a graded module
over $\mathcal{C}^{\infty}(\mathcal{M})$.
\end{Definition}

Let $\phi,\psi \in \mathcal{X}(\mathcal{M})$ be homogeneous elements of degree $|\phi|$ and $|\psi|$ respectively.
Then
\begin{align*}
[\phi,\psi]:=\phi \circ \psi - (-1)^{|\phi||\psi|}\psi \circ \phi
\end{align*}
defines the structure of a graded Lie algebra on $\mathcal{X}(\mathcal{M})$.
This bracket can be extented to a graded Lie algebra bracket $[-,-]_{SN}$ (SN stands for Schouten-Nijenhuis) on $\mathcal{V}(\mathcal{M})[1]$ by imposing the condition that $[-,-]_{SN}$ is a graded biderivation of $\mathcal{V}(\mathcal{M})$.

Assume that the smooth graded manifold $\mathcal{M}$ is represented by the graded vector bundle $E_{\bullet}\to M$.
Using connections on the components of $E_{\bullet}$ one sees that there is an isomorphism
between $\mathcal{V}(\mathcal{M})$ and $\mathcal{C}^{\infty}(T^{*}[1]M\oplus E_{\bullet} \oplus E^{*}_{\bullet}[1])$
where $E^{*}_{\bullet}$ refers to the graded vector bundle $\{E^{*}_{-i}\}_{i\in\mathbb{Z}}$. Hence:

\begin{Lemma}\label{l:mvfs}
Let $\mathcal{M}$ be a smooth graded manifold. Then the graded commutative algebra of multivector fields
$\mathcal{V}(\mathcal{M})$ on $\mathcal{M}$ defines a smooth graded manifold.
\end{Lemma}

Let $Z\in \mathcal{V}(\mathcal{M})$ be a bivector field on $\mathcal{M}$ of total degree $0$ (i.e. an element of $S^{2}_{\mathcal{C}^{\infty}(\mathcal{M})}(\mathcal{X}(\mathcal{M})[-1])$). The algebra $\mathcal{C}^{\infty}(\mathcal{M})[1]$ is an abelian Lie subalgebra of $(\mathcal{V}(\mathcal{M})[1],[-,-]_{SN})$ hence we can construct the derived brackets $(D^{n}_{Z})$ associated to $Z$, see subsection \ref{s:derivedbrackets}. The only possible non-vanishing term is $D^{2}_{Z}$. Using the d\'ecalage-isomorphism we obtain am map $\bigwedge^{2}(\mathcal{C}^{\infty}(\mathcal{M})) \to \mathcal{C}^{\infty}(\mathcal{M})$ which we denote by $[-,-]_{Z}$. According to
Theorem \ref{Voronov} in subsection \ref{s:derivedbrackets}, $[-,-]_{Z}$ equips $\mathcal{C}^{\infty}(\mathcal{M})$ with the structure of a graded Lie algebra if $Z$ satisfies $[Z,Z]_{SN}=0$. It can be checked that in this case $(\mathcal{C}^{\infty}(\mathcal{M}),[-,-]_{Z})$ is a graded Poisson algebra.

\subsection{Poisson Geometry}\label{s:Poisson}
Let $M$ be a smooth finite dimensional manifold. In subsection \ref{s:sgmfs} the Schouten-Nijenhuis bracket $[-,-]_{SN}$ was introduced: it equips $\mathcal{V}(M)[1]$ with the structure of a graded Lie algebra. A
Poisson bivector field $\Pi$ on $M$ is a MC-element of
$(\mathcal{V}(M)[1],[-,-]_{SN})$, i.e. $\Pi$ is a bivector field satifying $[\Pi,\Pi]_{SN}=0$.

Associated to any
Poisson bivector field $\Pi$ on $M$ there is a vector bundle morphism
$\Pi^{\#}: T^{*}M \to TM$ given by contraction. Denote the natural pairing
between $TM$ and $T^{*}M$ by $<\negthickspace-,-\negthickspace>$. The bracket on  $\mathcal{C}^{\infty}(M)$
defined by $[f,g]_{\Pi}:=<\negthickspace\Pi^{\#}(df),dg\negthickspace>$ is $\mathbb{R}$-bilinear, skew-symmetric,
satisfies the Jacobi-identity and is a biderivation for the multiplication on $\mathcal{C}^{\infty}(M)$. Hence $(\mathcal{C}^{\infty}(M),[-,-]_{\Pi})$
is a Poisson algebra.

Every Poisson manifold comes along with a singular foliation $\mathcal{F}_{\Pi}$, given
by 
\begin{align*}
\Pi^{\#}(T^{*}M) \hookrightarrow TM. 
\end{align*}
Locally this foliation is spanned
by elements of the form $\Pi^{\#}(df)$ for $f\in \mathcal{C}^{\infty}(M)$. The identity
\begin{align*}
[\Pi^{\#}(df),\Pi^{\#}dg]_{SN}=-\Pi^{\#}(d[f,g]_{\Pi})
\end{align*}
is satisfied which implies
that $\mathcal{F}_{\Pi}$ is involutive. By a generalization of the classical theorem of Frobenius due to Stefan and Sussman (see \cite{Stefan}, \cite{Sussman}) the
integrability of $\mathcal{F}_{\Pi}$ follows. The integrating leaves all carry a natural symplectic
structures induced from $\Pi$.

There is another interesting structure associated to every Poisson manifold $(M,\Pi)$.
Consider the binary operation on $\Gamma(T^{*}M)=\Omega^{1}(M)$ given by
\begin{align*}
[\alpha,\beta]_{K}:=\mathcal{L}_{\Pi^{\#}(\alpha)}(\beta) - \mathcal{L}_{\Pi^{\#}(\beta)}(\alpha)+d\Pi(\alpha,\beta)
\end{align*}
called the Koszul bracket. One can check that it is a Lie bracket on $\Omega^{1}(M)$ and that
the vector bundle morphism 
$\Pi^{\#}: T^{*}M \to TM$ induces an morphism of Lie algebras $(\Omega^{1}(M),[-,-]_K)\to (\mathcal{X}(M),[-,-]_{SN})$. Moreover the so called Leibniz identity holds: 
\begin{align*}
([\alpha,f\beta]_{K})=f[\alpha,\beta]_{K}+\Pi^{\#}(\alpha)(f)\cdot \beta 
\end{align*}
for all $\alpha,\beta \in \Omega^{1}(M)$
and $f\in \mathcal{C}^{\infty}(M)$. The triple
$(T^{*}M,[-,-]_{K},\Pi^{\#})$ is an example of a Lie algebroid over $M$. Associated to any Lie algebroid
is a cocomplex, called the Lie algebroid cocomplex. In fact this cocomplex encodes exactly the same information as the original Lie
algebroid data. In the case of the Lie algebroid $(T^{*}M,[-,-]_{K},\Pi^{\#})$ the Lie algebroid cocomplex
is $(\mathcal{V}(M),[\Pi,-]_{SN})$.

Consider a submanifold $S$ of $M$. The annihilator $N^{*}S$ of $TS$ is a natural subbundle of
$T^{*}M$. This
subbundle fits into a short exact sequence of vector bundles:
$$
\xymatrix{
0 \ar[r] & N^{*}S \ar[r] & T^{*}_{S}M \ar[r] & T^{*}S \ar[r] & 0}.
$$

\begin{Definition}\label{coisotropic}
A submanifold $S$ of a smooth finite dimensional Poisson manifold $(M,\Pi)$
is called \textsf{coisotropic} if the restriction of $\Pi^{\#}$ to $N^{*}S$
has image in $TS$.
\end{Definition}

Consequently any coisotropic submanifold $S$ is equipped with a natural singular foliation
$\mathcal{F}_{S}:=\Pi^{\#}(N^{*}S)$ which is involutive.
Involutivity of $\mathcal{F}_{S}$ follows from another equivalent characterization
of coisotropic submanifolds: define the vanishing ideal of $S$ by 
\begin{align*}
\mathcal{I}_{S}:=\{f\in \mathcal{C}^{\infty}(M): f|_{S}=0\}. 
\end{align*}
A submanifold $S$ is coisotropic if and only if $\mathcal{I}_{C}$ is a Lie subalgebra
of $(\mathcal{C}^{\infty}(M),[-,-]_{\Pi})$. Observe that 
$\Pi^{\#}(N^{*}S)$ is locally
spanned by $\Pi^{\#}(df)$ for $f \in \mathcal{I}_{S}$. For $f,g \in \mathcal{I}_{S}$ one
has
$[\Pi^{\#}(df),\Pi^{\#}(dg)]_{SN}=-\Pi^{\#}(d[f,g]_{\Pi})$. Since
$[f,g]_{\Pi} \in \mathcal{I}_{S}$ the foliation $\mathcal{F}_{S}$ is involutive. 
We denote the corresponding leaf space by $\underline{S}:=S/\negthickspace\sim_{\mathcal{F}_S}$. This space is usually very ill-behaved
(non-smooth, non-Hausdorff, etc.). In particular there might not be a meaningful
way to define $\mathcal{C}^{\infty}(\underline{S})$ using the topological
space $\underline{S}$. Instead one can define $\mathcal{C}^{\infty}(\underline{S})$
as the space of functions on $S$ which are invariant under $\mathcal{F}_{S}$, i.e.
\begin{align*}
\mathcal{C}^{\infty}(\underline{S}):=\{f\in\mathcal{C}^{\infty}(S): X(f)=0
\mbox{ for all } X \in \Gamma(\mathcal{F}_S)\}.
\end{align*}
This is a subalgebra of
$\mathcal{C}^ {\infty}(S)$.

Fix an embedding $\phi: NS \hookrightarrow M$ of the normal bundle of $S$ into $M$. Via the identification of $NS$ with an open neighbourhood
of $S$ in $M$ the vector bundle $NS$ inherits a Poisson bivector field $\Pi_{\phi}$. Hence we can assume without loss of generality that $M$ is the total space of a vector bundle $E \to S$. We will do so in the rest of the paper.
Observe that under the above assumptions there is a natural isomorphism $E \cong NS$.

With help of this assumption one sees that $\mathcal{C}^ {\infty}(\underline{S})$
comes quipped with a Poisson bracket $[-,-]_{\underline{S}}$ inherited from $(E,\Pi)$:
the algebra $\mathcal{C}^{\infty}(S)$ is the quotient of $\mathcal{C}^{\infty}(E)$ by
$\mathcal{I}_{S}$. There is an Lie algebra action of $(\mathcal{I}_{S},[-,-]_{\Pi})$
on this quotient. The algebra $\mathcal{C}^{\infty}(\underline{S})$ is given
by the invariants under this action, i.e.
\begin{align*}
\mathcal{C}^{\infty}(\underline{S})\cong (\mathcal{C}^{\infty}(E)/\mathcal{I}_{S})^{\mathcal{I}_{S}}. 
\end{align*}
This algebra is isomorphic to the quotient of
\begin{align*}
\mathcal{N}(\mathcal{I}_{S}):=\{f \in \mathcal{C}^{\infty}(E): [f,\mathcal{I}_{S}]_ {\Pi}\subset \mathcal{I}_{S}\}
\end{align*} by $\mathcal{I}_{S}$. It is straightforward to check that the
Poisson bracket on $\mathcal{C}^{\infty}(E)$ descends to this quotient.

The Lie algebroid structure $(T^{*}M,[-,-]_{K},\Pi^{\#})$ also restricts to co\-iso\-tro\-pic submanifolds:
the bundle map $\Pi^{\#}: T^{*}E \to TE$ restricts to a bundle map $E^{*} \to TS$ by definition 
and the Koszul bracket can also be restricted to $\Gamma(E^{*})$. The triple $(E^{*},[-,-]_{K},\Pi^{\#}|_{E^{*}})$ satisfies the same identities as $(T^{*}E,[-,-]_{K},\Pi^{\#})$ and hence is a Lie algebroid, see \cite{Weinstein} for details. The
easiest way to describe this Lie algebroid over $S$ is via its associated Lie algebroid cocomplex. Define
a projection $pr: \mathcal{V}(M) \to \Gamma(\bigwedge E)$ as the unique algebra
morphism extending the restriction $\mathcal{C}^{\infty}(E) \to \mathcal{C}^{\infty}(S)$
and $\mathcal{X}(E)=\Gamma(TE) \to \Gamma(T_{S}E) \to \Gamma(E)$.
The
graded algebra $\Gamma(\bigwedge E)$ is equipped with the differential given by
\begin{align*}
\partial_{S}(X):= pr([\Pi,\tilde{X}]_{SN}|_{S})
\end{align*}
where $\tilde{X}$ is any extension of $X \in \Gamma(\bigwedge E)$ to a multivector field on $E$. The cohomology of the cocomplex $(\Gamma(\bigwedge E),\partial_{S})$ is called the
\textsf{Lie algebroid cohomology} of $S$. It is well-known that

\begin{Lemma}
Let $S$ be a coisotropic submanifold of a smooth finite dimensionam Poisson manifold $(M,\Pi)$.
The algebra $\mathcal{C}^{\infty}(\underline{S})$ is isomorphic to the degree zero Lie algebroid cohomology $H^{0}(\bigwedge(NS),\partial_{S})$.
\end{Lemma}

Moreover it is possible to show that the Lie algebroid differential $\partial_{S}$ is independent of the embedding $NS \hookrightarrow M$ as is the Poisson bracket on $\mathcal{C}^{\infty}(\underline{S})$.

%% file: BFV.tex
Consider a finite rank vector bundle $E \to S$ that is equipped with a Poisson bivector field, i.e. $\Pi \in \mathcal{V}^2(E)$ satisfying
$[\Pi,\Pi]_{SN}=0$. Let $S$ be a coisotropic submanifold of $E$.

The aim of this section is to describe the construction of a homological
resolution of the Poisson algebra
$(\mathcal{C}^{\infty}(\underline{S}),[-,-]_{\underline{S}})$ (introduced in subsection \ref{s:Poisson}) in terms of a differential graded Poisson algebra 
\begin{align*}(BFV(E,\Pi),D_{BFV},[-,-]_{BFV}).
\end{align*}
$BFV(E,\Pi)$ can be described
as the space of smooth functions on some smooth graded manifold. The degree zero component
of the cohomology $H(BFV(E,\Pi),D_{BFV})$ is isomorphic to
$\mathcal{C}^{\infty}(\underline{S})$ and the induced bracket coincides with
$[-,-]_{\underline{S}}$.

The basic ideas of the construction of $(BFV(E,\Pi),D_{BFV},[-,-]_{BFV})$ were invented by Batalin, Fradkin and Vilkovisky
(\cite{BatalinFradkin}, \cite{BatalinVilkovisky}) with applications
to physics in mind. Later it was reinterprated by Stasheff in terms of homological algebra (\cite{Stasheff}). The proper globalization to the smooth setting was
presented by Bordemann and Herbig (\cite{Bordemann}, \cite{Herbig}). We essentially
follow \cite{Stasheff},\cite{Bordemann} and \cite{Herbig} in this exposition.
The only deviation will be a new conceptual approach to the Rothstein-bracket (\cite{Rothstein}) and its extension to the Poisson setting (\cite{Herbig}) in terms of higher homotopy structures given in section \ref{s:lifting}.

The construction of the BFV-complex relies on the following input data: $1.$ a choice of embedding of the normal bundle of $S$
as a tubular neighourhood (in order to obtain an appropriate vector bundle $E\to S$, see subsection \ref{s:Poisson}), $2.$ a connection on $E \to S$ and $3.$ a distinguished element $\Omega \in BFV(E,\Pi)$ satisfying
$[\Omega,\Omega]_{BFV}=0$. The dependence of the BFV-complex on these data will be clarified elsewhere (\cite{Schaetz}).

\subsection{The Ghost/Ghost-Momentum Bundle}\label{s:GGM_bundle}
Let $E\to S$ be a finite rank vector bundle over a smooth finite dimensional manifold.
Using the projection map of the vector bundle $E \to S$ we can pull back
the graded vector bundle $E^{*}[1]\oplus E[-1] \to S$ to a graded
vector bundle over $E$ which we denote by $\mathcal{E}^{*}[1]\oplus \mathcal{E}[-1] \to E$. The situation is summarized by the following Cartesian square:
$$
\xymatrix{
\mathcal{E}^{*}[1]\oplus\mathcal{E}[-1] \ar[r] \ar[d]_{P} & E^{*}[1] \oplus E[-1] \ar[d] &\\
E \ar[r] & S.}
$$

We define $BFV(E,\Pi)$ to be the space of smooth functions on the graded manifold
which is represented by the graded vector bundle $\mathcal{E}^{*}[1]\oplus \mathcal{E}[-1]$ over $E$.
In terms of sections one has $BFV(E,\Pi)=\Gamma(\bigwedge(\mathcal{E})\otimes\bigwedge(\mathcal{E^{*}}))$.
This algebra carries a bigrading given by 
\begin{align*}
BFV^{(p,q)}(E,\Pi):=\Gamma(\wedge^{p}(\mathcal{E})\otimes\wedge^{q}(\mathcal{E^{*}})).
\end{align*}
In physical terminology $p$ / $q$ is referred to as the
\textsf{ghost degree} / \textsf{ghost-momentum degree} respectively. One defines
\begin{align*}
BFV^{k}(E,\Pi):=\bigoplus_{p-q=k}BFV^{(p,q)}(E,\Pi)
\end{align*}
and calls $k$ the \textsf{total degree} (in physical terminology this is the ``ghost number''). There is yet another decomposition of $BFV(E,\Pi)$ that will be useful later: set 
$BFV_{r}(E,\Pi):=\Gamma(\bigwedge(\mathcal{E})\otimes \bigwedge^{r}(\mathcal{E^{*}}))$. Moreover we define $BFV_{\ge r}(E,\Pi)$ to be the ideal generated by
$BFV_{r}(E,\Pi)$.

The smooth graded manifold $\mathcal{E}^{*}[1]\oplus \mathcal{E}[-1]$ comes equipped with a Poisson
bivector field $G$ given by the natural fibre pairing between $\mathcal{E}$ and $\mathcal{E}^{*}$, i.e. it is defined to be the natural contraction on $\Gamma(\mathcal{E})\otimes \Gamma(\mathcal{E}^{*})$ and it extends uniquely
to a graded skew-symmetric biderivation of $BFV(E,\Pi)$.

\subsection{Lifting the Poisson Bivector Field}\label{s:lifting}
We want to equip $BFV(E,\Pi)$ with the structure of a graded Poisson algebra which
essentially combines the Poisson bivector field $\Pi$ on $E$ and the Poisson bivector field $G$ which encodes the natural
fibre pairing between $\mathcal{E}^{*}[1]$ and $\mathcal{E}[-1]$.

First we lift $\Pi$ from the base $E$ to the graded vector bundle
\begin{align*}
\mathcal{E}^{*}[1]\oplus \mathcal{E}[-1] \xrightarrow{P} E.
\end{align*}
For this purpose we choose
a connection $\nabla$ on the vector bundle $E \to S$. This yields a connection on $E^{*}[1]\oplus E[-1]$. Pulling back this connection along $E\to S$ gives a connection on $\mathcal{E}^{*}[1]\oplus \mathcal{E}[-1] \to E$
that is metric with respect to the natural fibre pairing. Fix
such a connection on $\mathcal{E}^{*}[1]\oplus \mathcal{E}[-1] \to S$ and consider the horizontal lift with respect to that connection,
i.e. we obtain a map
$\iota_{\nabla}: \mathcal{X}(E) \hookrightarrow \mathcal{X}(\mathcal{E}^{*}[1]\oplus \mathcal{E}[-1])$. Setting $\iota_{\nabla}(f):=f\circ P$ for $f\in \mathcal{C}^{\infty}(E)$ we can uniquely extend $\iota_{\nabla}$
to a morphism of algebras 
\begin{align*}
\iota_{\nabla}:  \mathcal{V}(E) \hookrightarrow \mathcal{V}(\mathcal{E}^{*}[1]\oplus \mathcal{E}[-1]).
\end{align*}

Since $\iota_{\nabla}[1]$ fails in general to be a morphism of graded Lie algebras,
the horizontal lift $\iota_{\nabla}(\Pi)$ of the Poisson bivector field $\Pi$ does not to satisfy the MC-equation in
$(\mathcal{V}(\mathcal{E}^{*}[1]\oplus \mathcal{E}[-1])[1],[-,-]_{SN})$. The same is true for
the sum $G+\iota_{\nabla}(\Pi)$, hence this bivector field does not define
the structure of a graded Poisson algebra on $BFV(E,\Pi)$. We will show that an
appropriate correction term $\triangle$ can be found such that
$G+\iota_{\nabla}(\Pi)+\triangle$ is a MC-element.
The existence of such a $\triangle$ is the straightforward consequence of the following
proposition:

\begin{Proposition}\label{prop:lift}
Let $\mathcal{E}$ be a finite rank vector bundle with connection $\nabla$ over a finite dimensional smooth manifold $E$.
Consider the smooth graded manifold $\mathcal{E}^ {*}[1]\oplus \mathcal{E}[-1] \to E$ and denote the Poisson
bivector field on it coming from the natural fibre pairing between $\mathcal{E}$ and $\mathcal{E}^{*}$ by $G$.

Then there is an $L_{\infty}$ quasi-isomorphism $\mathcal{L}_{\nabla}$ between the graded Lie algebra $(\mathcal{V}(E)[1],[-,-]_{SN})$ and the differential graded Lie algebra \newline $(\mathcal{V}(\mathcal{E}^{*}[1]\oplus \mathcal{E}[-1])[1],[G,-]_{SN},[-,-]_{SN})$.
\end{Proposition}

Observe that it is not assumed in the Proposition that $E$ is a vector bundle or that $\mathcal{E} \to E$ is a pull back bundle.

\begin{proof}
Consider the induced connection $\nabla$ on
$\mathcal{E}^{*}[1]\oplus \mathcal{E}[-1] \to E$ (by a slight abuse of notation we denote this connection again by $\nabla$). It is metric with respect to the natural fibre
pairing. The algebra morphism $\iota_{\nabla}: \mathcal{V}(E) \hookrightarrow
\mathcal{V}(\mathcal{E}^{*}[1]\oplus \mathcal{E}[-1])$ (given by the horizontal lift) is a section
of the natural projection 
\begin{align*}
Pr: \mathcal{V}(\mathcal{E}^{*}[1]\oplus \mathcal{E}[-1]) \to \mathcal{V}(E).
\end{align*}
Obviously $Pr \circ \iota_{\nabla} = id$ holds on $\mathcal{V}(E)$.

Consider the complexes $(\mathcal{V}(\mathcal{E}^{*}[1]\oplus \mathcal{E}[-1]), Q:=[G,-]_{SN})$ and
$(\mathcal{V}(E),0)$. It is easy to check that $Pr$ and $\iota_{\nabla}$
are chain maps. Here it is crucial that the induced connection on $\mathcal{E}^{*}[1]\oplus \mathcal{E}[-1]$
is metric with respect to the natural fibre pairing). We will construct a homotopy
\begin{align*}
H_{\nabla}:=\mathcal{V}(\mathcal{E}^{*}[1]\oplus \mathcal{E}[-1]) \to \mathcal{V}(\mathcal{E}^{*}[1]\oplus \mathcal{E}[-1])[-1]
\end{align*}
such that $Q\circ H_{\nabla} + H_{\nabla}\circ Q = 1 - \iota_{\nabla}\circ Pr$, i.e. $\iota_{\nabla}$
and $Pr$ are inverses up to homotopy and it follows that $Pr$
induces an isomorphism $H(\mathcal{V}(E^{*}[1]\oplus E[-1]),Q)\cong \mathcal{V}(M)$.

To construct an appropriate homotopy $H_{\nabla}$ we extend $\iota_{\nabla}$ to
an algebra isomorphism
\begin{align*}
\varphi_{\nabla}: \mathcal{A}:=\mathcal{C}^{\infty}(T^{*}[1]E\oplus \mathcal{E}^{*}[1]\oplus \mathcal{E}[-1] \oplus \mathcal{E}[0] \oplus \mathcal{E}^{*}[2]) \to \mathcal{V}(\mathcal{E}^{*}[1]\oplus \mathcal{E}[-1]),
\end{align*}
see Lemma \ref{l:mvfs} in subsection \ref{s:sgmfs}.
Via this identification we equip $\mathcal{A}$ with a Gerstenhaber bracket $[-,-]_{\nabla}$ and a differential $\tilde{Q}:=\varphi_{\nabla}^{-1}\circ Q \circ \varphi_{\nabla}$. Define $\tilde{H}$ to be the sum of the pullbacks by the maps $-id_{\mathcal{E}^{*}[1]}[1]: \mathcal{E}^{*}[1] \to \mathcal{E}^{*}[2]$
and $-id_{\mathcal{E}[-1]}[1]: \mathcal{E}[-1] \to \mathcal{E}[0]$ on $\mathcal{A}$. It is straightforward to check
that $\tilde{H}$ is a differential and that
$(\tilde{Q}\circ\tilde{H} + \tilde{H}\circ \tilde{Q})(X)$ is equal to the total polynomial degree of $X$ in all of the fibre components $\mathcal{E}^{*}[1], \mathcal{E}^{*}[2], \mathcal{E}[-1]$ and $\mathcal{E}[0]$. 
Normalising $\tilde{H}$ and using the idenfitication $\varphi_{\nabla}$ leads to a homotopy $H_{\nabla}$ on $\mathcal{V}(\mathcal{E}^{*}[1]\oplus \mathcal{E}[-1])$. It is straightforward to check that the side-conditions $H_{\nabla}\circ H_{\nabla}=0$, $H_{\nabla}\circ \iota_{\nabla}=0$ and $Pr\circ H_{\nabla}=0$ hold.

We summarize the situation in the following diagram:
\begin{align*}
\xymatrix{(\mathcal{V}(E),0) \ar@<2pt>[r]^<<<<<{\iota_{\nabla}}& \ar@<2pt>[l]^>>>>>{Pr} (\mathcal{V}(\mathcal{E}^{*}[1]\oplus \mathcal{E}[-1]),Q)},H_{\nabla}.
\end{align*}
According to subsection \ref{s:homotopytransfer} these data can be used to perform homological transfer of $L_{\infty}$-algebra structures along
the contraction $Pr$. Starting with the differential graded Lie algebra
$(\mathcal{V}(\mathcal{E}^{*}[1]\oplus \mathcal{E}[-1])[1], Q=[G,-]_{SN},[-,-]_{SN})$ one
constructs an $L_{\infty}$ quasi-isomorphic $L_{\infty}$-algebra structure
on $\mathcal{V}(M)[1]$ (with zero differential) together with an $L_{\infty}$ quasi-isomorphism $\mathcal{L}_{\nabla}$. The binary operation of this
structure will simply be  given by 
\begin{align*}
Pr([\iota_{\nabla}(-),\iota_{\nabla}(-)]_{SN})=[-,-]_{SN}.
\end{align*}
All potential higher operations can be checked to vanish as follows:
As described in \ref{s:homotopytransfer} one considers all trivalent oriented trees.
On the leaves (i.e. exterior vertices with edges oriented away from them) one places $\iota_{\nabla}$, on each
interior trivalent vertex one places $[-,-]_{SN}$, on the root (i.e.
the unique exterior vertex with edge oriented towards it) one places $Pr$ and on
interior edges (those not connected
to any leaf or to the root) one places $-H_{\nabla}$. Then
one composes these maps in the order given by the orientation of the tree.

To prove that no higher order operations occur we introduce a decomposition of $\mathcal{V}(\mathcal{E}^{*}[1]\oplus \mathcal{E}[-1])$.
By definition this is the space of multiderivations of the graded unital algebra $\mathcal{C}^{\infty}(\mathcal{E}^{*}[1]\oplus \mathcal{E}[-1])$.
The algebra of smooth functions is bigraded which induces a bigrading on its tensor algebra (just take the sum of the bidegrees of all tensor
components) which in turn induces a bigrading on the space of multivector fields, i.e. an element of bidegree
$(m,n)$ is one that maps a tensor product of function of total bidegree $(p,q)$ to a function of bidegree $(p+m,q+n)$. This bidegree
is obviously bounded from above. We denote the ideal generated by the component of bidegree $(M,N)$ with $M\ge m$ and $N \ge n$ by
$\mathcal{V}^{(m,n)}(\mathcal{E}^{*}[1]\oplus \mathcal{E}[-1])$.

Consider a tree as above and forget about $Pr$ at the root. One
can inductively show that the corresponding operation maps
tensor products of elements of $\mathcal{V}(E)$ to $\mathcal{V}^{(e-1,e-1)}(\mathcal{E}^{*}[1]\oplus \mathcal{E}[-1])$
where $e$ is the number of trivalent vertices of the tree. This relies on the following 

\begin{Lemma}
Denote by $R_{\nabla}$ the curvature of the connection $\nabla$ on $\mathcal{E} \to E$. We interpret $R_{\nabla}$ as an element of
$\Omega^2(E,End(\mathcal{E}))=\Omega^2(E,\mathcal{E}\otimes \mathcal{E}^{*})$. Then
\begin{align*}
R_{\nabla}(-,-)=H_{\nabla}([\iota_{\nabla}(-),\iota_{\nabla}(-)]_{SN})
\end{align*}
holds.
\end{Lemma}

\begin{proof}[Proof of the Lemma]
The left hand side of the claimed equality can be checked to be $\mathcal{C}^{\infty}(E)$-bilinear and multiplicative in both slots with respect to
the algebra structure on $\mathcal{V}(E)$. Hence it is determined by its values on a pair of vector fields and can be interprated
as a two-form on $E$ with values in a vector bundle. Consequently it is enough to prove the equality locally which is a
straightforward computation in coordinate charts.
\end{proof}

So all operations vanish identically after applying $Pr$ except for the
case of the tree with only one trivalent edge (which correspondes to the binary
operation $[-,-]_{SN}$).
\end{proof}

\begin{Corollary}\label{cor:lift}
Let $\mathcal{E } \to E$ be a finite rank vector bundle with connection $\nabla$ over a smooth finite dimensional Poisson manifold
$(E,\Pi)$. Consider the smooth graded manifold $\mathcal{E}^{*}[1]\oplus \mathcal{E}[-1] \to E$ and denote the Poisson
bivector field on it coming from the natural fibre pairing between $\mathcal{E}$ and $\mathcal{E}^{*}$ by $G$.

Then there is a Poisson bivector field $\hat{\Pi}$ on $\mathcal{E}^{*}[1]\oplus \mathcal{E}[-1]$ such that
\newline $\hat{\Pi}= G + \iota_{\nabla}(\Pi) + \triangle$ for $\triangle \in \mathcal{V}^{(1,1)}(\mathcal{E}^{*}[1]\oplus \mathcal{E}[-1])$. 
\end{Corollary}

Recall that $\mathcal{V}^{(1,1)}(\mathcal{E}^{*}[1]\oplus \mathcal{E}[-1])$ is the ideal of $\mathcal{V}(\mathcal{E}^{*}[1]\oplus \mathcal{E}[-1])$
generated by multiderivations which map a tensor product of functions of total bidegree $(p,q)$ to a function of bidegree $(P,Q)$ where $P>p$ and $Q>q$.

This Corollary was originally proven by Rothstein in \cite{Rothstein} for $(N,\Pi)$ symplectic with the help
of a concrete formula for $\hat{\Pi}$. Herbig
showed that Rothstein's formula holds also in the Poisson case (\cite{Herbig}).

\begin{proof}
The general theory of $L_{\infty}$-algebras implies that given two $L_{\infty}$ quasi-isomorphic
$L_{\infty}$-algebras and a formal MC-element of one of these $L_{\infty}$-algebras, one can construct
a formal MC-element of the other one. We apply this to the Poisson bivector field
$\Pi$ seen as a MC-element in $(\mathcal{V}(E)[1],[-,-]_{SN})$ which is
$L_{\infty}$ quasi-isomorphic to $(\mathcal{V}(\mathcal{E}^{*}[1]\oplus\mathcal{E}[-1]),[G,-]_{SN},[-,-]_{SN})$ according to
Proposition \ref{prop:lift}.

The unary operation from $\mathcal{V}(E)$ to $\mathcal{V}(\mathcal{E}^{*}[1]\oplus \mathcal{E}[-1])$ is given by $\iota_{\nabla}$. The higher structure maps of the $L_{\infty}$-morphism between $\mathcal{V}(E)$ and $\mathcal{V}(\mathcal{E}^{*}[1]\oplus \mathcal{E}[-1])$
are given in terms of trivalent oriented trees. One places $\iota_{\nabla}$ at leaves (i.e. exterior vertices with
edges oriented away from them), $[-,-]_{SN}$ at trivalent interior vertices and the homotopy $-H_{\nabla}$ at all interior edges
(all edges not connected to a leaf or root) and at the edge connected to the root (the unique exterior
vertex with the edge oriented towards it). There is an estimate similar to the one in the proof of Proposition
\ref{prop:lift}:  the operation corresponding to a tree with $e$ trivalent edges maps
elements of $\mathcal{V}(E)$ to $\mathcal{V}^{(e,e)}(\mathcal{E}^{*}[1]\oplus \mathcal{E}[-1])$.

This implies 1. that
we do not have to care about convergence since the filtration of $\mathcal{V}(\mathcal{E}^{*}[1]\oplus \mathcal{E}[-1])$
by the ideals $\mathcal{V}^{(k,l)}(\mathcal{E}^{*}[1]\oplus \mathcal{E}[-1])$ is bounded from above, so only finitely many
trees will contribute. And 2. by applying the $L_{\infty}$ quasi-isomorphism to $\Pi$ one obtains a Maurer-Cartan element
of $(\mathcal{V}(\mathcal{E}^{*}[1]\oplus \mathcal{E}[-1])[1],[G,-]_{SN},[-,-]_{SN})$ of the form $\iota_{\nabla}(\Pi)+\triangle$
with $\triangle \in \mathcal{V}^{(1,1)}(\mathcal{E}^{*}[1]\oplus \mathcal{E}[-1])$. This is equivalent to the statement
that $G+\iota_{\nabla}(\Pi)+\triangle$ is a Maurer-Cartan element of $\mathcal{V}(\mathcal{E}^{*}[1]\oplus \mathcal{E}[-1])[1],[-,-]_{SN})$ of the desired form.
\end{proof}

By definition such an element yields the structure of a graded Poisson algebra on
$\mathcal{C}^{\infty}(\mathcal{E}^{*}[1]\oplus \mathcal{E}[-1])=:BFV(E,\Pi)$:

\begin{Corollary}\label{cor:BFV-bracket}
Let $E\to S$ be a finite rank vector bundle over a smooth finite dimensional manifold $S$. Assume $(E,\Pi)$ is a Poisson manifold.
Consider the associated ghost/ghost-momentum bundle $\mathcal{E}^{*}[1]\oplus \mathcal{E}[-1] \xrightarrow{P} E$
with the embedding $j: E \hookrightarrow \mathcal{E}^{*}[1]\oplus \mathcal{E}[-1]$ as the zero section.
The natural fibre paring between $\mathcal{E}^{*}$ and $\mathcal{E}$ gives rise to
a Poisson bivector field $G$.

Then there is a well-defined graded Poisson bracket $[-,-]_{BFV}$ on $BFV(E,\Pi)$ such that:
\begin{enumerate}
\item $[-,-]_{\Pi}=j^{*}([P^{*}(-),P^{*}(-)]_{BFV})$ and
\item denoting the projection $BFV^{0}(E,\Pi)\to BFV^{(0,0)}(E,\Pi)$ by $proj$ the composition
\begin{align*}
BFV^{(1,0)}(E,\Pi)\otimes BFV^{(0,1)}(E,\Pi) \xrightarrow{[-,-]_{BFV}} BFV^{0}(E,\Pi) \xrightarrow{proj} BFV^{(0,0)}(E,\Pi)
\end{align*}
coincides with the natural fibre pairing between $\mathcal{E}$ and $\mathcal{E}^{*}$.
\end{enumerate}
\end{Corollary}

\subsection{The BFV-Charge}\label{s:BFV_charge}
Next we construct a differential $D_{BFV}$ on
$BFV(E,\Pi)$ with special properties.

\begin{Proposition}\label{prop:charge}
Let $E\to S$ be a finite rank vector bundle over a smooth finite dimensional manifold $S$. Assume $(E,\Pi)$ is a Poisson manifold such that $S$ is a coisotropic submanifold.
Consider the graded Poisson algebra $BFV(E,\Pi):=\mathcal{C}^{\infty}(\mathcal{E}^{*}[1]\oplus \mathcal{E}[-1]),[-,-]_{BFV})$ with a bracket as in Corollary \ref{cor:BFV-bracket}.

Then there is an element $\Omega \in BFV(E,\Pi)$ of degree $+1$ such that
\begin{enumerate}
\item $[\Omega,\Omega]_{BFV}=0$ and
\item $\Omega \mbox{ mod } BFV_{\ge 1}(E,\Pi)$ is given by the tautological section of $\mathcal{E}\to E$.
\end{enumerate}
\end{Proposition}

Recall that $\mathcal{E}$ is the pullback bundle of $E \to S$ under $E \to S$ which admits
a tautological section. By the inclusions
\begin{align*}
\Gamma(\mathcal{E}) \hookrightarrow \Gamma(\wedge(\mathcal{E})) \hookrightarrow \Gamma(\wedge(\mathcal{E})\otimes \wedge(\mathcal{E^{*}})) = BFV(E,\Pi)
\end{align*}
the tautological section
can be seen as an element of $BFV^{(1,0)}(E,\Pi)$ which we denote by $\Omega_{0}$.

The proof we give is a slight adaptation of the arguments in \cite{Stasheff}:

\begin{proof}
It is convenient to work in local coordinates: fix local coordinates $(x^{\beta})_{\beta=1,\dots ,s}$ on $S$, linear
fibre coordinates $(y^{j})_{j=1,\dots,e}$ along $E$, $(c_{j})_{j=1,\dots,e}$ along $\mathcal{E}^{*}[1]$
and $(b^{j})_{j=1,\dots,e}$ along $\mathcal{E}[-1]$. In local coordinates the tautological section reads 
\begin{align*}
\Omega_{0}:=\sum_{j=1}^{e}y^{j}c_{j}. 
\end{align*}
Since $[\Omega_{0},\Omega_{0}]_{G}=0$ -- $G$ being the
Poisson bivector field given by the natural fibre pairing between $\mathcal{E}^{*}[1]$ and $\mathcal{E}[-1]$ --
we obtain a differential
\begin{align*}
\delta:=[\Omega_{0},-]_{G}=\sum_{j=1}^{e}y^{j}\frac{\vec{\partial}}{\partial b^{j}}.
\end{align*}

\underline{claim:} $H(BFV(E,\Pi),\delta)\cong \mathcal{C}^{\infty}(E^{*}[1])=\Gamma(\bigwedge E)$
\newline 
There are natural maps
\begin{eqnarray*} 
i:E^{*}[1] \hookrightarrow \mathcal{E}^{*}[1]\oplus \mathcal{E}[1] \mbox{ and}\\
p: \mathcal{E}^{*}[1]\oplus \mathcal{E}[-1] \to E^{*}[1].
\end{eqnarray*}
Define $h: BFV(E,\Pi) \to BFV(E,\Pi)[-1]$ by setting
\begin{multline*}
h(f_{j_{1}\dots j_{k}}(x,y,c)b^{j_{1}}\cdots b^{j_{k}}):=\\
\sum_{1\le \mu \le e}b^{\mu}\left( \int_{0}^{1} \frac{\partial f_{j_{1}\dots j_{k}}}{\partial y^{\mu}}(x,t\cdot y, c)t^{k}dt \right) b^{j_{1}}\cdots b^{j_{k}}
\end{multline*}
which is globally well-defined. It is straight-forward to check $i^{*}\circ \delta=0$, $\delta \circ p^{*}=0$, $h\circ h = 0$, $i^{*}\circ h =0$, $h \circ p^{*}=0$ and
$\delta \circ h + h \circ \delta = id - p^{*}\circ i^{*}$. It follows that
$i^{*}: BFV(E,\Pi) \to \mathcal{C}^{\infty}(E^{*}[1])$ induces and isomorphism on cohomology.

First note that $[\Omega_{0},\Omega_{0}]_{BFV} \mbox{ mod } BFV_{\ge 1}(E,\Pi)=[\Omega_{0},\Omega_{0}]_{\iota_{\nabla}(\Pi)}=:2 R_0$.
Using the biderivation property of $[-,-]_{\iota_{\nabla}(\Pi)}$ one sees that
\begin{align*}
[\Omega_{0},\Omega_{0}]_{\iota_{\nabla}(\Pi)}=[y^i,y^j]_{\iota_{\nabla}(\Pi)}c_i c_j + 2 y^i [c_i,y^j]_{\iota_{\nabla}(\Pi)} c_j + y^i y^j [c_i,c_j]_{\iota_{\nabla}(\Pi)}.
\end{align*}
Because $[y^i,y^j]_{\iota_{\nabla}(\Pi)}$ is equal to the pull back of $[y^i,y^j]_{\Pi}$ along the projection
$\mathcal{E}^{*}[1]\oplus \mathcal{E}[-1]\to E$, the condition that $[y^i,y^j]_{\Pi}$ is ontained in the vanishing ideal $\mathcal{I}_S$
of $S$ for abitrary $i,j=1,\dots,e$ is equivalent to the condition that $R_0$ vanishes when evaluated on $S$.
Hence the fact that $R_{0}$ vanishes along $S$ is equivalent to the fact that $S$
is coisotropic, see subsection \ref{s:Poisson}.

Using
$\delta([\Omega_{0},\Omega_{0}]_{\iota_{\nabla}(\Pi)})=0$ we obtain a cohomology class $[R_{0}] \in H(BFV(E,\Pi))\cong \mathcal{C}^{\infty}(N^{*}[1]S)$. Since the isomorphism between the two cohomologies is induced by setting the fibre coordinates $(y^{j})_{j=1,\dots,e}$ and $(b^{j})_{j=1,\dots,e}$ to zero
one sees that $[R_{0}]=0$. Hence $R_{0}= - \delta(\Omega_{1})$ for some $\Omega_{1} \in BFV_{1}(E,\Pi)$.
Consquently 
\begin{eqnarray*}
[\Omega_{0}+\Omega_{1},\Omega_{0}+\Omega_{1}]_{BFV} \mbox{ mod } BFV_{\ge 1}(E,\Pi) &=&
[\Omega_0,\Omega_0]_{\iota_{\nabla}(\Pi)} + [\Omega_0,\Omega_1]_G \\ &=& 2R_0 + \delta(\Omega_1)=0.
\end{eqnarray*}

\underline{claim:} Given $k>0$ and $\Omega(k):=\sum_{1\le i \le k}\Omega_{k}$ with $\Omega_{0}$ as above, $\Omega_{i} \in \Gamma(\bigwedge^{(i+1)}(\mathcal{E})\otimes \bigwedge^{i}(\mathcal{E}^{*}))$ and
\begin{align*}
[\Omega(k),\Omega(k)]_{BFV}=0 \mbox{ mod } BFV_{\ge k}(E,\Pi),
\end{align*}
there is an $\Omega_{k+1} \in BFV_{k+1}(E,\Pi)$
of total degree $+1$ such that $\Omega(k+1):=\Omega(k) + \Omega_{k+1}$ satisfies
\begin{align*}
[\Omega(k+1),\Omega(k+1)]_{BFV}=0 \mbox{ mod } BFV_{\ge(k+1)}(E,\Pi).
\end{align*}
Set $2R_{k}:=[\Omega(k),\Omega(k)]_{BFV} \mbox{ mod } BFV_{\ge(k+1)}(E,\Pi)$, hence $R_{k} \in BFV_{k}(E,\Pi)$.
By the graded Jacobi identity we know that $[\Omega(k),[\Omega(k),\Omega(k)]_{BFV}]_{BFV}=0$.
Moreover $[\Omega(k),\Omega(k)]_{BFV}=2R_{k} \mbox{ mod } BFV_{\ge k+1}(E,\Pi)$ implies that
\begin{eqnarray*}
0=[\Omega(k),[\Omega(k),\Omega(k)]_{BFV}]_{BFV}&=&
[\Omega_{0},2R_{k}]_{BFV} \mbox{ mod } BFV_{\ge k}(E,\Pi)\\
&=&\delta(2R_{k}).
\end{eqnarray*}
So $R_{k}$ is $\delta$-closed and using $H(BFV(E,\Pi),\delta)\cong \mathcal{C}^{\infty}(N^{*}[1]S)$
we can conclude that there is an element 
$\Omega_{k+1} \in BFV_{k+1}(E,\Pi)$ of total degree $+1$ such that $R_{k}=-\delta(\Omega_{k+1})$. It is easy to check that this element satisfies the conditions
of the claim.

After finitly many steps this procedure is finished thanks to the boundedness of the filtration $BFV_{\ge k}(E,\Pi)$.
The (well-defined) element 
\begin{align*}
\Omega:=\sum_{k \ge 0}\Omega_{k}
\end{align*}
satisfies properties 1. and 2. of the
Proposition by construction.
\end{proof}

\begin{Definition}\label{def:BFV}
Let $E\to S$ be a finite rank vector bundle over a smooth finite dimensional manifold. Assume $(E,\Pi)$ is a Poisson manifold and $S$ is a coisotropic submanifold.

Then a differential graded Poisson algebra algebra
$(BFV(E,\Pi),D_{BFV}:=[\Omega,-]_{BFV},[-,-]_{BFV})$ as constructed above is referred to as a \textsf{BFV-complex} associated to $(E,\Pi)$.
\end{Definition}

We remark that there are several BFV-complexes associated to $(E,\Pi)$. However in \cite{Schaetz} it is shown that different choices of a connection
on $E\to S$ and of the BFV-charge $\Omega$ yield isomorphic differential graded Poisson algebras.

\begin{Corollary}\label{cor:BFV-cohomology}
Let $E\to S$ be a finite rank vector bundle over a smooth finite dimensional manifold. Assume $(E,\Pi)$ is a Poisson manifold and $S$ a coisotropic submanifold.

The cohomology of $(BFV(E,\Pi),D_{BFV})$ is naturally isomorphic to the Lie algebroid
cohomology of $S$ introduced in subsection \ref{s:Poisson}.
\end{Corollary}

\begin{proof}
We use the filtration of $(BFV(E,\Pi),D_{BFV})$ given by $BFV^{(\ge q,\bullet)}(E,\Pi)$ to obtain a spectral sequence. Decomposing $D_{BFV}$ with respect to the degree $q$ yields
$\sum_{k\ge 0}\delta_{k}$ with  $\delta_{0}=\delta=[\Omega_{0},-]_{G}$. In the proof
of Proposition \ref{prop:charge} the isomorphism $H(BFV(E,\Pi),\delta)\cong \mathcal{C}^{\infty}(E^{*}[1])$
was established. This means that the spectral sequence under consideration collapses
after one step and so $H(BFV(E,\Pi),D_{BFV})$ is naturally isomorphic to the next sheet
of the spectral sequence. Hence we have to compute the cohomology of $\mathcal{C}^{\infty}(E^{*}[1]S)$
with respect to the induced differential to obtain $H(BFV(E,\Pi),D_{BFV})$.

It is straightforward to check that
the induced differential does not depend on the particular choice of $\Omega$ and that it
is given by the restriction of $\delta_{1}:=[\Omega_{0},-]_{\iota_{\nabla}(\Pi)}+[\Omega_{1},-]_G$ to 
$\mathcal{C}^{\infty}(E^{*}[1])=\Gamma(\bigwedge E)$.
A possible choice of $\Omega_{1}$ is given by $-h(1/2[\Omega_{0},\Omega_{0}]_{\iota_{\nabla}(\Pi)})$
with $h$ being the homotopy defined in the proof of Proposition \ref{prop:charge}.
In the local coordinates used in the proof of Proposition \ref{prop:charge} the induced differential
is given by 
\begin{align}\label{delta1}
[\delta_{1}]=c_{i}\left(\Pi^{i \beta}|_{S}\right) \frac{\vec{\partial}}{\partial x^{\beta}}-
\frac{1}{2}\left(\frac{\partial \Pi^{i j}}{\partial y^{k}}|_{S}\right)c_{i}c_{j}\frac{\vec{\partial}}{\partial c_{k}}
\end{align}
which coincides with the Lie algebroid
differential $\partial_{S}$. Hence the second sheet of the collapsing spectral sequence
associated to $(BFV(E,\Pi),D_{BFV})$ is equal to the Lie algebroid cocomplex $(\Gamma(\bigwedge E),\partial_{S})$ associated to $S$.
Consequently there is an isomorphism between $H(BFV(E,\Pi),D_{BFV})$ and the Lie algebroid cohomology
of $S$.
\end{proof}

In particular one obtains 
\begin{eqnarray*}
&& H^{0}(BFV(E,\Pi),D_{BFV})\cong \mathcal{C}^{\infty}(\underline{S})=\\
&& \quad \{f \in \mathcal{C}^{\infty}(S): X(f)=0 \mbox{ for all } X \in \Gamma(\mathcal{F}_{S})\}.
\end{eqnarray*}
Due to the compatibility between $D_{BFV}$ and the $BFV$-bracket $[-,-]_{BFV}$,
the cohomology $H(BFV(E,\Pi),D_{BFV})$ carries the structure of a graded Poisson algebra.
This structure restricts to the structure of a Poisson algebra on $H^{0}(BFV(E,\Pi),D_{BFV})\cong \mathcal{C}^{\infty}(\underline{S})$. It is easy to show that

\begin{Lemma}\label{l:resolution1}
The algebra isomorphism $H^{0}(BFV(E,\Pi),D_{BFV}) \cong \mathcal{C}^{\infty}(\underline{S})$ interwines the Poisson bracket induced from $[-,-]_{BFV}$ with $\{-,-\}_{\underline{S}}$ defined in subsection \ref{s:Poisson}.
\end{Lemma} 

Hence the BFV-complex
$(BFV(E,\Pi),D_{BFV},[-,-]_{BFV})$ is a resolution of the Poisson algebra
$(\mathcal{C}^{\infty}(\underline{S}),\{-,-\}_{\underline{S}})$.

%% file: Liealgebroid.tex
Let $S$ be a coisotropic submanifold of a smooth finite dimenional Poisson manifold $(M,\Pi)$. In section \ref{s:BFV} a differential graded Poisson algebra $(BFV(E,\Pi),D_{BFV},[-,-]_{BFV})$ was constructed such that the degree zero cohomology $H^{0}(BFV(E,\Pi),D_{BFV})$ is isomorphic to  $\mathcal{C}^{\infty}(\underline{S})$ as an algebra and the Poisson bracket induced from $[-,-]_{BFV}$ coincides with $\{-,-\}_{\underline{S}}$.

There is another resolution of the Poisson algebra
$(\mathcal{C}^{\infty}(\underline{S}),\{-,-\}_{\underline{S}})$ given by
the Lie algebroid cocomplex associated to $S$, enriched with compatible higher operations. 
This structure was found by Oh and Park (\cite{OhPark}) in the symplectic setting and called ``strong homotopy Lie algebroid'' there. It can also be derived
as the classical limit of the Poisson Sigma model with boundary conditions given by
$S$ (\cite{CattaneoFelder}). Our main aim is to show that the resolution
constructed by Oh and Park
is equivalent to $(BFV(E,\Pi),D_{BFV},[-,-]_{BFV})$ in the appropriate sense:
they are $L_{\infty}$ quasi-isomorphic, see Theorem \ref{thm:big} in subsection \ref{s:main}.
We remark that there is a connection between these resolutions and deformations of
$S$, see \cite{OhPark} and section \ref{s:deformation}. Moreover Kieserman showed in \cite{Kieserman} that the resolutions capture very
subtle properties of the foliation $\mathcal{F}_{S}:=\Pi^{\#}(N^{*}S)$ associated to $S$.

\subsection{The strong homotopy Lie Algebroid}\label{s:shLA}

We follow the presentation in \cite{CattaneoFelder} and \cite{Cattaneo} where the connection to the
derived bracket formalism (\cite{Voronov}) was made explicite.

Let $S \hookrightarrow M$ be a submanifold of a smooth finite dimensional Poisson manifold $(M,\Pi)$. By chosing an embedding of the normal bundle of $S$ as a tubular neighbourhood inside $M$ we obtain a finite rank vector bundle $E \xrightarrow{p} S$ equipped with a Poisson bivector field. We denote the embedding of $S$ into $E$ as the zero section by $i$. Abusing notation we denote the Poisson bivector field on $E$ by $\Pi$. We remark that there is a natural identification $E \cong NS$ ($NS$ being the normal bundle of $S$ in $E$).

There is a natural projection $pr: \mathcal{V}(E) \to \Gamma(\bigwedge E)$ given by the unique algebra morphism extending $f\mapsto f\circ i$ on
$\mathcal{C}^{\infty}(E)$ and 
\begin{align*}
\Gamma(TE) \to \Gamma(T_{S}E) \to \Gamma(E)
\end{align*}
where $E\to S$ is identified with the vertical part of $T_{S}E \to S$.
This projection admits a section
$s: \Gamma(\bigwedge E) \to \mathcal{V}(E)$: on functions $g\in \mathcal{C}^{\infty}(S)$ it is given by
$s(g):=g\circ p$ and on elements $X \in \Gamma(E)$ one defines $s(X)$ to be the unique vertical extension of $X$ that is constant along fibres of $E \to S$.

One checks that $s(\Gamma(\bigwedge E))\hookrightarrow \mathcal{V}(E)$ is an abelian Lie subalgebra of the
graded Lie algebra $(\mathcal{V}(E)[1],[-,-]_{SN})$. Moreover $ker(pr)[1]$ is a Lie subalgebra and
$\mathcal{V}(E)= ker(pr) \oplus s(\Gamma(\bigwedge E))$. Consequently 
\begin{align*}
(\mathcal{V}(E)[1],\Gamma(\wedge E)[1],pr[1])
\end{align*}
is a V-algebra (Definition \ref{V-algebra}).

The Poisson bivector field $\Pi$ on $E$ can be interpreted as a Maurer-Cartan element of
$(\mathcal{V}(E)[1],[-,-]_{SN})$. By Theorem \ref{Voronov} the derived brackets associated to the Poisson bivector field
\begin{align}\label{mu}
\hat{\mu}^{k}:=D^{k}_{\Pi}: (\Gamma(\wedge E)[1])^{\otimes k} \to \Gamma(\wedge E)[2]
\end{align}
define the structure of a (possibly non-flat) $L_{\infty}[1]$-algebra on $\Gamma(\bigwedge E)[1]$.
This corresponds to the structure of a (possibly non-flat) $L_{\infty}$-algebra on $\Gamma(\bigwedge E)$. We denote the structure maps of the $L_{\infty}$-algebra by $(\mu^{k})_{k\in\mathbb{N}}$.

The submanifold $S$ is coisotropic if and only if $pr(\Pi)=0$. In this case the $L_{\infty}$-algebra
is flat (i.e. the zero order component $\mu^{0} \in \Gamma(\bigwedge^{2}E)$ vanishes) and $\mu^{1}$ coincides with 
the Lie algebroid differential $\partial_{S}$ associated to $S$ (see subsection \ref{s:Poisson}). Hence:

\begin{Theorem}\label{extentedLiealgebroid}
Let $S$ be a coisotropic submanifold of a smooth finite dimensional Poisson manifold $(M,\Pi)$. Then
$(\Gamma(\bigwedge E),\partial_{S}=\mu^{1},\mu^{2},\cdots)$ constructed as above is
an $L_{\infty}$-algebra extending the Lie algebroid cocomplex associated to $S$.
\end{Theorem}

This Theorem first appeared in \cite{OhPark} in the symplectic setting.

\begin{Definition}\label{def:shLA}
Let $S$ be a coisotropic submanifold of a smooth finite dimensional Poisson manifold $(M,\Pi)$.
The \textsf{strong homotopy Lie algebroid} associated to $S$ is the $L_{\infty}$-algebra
$(\Gamma(\bigwedge NS),\partial_{S}=\mu^{1},\mu^{2},\cdots)$.
\end{Definition}

Since $(\mathcal{V}(E)[1],[-,-]_{SN})$ is a Gerstenhaber algebra and $pr$ and $s$ are morphisms of algebras one can
check that the structure maps $\mu^{k}$ are graded multiderivations with respect to the graded algebra structure: i.e.
\begin{eqnarray}\label{P_infty}
\begin{split}
\mu^{k}(a_{1}\otimes\cdots\otimes a_{k-1} \otimes a\cdot b) \ = \ \mu^{k}(a_{1}\otimes\cdots\otimes a_{k-1} \otimes a)\cdot b +&&\\
(-1)^{(|a_{1}|+\cdots+|a_{k-1}|+2-n)|a|}a\cdot\mu^{k}(a_{1}\otimes\cdots\otimes a_{k-1} \otimes b)&&
\end{split}
\end{eqnarray}
holds for all $k$ and arbitrary homogeneous elements $a_{1},\dots,a_{k-1},a,b$ of $\Gamma(\bigwedge E)$.
In \cite{CattaneoFelder} $L_{\infty}$-algebras on graded algebras with this property where called \textsf{$P_{\infty}$-algebras}.

We remark that the derived brackets $\mu^{k}$ depend in general on the choice of embedding $\phi: E \hookrightarrow M$. However it was proved in \cite{OhPark} in the symplectic case and in \cite{CattaneoSchaetz} in the Poisson case and for arbitrary submanifolds (not necessary coisotropic) that different choices lead to $L_{\infty}$-isomorphic $L_{\infty}$-algebras:

\begin{Theorem}
The $L_{\infty}$-algebra structures constructed on $\Gamma(\bigwedge(NS))$ with the help of two different
embeddings of $NS$ into $M$ as tubular neighbourhoods of $S$ are $L_{\infty}$-isomorphic.
\end{Theorem}

Let $S$ be a coisotropic submanifold of $(M,\Pi)$. By Theorem \ref{extentedLiealgebroid} there is a nontrivial
extension of the Lie algebroid cocomplex $(\Gamma(\bigwedge NS),\partial_{S})$ associated to $S$ to an
$L_{\infty}$-algebra. As observed in subsection \ref{s:Poisson} the zero Lie algebroid cohomology
$H^{0}(\Gamma(\bigwedge NS),\partial_{S})$ is given by $\mathcal{C}^{\infty}(\underline{S})$.
The binary operation $\mu^{2}$ descends to cohomology where it induces a Lie bracket. Since
$\mu^{2}$ is a graded biderivation with respect to the graded algebra structure the induced Lie bracket
will be a biderivation, i.e. $\mathcal{C}^{\infty}(\underline{S})$ inherits a Poisson bracket. A computation shows that

\begin{Lemma}\label{l:resolution2}
The algebra isomorphism $H^{0}(\Gamma(\bigwedge NS),\partial_{S})\cong \mathcal{C}^{\infty}(\underline{S})$
interwines the Poisson bracket induced from $\mu^{2}$ with $\{-,-\}_{\underline{S}}$ as defined in subsection
\ref{s:Poisson}.
\end{Lemma}

Consequently the $P_{\infty}$-algebra $(\Gamma(\bigwedge NS),\partial_{S}=\mu^{1},\mu^{2},\dots)$ is a resolution of the
Poisson algebra $(\mathcal{C}^{\infty}(\underline{S}),\{-,-\}_{\underline{S}})$.

\subsection{Relation of the two Resolutions}\label{s:main}

Let $S$ be a coisotropic submanifold of a smooth finite dimensional Poisson manifold $(M,\Pi)$.
Lemma \ref{l:resolution1} in subsection \ref{s:BFV_charge} established that the differential graded Poisson algebra
$(BFV(E,\Pi),D_{BFV},[-,-]_{BFV})$ can be interpreted as a homological resolution of the Poisson algebra
$(\mathcal{C}^{\infty}(\underline{S}),\{-,-\}_{\underline{S}})$ introduced in subsection \ref{s:Poisson}.
The same is true for the strong homotopy Lie algebroid $(\Gamma(\bigwedge E),\partial_{S}=\mu^{1},\mu^{2},\dots)$ constructed in subsection \ref{s:shLA} (see Lemma \ref{l:resolution2}).
Moreover Corollary \ref{cor:BFV-cohomology} in subsection \ref{s:BFV_charge} established an isomorphism of graded algebras
$H^{\bullet}(BFV(E,\Pi),D_{BFV})\cong H^{\bullet}(\Gamma(\bigwedge E),\partial_{S})$.

A natural question to ask is how tight the connection between the BFV-complex $(BFV(E,\Pi),D_{BFV},[-,-]_{BFV})$ and the $P_{\infty}$-algebra $(\Gamma(\bigwedge E),\partial_{S}=\mu^{1},\mu^{2},\dots)$ actually is. We provide an answer to this question:

\begin{Theorem}\label{thm:big}
Let $E\to S$ be a finite rank vector bundle over a smooth finite dimensional manifold $S$. Assume $(E,\Pi)$ is a Poisson manifold such that $S$ is a coisotropic submanifold.

Then there is an $L_{\infty}$ quasi-isomorphism between the BFV-complex \newline $(BFV(E,\Pi),D_{BFV},[-,-]_{BFV})$ associated to $S$ (Definition \ref{def:BFV} in subsection \ref{s:BFV_charge}) and the strong homotopy Lie algebroid associated to $S$, i.e. $(\Gamma(\bigwedge E),\partial_{S}=\mu^{1},\mu^{2},\dots)$ (Definition \ref{def:shLA}).
\end{Theorem}

An immediate consequence of Theorem \ref{thm:big} is:

\begin{Corollary}\label{cor:formaldef}
Let $E\to S$ be a finite rank vector bundle over a smooth finite dimensional manifold $S$. Assume $(E,\Pi)$ is a Poisson manifold such that $S$ is a coisotropic submanifold.

Then the formal deformation problems associated to the BFV-complex $(BFV(E,\Pi),D_{BFV},[-,-]_{BFV})$
and to the strong homotopy Lie algebroid \newline $(\Gamma(\bigwedge E),\partial_{S}=\mu^{1},\mu^{2},\dots)$
are equivalent.
\end{Corollary}

Next we prove Theorem \ref{thm:big}:

\begin{proof}
The strategy of the proof is as follows: the starting point is the BFV-complex $(BFV(E,\Pi),D_{BFV},[-,-]_{BFV})$.
As a by-product of the proof of Proposition \ref{prop:charge} in subsection \ref{s:BFV_charge} we obtained an isomorphism of graded algebras $H^{\bullet}(BFV(E,\Pi),\delta)\cong \Gamma(\bigwedge^{\bullet} E)$. The $\bullet$
on the left hand side refers to the grading with respect to the total degree. Recall that $\delta$ is $[\Omega_{0},-]_{G}$ where $\Omega_{0} \in BFV(E,\Pi)$ is given by the tautological section of the bundle $\mathcal{E} \to E$ and
$G$ denotes the Poisson bivector field on $\mathcal{E}^{*}[1]\oplus \mathcal{E}[-1]$ representing the fibre pairing between $\mathcal{E}^{*}[1]$ and $\mathcal{E}[-1]$.

More explicitly we considered pullbacks $i^{*}$ and $p^{*}$ along $i:E^{*} \hookrightarrow \mathcal{E}^{*}[1]\oplus \mathcal{E}[1]$ and $p: \mathcal{E}^{*}[1]\oplus \mathcal{E}[-1] \to E^{*}[1]$ and a homotopy
$h: BFV(E,\Pi)\to BFV(E,\Pi)[-1]$ such that $h\circ h = 0$, $i^{*}\circ h =0$, $h \circ p^{*}=0$ and
$\delta \circ h + h \circ \delta = id - p^{*}\circ i^{*}$ hold. We summarize the situation in the following diagram:
\begin{align}\label{transferBFV}
\xymatrix{(\mathcal{C}^{\infty}(E^{*}[1]),0) \ar@<2pt>[r]^<<<<<{p^{*}}& \ar@<2pt>[l]^>>>>>{i^{*}} (BFV(E,\Pi),\delta),h.}
\end{align}
By Theorem \ref{thm:transfer} in subsection \ref{s:homotopytransfer} these data can be used for homological transfer of an $L_{\infty}$-algebra structure from
$(BFV(E,\Pi),\delta)$ to $\mathcal{C}^{\infty}(E^{*}[1])=\Gamma(\bigwedge E)$.

We will use these data to perform the homological transfer of the differential graded Lie algebra $(BFV(E,\Pi),D_{BFV},[-,-]_{BFV})$ to $\Gamma(\bigwedge E)$
in terms of diagrams as described in subsection \ref{s:homotopytransfer}. It will turn out that no convergence issues arise
and that the induced $L_{\infty}$-algebra structure on $\Gamma(\bigwedge E)$ is a $P_{\infty}$-algebra, i.e. the structure maps are graded multiderivations. Hence we have two $P_{\infty}$-algebra structures on $\Gamma(\bigwedge E)$: one is induced from the BFV-complex and the second onea \textsf{BFV-complex} associated to $S$ comes from the strong homotopy Lie algebroid associated to $S$. Since $\Gamma(\bigwedge E)$ is generated by $\mathcal{C}^{\infty}(S)$ and $\Gamma(E)$ as a graded algebra it suffices to know the structure maps of the $P_{\infty}$-algebra structures restricted to $\mathcal{C}^{\infty}(S)$ and $\Gamma(E)$ respectively in order to be able to reconstruct them completely. We will check that the restricted structure maps of the two $P_{\infty}$-algebras coincide, hence so do the full $P_{\infty}$-algebras.

\vspace{0.2cm}
\underline{Step 1) homological transfer in terms of trees:}
\newline
We perform the homological transfer of the differential graded Lie algebra structure on $BFV(E,\Pi)$ along the diagram \eqref{transferBFV}. What does the induced $L_{\infty}$-algebra structure on $\Gamma(\bigwedge E)$ look like?

The BFV-differential $D_{BFV}=[\Omega,-]_{BFV}$ can be decomposed as $D_{BFV}=\delta + D_{R}$ satisfying $\delta \circ \delta = 0$
and $\delta \circ D_{R} + D_{R} \circ \delta + D_{R} \circ D_{R}=0$. Recall that $BFV(E,\Pi)$ carries a bigrading given by $BFV^{(p,q)}(E,\Pi):=\Gamma(\bigwedge^{p}(\mathcal{E})\otimes \bigwedge^{q}(\mathcal{E}^{*}))$. We have
\begin{eqnarray*}
&& h: BFV^{(p,q)}(E,\Pi) \to BFV^{(p,q+1)}(E,\Pi),\\
&& D_{R}: BFV^{(p,q)}(E,\Pi) \to \bigoplus_{(p'>p,q'\ge q, p'-q'=p-q)}BFV^{(p',q')}(E,\Pi) \mbox{ and}\\
&& [-,-]_{BFV}=[-.-]_{G}+[-,-]_{\iota_{\nabla}(\Pi)}\mbox{ mod }BFV_{\ge 1}(E,\Pi).
\end{eqnarray*}

Following subsection \ref{s:homotopytransfer} the induced structure maps are given in terms of oriented trees with edges decorated by non-negative integers. The set of exterior vertices decomposes into the set of leaves (with edges pointing away from them) and a unique root (with an edge pointing towards it). To each such decorated tree $T$
a map 
\begin{align*}
m_{T}:=\Gamma(\bigwedge E)^{\otimes \#(leaves)} \to \Gamma(\bigwedge E)
\end{align*}
is associated by the following rule:
put $[-,-]_{BFV}$ at the trivalent vertices and $k$ copies of $D_{R}$
at edges decorated by the number $k$. Between consecutive operations $[-,-]_{BFV}$ or $D_{R}$ place a homotopy $-h$. We define $\tilde{m}_{T}: BFV(E,\Pi)^{\otimes\#(leaves)} \to BFV(E,\Pi)$ to be the composition of all these maps in the order given by the orientation of the tree $T$. Then we set $m_{T}:=i^{*} \circ \tilde{m}_{T} \circ (p^{*})^{\otimes \#(leaves)}$.

Because $p^{*}(\Gamma(\bigwedge E)) \subset BFV^{(\bullet,0)}(E,\Pi)$ and
$(i^{*})^{-1}(\Gamma(\bigwedge E)) \subset BFV^{(\bullet,0)}(E,\Pi)$ the operation $m_{T}$ associated to a 
decorated tree $T$ can only
be non-zero if the corresponding $\tilde{m}_{T}$ maps the subspace $(BFV^{(\bullet,0)}(E,\Pi))^{\otimes \#(leaves)}$ to
a subspace having nonvanishing intersection with $BFV^{(\bullet,0)}(E,\Pi)$.

Since the homotopy $h$ increases the ghost-momentum degree by $1$ and $[-,-]_{G}$ is the only operation the decreases it by $1$, there must be at least as many trivalent vertices decorated by $[-,-]_{G}$ as there are $h$. From
\begin{align*}
\#([-,-]_{G})\ge \#(h)=\#(D_{R})+\#(\mbox{trivalent vertices}) - 1
\end{align*}
it follows that 
\begin{align*}
\#(D_{R}) + \#(\mbox{trivalent vertices decorated by } ([-,-]_{BFV})-([-,-]_{G})) \le 1.
\end{align*}
One can easily exclude the sharp inequality so there
are two remaining cases: either $1)$ all the edges of the tree are decorated by zeros. In this case exactly one of the trivalent vertices is decoratd by $([-,-]_{BFV})-([-,-]_{G})$ and the other trivalent vertices are decorated by $[-,-]_{G}$.
Or $2)$ exactly one edge is decorated by $1$ and all the others by zero. In this case all of the trivalent vertices are decorated by $[-,-]_{G}$.

Observe that in both case $1)$ and $2)$ the part of the ``exceptional'' operation $D_{R}$ and $([-,-]_{BFV})-([-,-]_{G})$ respectively that actually contributes to $m_{T}$ is the part of ghost-momentum degree $0$. We decompose $D_{BFV}$ with respect to the ghost degree $D_{BFV}=\sum_{k \ge 0} \delta_{k}$ with $\delta_{0}=\delta$, and hence $D_{R}=\sum_{k \ge 1}\delta_{k}$. The fact that $D_{R}$ is of total degree $1$ implies that its component of ghost-momentum degree $0$ is given by $\delta_{1}$. The ghost-momentum degree $0$ component of $([-,-]_{BFV})-([-,-]_{G})$ is $[-,-]_{\iota_{\nabla}(\Pi)}$. 

Moreover the identity $[p^{*}(-),p^{*}(-)]_{G}=0$ holds because $p^{*}(\Gamma(\bigwedge E))\subset BFV^{(\bullet,0)}(E,\Pi)$ and because $[-,-]_{G}$ is the graded Poisson bracket induced from the fibre pairing between $\mathcal{E}^{*}[1]$ and $\mathcal{E}[-1]$. Hence the only two types of trees that contribute to the induced
$L_{\infty}$-algebra structure on $\Gamma(\bigwedge E)$ are the following:
\vspace{0.2cm}
\newline
\begin{minipage}[c]{\textwidth}
\centering \includegraphics[width=10cm]{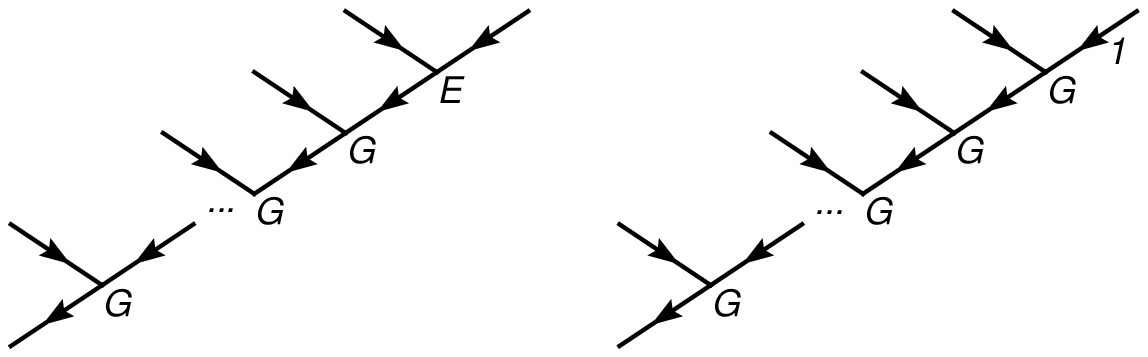}
\vspace{0.2cm}
\end{minipage}
Here the decoration $E$ refers to $[-,-]_{\iota_{\nabla}(\Pi)}$, $G$ refers to $[-,-]_{G}$ and the decoration of the edges was left out whenever it is zero. We denote the maps from $(\Gamma(\bigwedge E))^{\otimes n}$ to $\Gamma(\bigwedge E)$ associated to the trees on the left / right hand side with $n$ leaves by $L_{n}$ and $R_{n}$ respectively. Up to skew-symmetrization and sign issues these two families of maps define the induced $L_{\infty}$-algebra structure on $\Gamma(\bigwedge E)$.

\vspace{0.2cm}
\underline{Step 2) $P_{\infty}$-property:}
\newline
The $L_{\infty}$-algebra structure $(\Gamma(\bigwedge E), \partial_{S}=\mu^{1},\mu^{2},\dots)$ satisfies the
$P_{\infty}$ property \eqref{P_infty} as remarked before. Furthermore

\begin{Lemma}\label{inducedP}
The $L_{\infty}$-algebra structure on $\Gamma(\bigwedge E)$ induced from the differential graded
Poisson algebra $(BFV(E,\Pi),D_{BFV},[-,-]_{BFV})$ satisfies the $P_{\infty}$ property \eqref{P_infty}.
\end{Lemma}

\begin{proof}
We first prove that the result of the evaluation of
$L_{n}$ ($R_{n}$) on elements of the form
\begin{align*}
a_{1}\otimes\cdots \otimes a_{k-1}\otimes a\cdot b \otimes a_{k} \cdots \otimes a_{n-1} \in \Gamma(\bigwedge E)^{\otimes n}.
\end{align*}
can be expressed using $L_{n}$ ($R_{n}$) evaluated on $a_{1}\otimes \cdots a \cdots \otimes a_{n-1}$ and on $a_{1}\otimes\cdots \otimes b \otimes a_{n-1}$ only.
Without loss of generality one may assume that $a_{1}, \cdots,a_{(n-1)},a,b$ are all homogeneous.

Consider the map $L_{n}$ first and assume that $k < (n-1)$. By the graded Leibniz identity for $[-,-]_{G}$ we have
\begin{multline*}
[p^{*}(a\cdot b),\bullet]_{G}=[p^{*}(a)\cdot p^{*}(b),\bullet]_{G}=\\
p^{*}(a)\cdot[p^{*}(b),\bullet]_{G} + (-1)^{|a||b|}p^{*}(b)\cdot [p^{*}(a),\bullet]_{G}. 
\end{multline*}
Recall the definition
of the homotopy $h$ given during the proof of Proposition \ref{prop:charge} in subsection \ref{s:BFV_charge}:
\begin{multline*}
h(f_{\mu_{1}\dots \mu_{k}}(x,y,c)b^{\mu_{1}}\cdots b^{\mu_{k}}):=\\
\sum_{1\le \mu \le s}b^{\mu}\left( \int_{0}^{1} \frac{\partial f_{\mu_{1}\dots \mu_{k}}}{\partial y^{\mu}}(x,t\cdot y, c)t^{k}dt \right) b^{\mu_{1}}\cdots b^{\mu_{k}}.
\end{multline*}
Hence $h(p^{*}(X)\cdot Y)= (-1)^{|X|}p^{*}(X) \cdot h(Y)$ because $p^{*}X$ does not depend on the coordinates $y^{\mu}$ and $b^{\mu}$. So
\begin{multline}\label{ex1}
h([p^{*}(a\cdot b),\bullet]_{G})=\\
(-1)^{|a|}p^{*}(a)\cdot h([p^{*}(b),\bullet]_{G})+(-1)^{(|a|+1)|b|}p^{*}(b)\cdot h([p^{*}(a),\bullet]_{G})
\end{multline}
holds. Applying consecutivly 
\begin{enumerate}
\item $[p^{*}(-),-]_{G}$ and using the graded Leibniz identity together with $[p^{*}(-),p^{*}(-)]_{G}=0$; and 
\item $h$ and using $h(X\cdot p^{*}(Y))=h(X)\cdot p^{*}(Y)$
\end{enumerate}
leads to
\begin{multline*}
L_{n}(a_{1}\otimes\cdots \otimes a_{k-1}\otimes a\cdot b \otimes a_{k} \cdots \otimes a_{n-1})=\\
\begin{gathered}
(-1)^{(|a_{1}|+\cdots+|a_{k-1}|+k)|a|}a\cdot L_{n}(a_{1}\otimes\cdots \otimes a_{k-1}\otimes b \otimes a_{k} \cdots \otimes a_{n-1}) + \hfill \\
(-1)^{(|a_{1}|+\cdots+|a_{k-1}|+|a|+k)|b|}b\cdot L_{n}(a_{1}\otimes\cdots \otimes a_{k-1}\otimes a \otimes a_{k} \cdots \otimes a_{n-1})
\end{gathered}
\end{multline*}
for $k < (n-1)$. By similar reasoning this formula can be extented to the cases $k=(n-1)$ and $k=n$.

We claim that
\begin{multline*}
R_{n}(a_{1}\otimes\cdots \otimes a_{k-1}\otimes a\cdot b \otimes a_{k} \cdots \otimes a_{n-1})=\\
\begin{gathered}
(-1)^{(|a_{1}|+\cdots+|a_{k-1}|+k)|a|}a\cdot R_{n}(a_{1}\otimes\cdots \otimes a_{k-1}\otimes b \otimes a_{k} \cdots \otimes a_{n-1}) + \hfill \\
(-1)^{(|a_{1}|+\cdots+|a_{k-1}|+|a|+k)|b|}b\cdot R_{n}(a_{1}\otimes\cdots \otimes a_{k-1}\otimes a \otimes a_{k} \cdots \otimes a_{n-1})
\end{gathered}
\end{multline*}
holds as well. For $k < n$ the arguments previously applied to $L_{n}$ go through. For the case $k=n$ we make use of the
explicite formula for $\delta_{1}$ which was derived in the proof of Corollary \ref{cor:BFV-cohomology} in subsection \ref{s:BFV_charge}:
\begin{align*}
\delta_{1}=[\Omega_{0},-]_{\iota_{\nabla}(\Pi_{\phi})}+[\Omega_{1},-]_G.
\end{align*}
Hence $\delta_{1}(p^{*}(a\cdot b))=\delta_{1}(p^{*}(a)) \cdot p^{*}(b) + (-1)^{|a|}p^{*}(a) \cdot \delta_{1}(p^{*}(b))$ and applying the established computation rules for $h$ and $[p^{*}(-),-]_{G}$ yields the claimed formula for $R_{n}$.

If one takes the signs arising from the d\'ecalage-isomorphism and graded symmetrization into account one obtains the signs as stated in \eqref{P_infty}.
\end{proof}

\vspace{0.2cm}
\underline{Step 3) localization:}
\newline
The graded commutative associative algebra $\mathcal{C}^{\infty}(E^{*}[1]) = \Gamma(\bigwedge E)$ is
generated by elements of degree $0$ and $1$, i.e. by $\mathcal{C}^{\infty}(S)$ and $\Gamma(E)$.
Hence it is enough to know $L_{n}$ and $R_{n}$ restricted
to $(\mathcal{C}^{\infty}(S)\oplus\Gamma(E))^{\otimes n}\subset (\Gamma(\bigwedge E))^{\otimes n}$ by Lemma \ref{inducedP}.
Since $\Gamma(\bigwedge E)$ is concentrated in non-negative degrees and the total degree of
$L_{n}$ and $R_{n}$ is $(2-n)$, it suffices to know $L_{n}$ and $R_{n}$ on elements of one of the following
types:
\renewcommand{\labelenumi}{\Alph{enumi}')}
\begin{enumerate}
\item $\gamma^{1}\otimes \cdots \otimes \gamma^{n}$ for $\gamma^{i} \in \Gamma(E)$,
\item $\gamma^{1}\otimes \cdots \otimes \gamma^{(k-1)} \otimes f \otimes \gamma^{k} \otimes \cdots \otimes \gamma^{(n-1)}$ for $\gamma^{i} \in \Gamma(E)$, $f\in \mathcal{C}^{\infty}(S)$,
\item $\gamma^{1} \otimes \cdots \otimes \gamma^{(k-1)} \otimes f \otimes \gamma^{k} \otimes \cdots \otimes \gamma^{(k+l-1)} \otimes g \otimes \gamma^{(k+l)} \otimes \cdots \otimes \gamma^{(n-2)}$ for
$\gamma^{i} \in \Gamma(E)$, $f,g \in \mathcal{C}^{\infty}(S)$.
\end{enumerate}
We choose a trivializing cover $\mathcal{U}:=\{U_{\alpha}\}_{\alpha \in A}$ for the vector bundle $E\to S$.
Let $\{\rho_{\alpha}\}_{\alpha \in A}$ be a partition of unity subordinated to $\mathcal{U}$, i.e.
$a)$ $\rho_{\alpha} \in \mathcal{C}^{\infty}(S)$,$b)$ $supp(\rho_{\alpha})\subset U_{\alpha}$ for every
$\alpha \in A$, $c)$ $\{\rho_{\alpha}\}_{\alpha \in A}$ is locally finite (for every $x \in S$ there is an open neighbourhood $U$ such that there are only
finitly many $\alpha \in A$ with $\rho_{\alpha}|_{U}\neq 0$) and $d)$ $\sum_{\alpha \in A} \rho_{\alpha}=1$.

For an arbitrary $f \in \mathcal{C}^{\infty}{S}$ we write
\begin{align*}
f=(\sum_{\alpha \in A}\rho_{\alpha})f=\sum_{\alpha\in A}(\rho_{\alpha}f)=:\sum_{\alpha\in A}f_{\alpha}
\end{align*}
where $f_{\alpha}$ is supported on $U_{\alpha}$. Similarly we get $\gamma=\sum_{\alpha \in A}\gamma_{\alpha}$ for
any section $\gamma \in \Gamma(E)$. Since $\mathcal{U}$ is a collection of trivializing neighbourhoods
of the vector bundle $E$ we can choose a local frame $(e^{\alpha}_{1},\dots,e^{\alpha}_{s})$ of $E$ restricted
to $U_{\alpha}$. The section $\gamma_{\alpha}$ is supported on $U_{\alpha}$ and hence there are local functions
$(w_{\alpha}^{1},\dots,w_{\alpha}^{r})$ such that
\begin{align*}
\gamma_{\alpha}=\sum_{j=1}^{s}w_{\alpha}^{j}e^{\alpha}_{j}.
\end{align*}
Using this decomposition of smooth functions and sections of $E$ on elements of
$(\mathcal{C}^{\infty}(S)\oplus \Gamma(E))^{\otimes n}$ of type A'), B') or C') shows that
$L_{n}$ and $R_{n}$ are totally determined by evaluating them for abitrary $\alpha \in A$ on elements of the form
\renewcommand{\labelenumi}{\Alph{enumi})}
\begin{enumerate}
\item $e^{\alpha}_{j_{1}}\otimes \cdots \otimes e^{\alpha}_{j_{n}}$,
\item $e^{\alpha}_{j_{1}}\otimes \cdots \otimes e^{\alpha}_{j_{(k-1)}} \otimes f \otimes e^{\alpha}_{j_{k}} \otimes \cdots \otimes e^{\alpha}_{j_{(n-1)}}$ with $f\in \mathcal{C}^{\infty}(U_{\alpha})$,
\item $e^{\alpha}_{j_{1}}\otimes \cdots \otimes e^{\alpha}_{j_{(k-1)}} \otimes f \otimes e^{\alpha}_{j_{k}} \otimes \cdots \otimes e^{\alpha}_{j_{(k+l-1)}} \otimes g \otimes e^{\alpha}_{j_{(k+l)}} \otimes \cdots \otimes e^{\alpha}_{j_{(n-2)}}$ with $f,g \in \mathcal{C}^{\infty}(U_{\alpha})$.
\end{enumerate}
\renewcommand{\labelenumi}{\arabic{enumi})}
Since we only used the $P_{\infty}$ property and the total degrees of the structure maps $L_{n}$ and $R_{n}$, the same is true for the structure maps $\mu^{n}$ of the strong homotopy Lie algebroid $(\Gamma(\bigwedge E),\partial_{S}=\mu^{1},\mu^{2},\dots)$ associated to $S$.

\vspace{0.2cm}
\underline{Step 4) comparison of the restricted structure maps:}
\newline
Let $U_{\alpha}$ be an open subset of the trivializing cover $\mathcal{U}$. The aim is to compute explicite coordinate expressions on $U_{\alpha}$ for the restricted structure maps of the strong homotopy Lie algebroid and the structure induced from the BFV-complex respectively.

Let $\{x^{\beta}\}_{\beta=1,\dots s}$ be coordinates for $S$ and $\{y^{j}\}_{j=1,\dots e}$ linear fibre coordinates along $E$. We have to consider the graded Lie algebra $\mathcal{V}(E|_{U_{\alpha}})[1]$ with the bracket given by
\begin{align*}
[\frac{\partial}{\partial x^{\alpha}},x^{\beta}]_{SN}=\delta_{\alpha}^{\beta}, \quad
[\frac{\partial}{\partial y^{i}},y^{j}]_{SN}=\delta_{i}^{j}.
\end{align*}
The Poisson bivector field $\Pi$ is given by
\begin{align*}
\frac{1}{2}\Pi^{\alpha\beta}\frac{\partial}{\partial x^{\alpha}}\frac{\partial}{\partial x^{\beta}} + \Pi^{\alpha j}\frac{\partial}{\partial x^{\alpha}}\frac{\partial}{\partial y^{j}} + \frac{1}{2}\Pi^{i j}\frac{\partial}{\partial y^{i}}\frac{\partial}{\partial y^{j}}.
\end{align*}
A straightforward computation of the restricted structure maps $\hat{\mu}^{k}$ of the $L_{\infty}[1]$-algebra structure on $\Gamma(\bigwedge E|_{U_{\alpha}})[1]$ yields
\begin{eqnarray*}
&&\hspace{-0.7cm}\hat{\mu}^{k}(\frac{\partial}{\partial y^{j_{1}}}\otimes \cdots \otimes \frac{\partial}{\partial y^{j_{k}}})=(-1)^{k}\frac{1}{2}
\left(\frac{\partial}{\partial y^{j_{1}}}\cdots \frac{\partial}{\partial y^{j_{k}}}\left(\Pi^{i l}\frac{\partial}{\partial y^{i}}\frac{\partial}{\partial y^{l}}\right)\right)|_{S}\\
&&\hspace{-0.7cm}\hat{\mu}^{k}(\frac{\partial}{\partial y^{j_{1}}}\otimes \cdots \otimes \frac{\partial}{\partial y^{j_{(k-1)}}}\otimes f(x))= (-1)^{k}\left(\frac{\partial}{\partial y^{j_{1}}}\cdots \frac{\partial}{\partial y^{j_{(k-1)}}}\left(\Pi^{\alpha l}\frac{\partial f(x)}{\partial x^{\alpha}}\frac{\partial}{\partial y^{l}}\right)\right)|_{S}\\
&&\hspace{-0.7cm}\hat{\mu}^{k}(\frac{\partial}{\partial y^{j_{1}}}\otimes \cdots \otimes \frac{\partial}{\partial y^{j_{(k-2)}}}\otimes f(x) \otimes g(x))=\\
&& \hspace{2cm} (-1)^{(k-1)}\left(\frac{\partial}{\partial y^{j_{1}}}\cdots \frac{\partial}{\partial y^{j_{(k-2)}}}\left(\Pi^{\alpha \beta }\frac{\partial f(x)}{\partial x^{\alpha}}\frac{\partial g(x)}{\partial x^{\beta}}\right)\right)|_{S}
\end{eqnarray*}
Only the last expression picks up a sign under the d\'ecalage-isomorphism: the exponent changes from $(k-1)$ to $k$.

To obtain concrete formulae for the induced $L_{\infty}$-algebra structure we first make some general observations on the induced structure maps. All the operations $D_{R}$, $h$, $[-,-]_{G}$ and $[-,-]_{\iota_{\nabla}(\Pi)}$ are (multi-)differential operators and the surjection from $BFV(E,\Pi)$ to its cohomology $\Gamma(\bigwedge E)$ involves the evaluation of sections at $S \hookrightarrow E$. It follows that the induced structure maps only depend on the jet-expansion of $\Pi$ in transversal directions and that the homotopy $h$ can be replaced by its jet-version. For convenience let us introduce the following local coordinates: $(x^{\beta})_{\beta=1,\dots ,s}$ on $S$, linear
fibre coordinates $(y^{j})_{j=1,\dots,e}$ along $E$, $(c_{j})_{j=1,\dots,e}$ along $\mathcal{E}^{*}[1]$
and $(b^{j})_{j=1,\dots,e}$ along $\mathcal{E}[-1]$. In these local coordinates the jet-version of the homotopy reads
\begin{eqnarray*}
&&\hat{h}(f_{j_{1}\dots j_{k}}(x,y,c)b^{j_{1}}\cdots b^{j_{k}}):=\\
&&\sum_{1\le \mu \le e}\frac{1}{N(f)+k}b^{\mu}\left(\frac{\partial f_{j_{1}\dots j_{k}}}{\partial y^{\mu}}(x, y, c) \right) b^{j_{1}}\cdots b^{j_{k}}
\end{eqnarray*}
where $N(f)$ is the polynomial degree of $f$ with respect to the transverse directions $\{y^{j}\}_{j=1,\dots e}$.

In local coordinates the horizontal lift $\iota_{\nabla}(\Pi)$ of $\Pi$ is given by
\begin{eqnarray*}
&&\frac{1}{2}\Pi^{\alpha\beta}\left(\frac{\partial}{\partial x^{\alpha}}+\Gamma_{\alpha r}^{s}c_{s}\frac{\partial}{\partial c_{r}}-\Gamma_{\alpha r}^{s}b^{r}\frac{\partial}{\partial b^{s}}\right)\left(\frac{\partial}{\partial x^{\beta}}+\Gamma_{\beta m}^{n}c_{n}\frac{\partial}{\partial c_{m}}-\Gamma_{\beta m}^{n}b^{m}\frac{\partial}{\partial b^{n}}\right)\\ 
&&+\Pi^{\alpha j}\left(\frac{\partial}{\partial x^{\alpha}}+\Gamma_{\alpha r}^{s}c_{s}\frac{\partial}{\partial c_{r}}-\Gamma_{\alpha r}^{s}b^{r}\frac{\partial}{\partial b^{s}}\right)\frac{\partial}{\partial y^{j}} + 
\frac{1}{2}\Pi^{i j}\frac{\partial}{\partial y^{i}}\frac{\partial}{\partial y^{j}}.
\end{eqnarray*}
Here $\Gamma$ denotes the Christoffel symbols of the pull back connection on $\mathcal{E}[1]\oplus \mathcal{E}^{*}[-1]$. Moreover the restriction of
\begin{align*}
\delta_{1}(-)=[\Omega_{0},-]_{\iota_{\nabla}(\Pi)}+[\Omega_{1},-]_{G}
\end{align*}
(with $\Omega_{1}:=-\frac{1}{2}h([\Omega_{0},\Omega_{0}]_{\iota_{\nabla}(\Pi)})$)
to $\Gamma(E) \hookrightarrow BFV(E,\Pi)$ reads
\begin{eqnarray*}
&&\Gamma_{\alpha r}^{s}y^{r}c_{s}\Pi^{\alpha \beta}\left(\frac{\partial}{\partial x^{\beta}}+\Gamma_{\beta m}^{n}c_{n}\frac{\partial}{\partial c_{m}}\right) + c_{m}\Pi^{m \alpha}\left(\frac{\partial}{\partial x^{\alpha}}+\Gamma_{\alpha r}^{s}c_{s}\frac{\partial}{\partial c_{r}}\right)\\
&&-\frac{\partial}{\partial b^{\mu}}\left(\hat{h}\left(\frac{1}{2}\Pi^{\alpha\beta}\Gamma_{\alpha r}^{s}y^{r}c_{s}\Gamma_{\beta m}^{n}y^{m}c_{n}+\Pi^{\alpha k}\Gamma_{\beta r}^{s}y^{r}c_{s}c_{k}+\frac{1}{2}\Pi^{ij}c_{i}c_{j}\right)\right)\frac{\partial}{\partial c_{\mu}}.
\end{eqnarray*}
A straightforward but lenghty calculation of the restricted structure maps of the induced $P_{\infty}$-algebra structure shows that all contributions involving Christoffel-symbols cancel each other and that the formulae reduce to the local expressions for $\mu^{k}$ we derived above.

\end{proof}

%% file: Deformation.tex
A relation between BFV-complexes (see Definition \ref{def:BFV} in subsection \ref{s:BFV_charge}) and so-called coisotropic graphs is presented. More precisely Theorem \ref{thm:BFV_MC} in subsection \ref{s:(n)MC} establishes a one-to-one correspondence between equivalence classes of normalized MC-elements of a BFV-complex and coisotropic graphs. Although the BFV-complex is $L_{\infty}$ quasi-isomorphic to the strong homotopy Lie algeberoid  according to Theorem \ref{thm:big} in subsection \ref{s:main} the two structures capture different information in the non-formal regime. As a demonstration of this phenomenon we provide a simple example of a coisotropic submanifold inside a Poisson manifold where the strong homotopy Lie algeberoid fails to detect obstructions to coisotropic deformations. In the formal setting the normalization condition on MC-elements introduced in subsection \ref{s:(n)MC} turns out to be superflous.
Furthermore we use the BFV-complex to treat an example which was also considered in \cite{OhPark} and \cite{Zambon} and recover some of the results derived there. 

\subsection{Deformations of coisotropic Submanifolds}\label{s:def_coisotropicsubmanifolds}

Let $S$ be a coisotropic submanifold of a smooth finite dimensional Poisson manifold $(M,\Pi)$. We
fix an embedding of the normal bundle of $S$ into $M$. Hence we obtain a vector bundle $E \to S$ such that
$E$ is equipped with a Poisson bivector field $\Pi$ for which $S \to E$ is coisotropic.

Consider all embedded submanifolds of $E$. These form a subset $\mathcal{S}(E)$ of the set $\mathcal{P}(E)$ of all subsets of $E$.
There is a map
\begin{eqnarray*}
\widetilde{graph}: \Gamma(E) &\to& \mathcal{S}(E)\\
\mu &\mapsto& S_{\mu}:=\{(x,-\mu(x))\in E : x \in S\}.
\end{eqnarray*}
We denote the intersection of the image of $\widetilde{graph}$ with the space of all coisotropic submanifolds of
$(E,\Pi)$ by $\mathcal{C}(E,\Pi)$, the set of \textsf{coisotropic graphs}.

Given the set $\mathcal{C}(E,\Pi)$ one can ask the question whether it is representable in an algebraic way.
The precise meaning of this is the following: consider a differential graded Lie algebra $(V,d,[-,-])$.
In subsection \ref{s:infty} the set of MC-elements of $(V,d,[-,-])$ was defined to be
\begin{align*}
MC(V):=\{\beta \in V_{1}: d(\beta)+\frac{1}{2}[\beta,\beta]=0\}.
\end{align*}
One can ask whether there is a differential graded Lie algebra (more generally an $L_{\infty}$-algebra) $V$ such that
$MC(V)=\mathcal{C}(E,\Pi)$. We will show in subsection \ref{s:(n)MC} that this is the case if one chooses the differential graded Poisson algebra $(BFV(E,\Pi),D_{BFV},[-,-]_{BFV})$ and imposes a normalization condition.

We remark that a very special case of this situation occurs when one considers Lagrangian
submanifolds of symplectic manifolds. Let $(M,\Pi)$ be symplectic, i.e. $\Pi^{\#}$ is assumed
to be an isomorphism of bundles. Consequently $dim(M)$ must be $2n$ for some $n\in \mathbb{N}$.
A coisotropic submanifold $L$ of $M$ is called \textsf{Lagrangian} if $dim(L)=n$. Using an extension of Darboux's Theorem
due to Weinstein (\cite{WeinsteinD}) one can show that there is an embedding of the normal bundle $E$ of $L$ into $M$ as a tubular neighbourhood
such that 
\begin{align*}
\mathcal{C}(E,\Pi)\cong \{\gamma \in \Omega^{1}(L): d_{DR}(\gamma)=0\}.
\end{align*}

A generalization of this statement to coisotropic submanifolds $S$ of symplectic manifolds $(M,\Pi)$ was investigated in 
\cite{OhPark}. It was shown that 
\begin{align*}
\mathcal{C}^{c}(E,\Pi)\cong MC^{c}(\Gamma(\wedge E))
\end{align*}
where
$\Gamma(\bigwedge E)$ is equipped with the structure of the strong homotopy Lie algebroid associated to
$S$, see Definition \ref{def:shLA} in subsection \ref{s:shLA}. The superscript $c$ stands for ``close'' and refers to the fact that only sections sufficient close to the zero section are taken into account.

The arguments in \cite{OhPark} heavily
rely on Gotay's study of coisotropic submanifolds inside symplectic manifolds, see \cite{Gotay}. Gotay showed that the pull back of the symplectic form to the submanifold determines the symplectic form on a tubular neighbourhood (up to neighbourhood equivalence). In particulary this implies that there is an embedding of the normal bundle of a coisotropic submanifold into the symplectic manifold such that the Poisson bivector field is polynomial in fibre directions. This fails in the Poisson case.

The following example shows that the results concerning the deformation problem of coisotropic submanifolds inside symplectic manifolds mentioned above do not carry over to the Poisson case:
Consider $\mathbb{R}^{2}$ equipped with the smooth Poisson bivector field
\begin{align*}
\Pi := \left\{ 
\begin{array}{c@{\quad \mbox{for} \quad}c} 
0 & (x,y)=(0,0) \\
exp(-\frac{1}{x^{2}+y^{2}})\frac{\partial}{\partial x}\wedge\frac{\partial}{\partial y} & (x,y)\neq (0,0).
\end{array}
\right.
\end{align*}
It vanishes to all orders at $(0,0)$ but is symplectic on $\mathbb{R}^{2}\setminus\{(0,0)\}$.
The point $(0,0)$ is a coisotropic submanifold and obviously
\begin{align*}
\mathcal{C}(\mathbb{R}^{2},\Pi)=\{(0,0)\}.
\end{align*}
However the strong homotopy Lie algebroid associated to $(0,0)$ is $(\mathbb{R}^{2},0,\dots)$, so
\begin{align*}
MC(\mathbb{R}^{2})\cong \mathbb{R}^{2}.
\end{align*}
Hence $\mathcal{C}(\mathbb{R}^{2},\Pi)$ is not isomorphic to $MC(\mathbb{R}^{2})$.

\subsection{(Normalized) MC-elements and the Gauge Action}\label{s:(n)MC}

Let $E \to S$ be a finite rank vector bundle over a smooth finite dimensional manifold $S$. Assume $(E,\Pi)$ is a Poisson manifold such that $S$ is a coisotropic submanifold. The aim is to study the set of MC-elements and the deformation problem associated to the BFV-complex
$(BFV(E,\Pi),D_{BFV},[-,-]_{BFV})$, see Definition \ref{def:BFV} in subsection \ref{s:BFV_charge}.

Recall that the BFV-differential $D_{BFV}$ is given by the adjoint action of a special degree one
element $\Omega$ which was constructed in subsection \ref{s:BFV_charge}. Consequently the MC-equation for the
BFV-complex can be written as
\begin{align}\label{BFV_MC}
[\Omega + \beta, \Omega + \beta]_{BFV}=0
\end{align}
for
 $\beta \in BFV^{1}(E,\Pi)$.

\begin{Definition}\label{def:BFV_MC}
Let $E \to S$ be a finite rank vector bundle over a smooth finite dimensional manifold $S$. Assume $(E,\Pi)$ is a Poisson manifold such that $S$ is a coisotropic submanifold.

The set of \textsf{algebraic Maurer-Cartan elements} associated to the BFV-complex $(BFV(E,\Pi),D_{BFV}(-)=[\Omega,-]_{BFV},[-,-]_{BFV})$ is given by
\begin{align*}
\mathcal{D}_{alg}(E,\Pi):=\{\beta \in BFV^{1}(E,\Pi): [\Omega +\beta,\Omega+\beta]_{BFV}=0\}.
\end{align*}
\end{Definition}

We remark that $\mathcal{D}_{alg}(E,\Pi)$ contains elements that do not possess a clear geometric interpretation. Moreover $-\Omega$ is an element of $\mathcal{D}_{alg}(E,\Pi)$ that corresponds to the fact that $E$ is a coisotropic submanifold of $(E,\Pi)$. However we would prefer to study coisotropic submanifolds of $E$
that are ``similar'' to $S$ only, so they should at least be of the same dimension as $S$.

These defects can be cured with the help of a normalization condition on $\beta$. By definition
\begin{align*}
BFV^{1}(E,\Pi):=\Gamma\left(\bigoplus_{k\ge 0}(\wedge^{(k+1)} \mathcal{E}\otimes \wedge^{k}\mathcal{E}^{*})\right)
\end{align*}
where $\mathcal{E}\to E$ is the pull back bundle of $E \to S$ under $E\to S$. Hence
$\beta \in BFV^{1}(E,\Pi)$ decomposes uniquely into
\begin{align*}
\beta=\sum_{k\ge 0}\beta_{k}
\end{align*}
with $\beta_{k} \in \Gamma(\bigwedge^{(k+1)}\mathcal{E} \otimes \bigwedge^{k}\mathcal{E}^{*})=:BFV^{k+1,k}(E,\Pi)$.
In particular we obtain a map
\begin{eqnarray*}
T: BFV^{1}(E,\Pi) &\to& \Gamma(\mathcal{E})\\
\beta &\mapsto& \beta_{0}
\end{eqnarray*}
which we call the \textsf{truncation map}. Furthermore there is a natural map $p^{!}: \Gamma(E) \to \Gamma(\mathcal{E})$ given by the pull back of sections.

\begin{Definition}\label{def:BFV_MC2}
Let $E \to S$ be a finite rank vector bundle over a smooth finite dimensional manifold $S$. Assume $(E,\Pi)$ is a Poisson manifold such that $S$ is a coisotropic submanifold.

The set of \textsf{normalized Maurer-Cartan elements} associated to the BFV-complex $(BFV(E,\Pi),D_{BFV}(-)=[\Omega,-]_{BFV},[-,-]_{BFV})$ is given by
\begin{align*}
\mathcal{D}_{nor}(E,\Pi):=  \mathcal{D}_{alg}(E,\Pi) \cap T^{-1}(p^{!}(\Gamma(E)).
\end{align*}
\end{Definition}

Assume that $\beta \in \mathcal{D}_{nor}(E,\Pi)$, consequently 
\begin{align*}
T(\Omega + \beta)=\Omega_{0}+p^{!}(\mu)
\end{align*}
for a unique $\mu\in \Gamma(E)$. It is straightforward to check that the set $zero(\Omega_{0}+p^{!}(\mu))$
of zeros of the section $\Omega_{0}+p^{!}(\mu)$ is given by the submanifold $\widetilde{graph}(\mu)=:S_{\mu}$ of $E$.
In conclusion we obtain a map
\begin{eqnarray*}
\mathcal{Z}: \mathcal{D}_{nor}(E,\Pi) &\to& \mathcal{S}(E)\\
\beta &\mapsto& zero(T(\Omega+\beta))
\end{eqnarray*}
with $\mathcal{S}(E)$ denoting the set of embedded submanifolds of $E$.

We consider the adjoint action of $BFV^{0}(E,\Pi)$ on $BFV(E,\Pi)$. The Poisson algebra $BFV^{0}(E,\Pi)$
comes equipped with a filtration by Poisson subalgebras $BFV^{0}_{\ge r}(E,\Pi):=BFV^{0}(E,\Pi)\cap BFV_{\ge r}(E,\Pi)$
where $BFV_{\ge r}(E,\Pi)$ was defined as $\Gamma(\bigwedge \mathcal{E} \otimes \bigwedge^{\ge r}\mathcal{E}^{*})$.
Let $\widetilde{BFV}(E,\Pi)$ be the space of smooth sections of the pull back bundle of
$\bigwedge \mathcal{E}\otimes \bigwedge \mathcal{E}^{*}$ under $E \times [0,1] \to E$. This graded algebra
inherits the structure of a graded Poisson algebra and all the gradings (by ghost degree, ghost-momentum degree, total degree) and the filtration by $BFV_{\ge r}(E,\Pi)$ from $BFV(E,\Pi)$. In particular we obtain
a Poisson algebra $\widetilde{BFV}^{0}(E,\Pi)$ which is filtered by Poisson subalgebras $\widetilde{BFV}^{0}_{\ge r}(E,\Pi)$. It acts on $BFV(E,\Pi)$ by time-dependent endomorphisms which are derivations for both the associative algebra structure and the graded Poisson bracket $[-,-]_{BFV}$. We denote the Lie algebra of such time-dependent endomorphisms given by elements of $\widetilde{BFV}^{0}(E,\Pi)$ by $\underline{\mathfrak{inn}}(BFV(E,\Pi))$. Such endomorphisms can be interprated as time-dependent vector fields on the smooth graded manifold $\mathcal{E}[1]\oplus \mathcal{E}^{*}[-1]$ that preserve the Poisson bivector field $\hat{\Pi}$, see Corollary \ref{cor:lift} in subsection \ref{s:lifting}.

The group of \textsf{automorphisms} $Aut(BFV(E,\Pi))$ is the space of all isomorphisms of the unital graded commutative
associative algebra $BFV(E,\Pi)$ that preserve the total degree and the graded Poisson bracket $[-,-]_{BFV}$. An automorphism
$\psi$ is called \textsf{inner} if it is generated by an element of $\underline{\mathfrak{inn}}(BFV(E,\Pi))$. More precisely we impose that 
\begin{itemize}
\item there is a family of automorphisms $\{\psi_{t}\}_{t\in [0,1]}$ with $\psi_{0}=id$ and
$\psi_{1}=\psi$,
\item there is a morphism of unital graded commutative associative algebras and Poisson algebras
$\hat{\psi}: BFV(E,\Pi) \to \widetilde{BFV}(E,\Pi)$
\end{itemize}
such that
\begin{itemize}
\item the composition of $\hat{\psi}$ with the pull back along the inclusion $E \times \{t\} \to E \times [0,1]$ coincides with $\psi_{t}$,
\item the time-dependent derivation of $BFV(E,\Pi)$ that maps $\beta$ to
\begin{align*}
(e,s) \mapsto \frac{d}{dt}|_{t=s}\left(\hat{\psi}(\beta)|_{e}\right)
\end{align*}
is an element of $\underline{\mathfrak{inn}}(BFV(E,\Pi))$.
\end{itemize}
We denote the subset of inner automorphisms of $BFV(E,\Pi)$ by $Inn(BFV(E,\Pi))$ which can be checked to be
a subgroup of $Aut(BFV(E,\Pi))$. Moreover the filtration of $\widetilde{BFV}^{0}(E,\Pi)$ by the Poisson subalgebras
$\widetilde{BFV}^{0}_{\ge r}(E,\Pi)$ yields a filtration of $Inn(BFV(E,\Pi))$ by subgroups which we denote by
$Inn^{\ge r}(BFV(E,\Pi))$.

The group $Aut(BFV(E,\Pi))$ acts on $\mathcal{D}_{alg}(E,\Pi)$ via
\begin{align*}
(\hat{\Theta},\alpha) \mapsto \hat{\Theta}(\Omega+\alpha) - \Omega
\end{align*} 
and consequently so do all the groups
$Inn^{\ge r}(BFV(E,\Pi))$. Observe that the action of $Inn^{\ge 2}(BFV(E,\Pi))$ on $\mathcal{D}_{alg}(E,\Pi)$ restricts to an action on $\mathcal{D}_{nor}(E,\Pi)$.

\begin{Theorem}\label{thm:BFV_MC}
Let $E \to S$ be a finite rank vector bundle over a smooth finite dimensional manifold $S$. Assume $(E,\Pi)$ is a Poisson manifold such that $S$ is a coisotropic submanifold.

Mapping elements of $BFV^{1}(E,\Pi)$ to the zero set of their truncation induces a bijection between
\begin{enumerate}
\item $\mathcal{D}_{nor}(E,\Pi)/Inn^{\ge 2}(BFV(E,\Pi))$ and
\item $\mathcal{C}(E,\Pi)$.
\end{enumerate}
\end{Theorem}

\begin{proof}
\ \linebreak
\underline{claim A:} An element $p^{!}(\mu) \in \Gamma(\mathcal{E})$ can be extented to a MC-element
of $(BFV(E,\Pi),D_{BFV},[-,-]_{BFV})$ if and only if 
\begin{align*}
S_{\mu}:=\{(x,-\mu(x))\in E| x \in S\}
\end{align*}
is a coisotropic submanifold of $(E,\Pi)$.
\newline 

Given an arbitrary $\mu \in \Gamma(E)$ we want to construct a $\beta \in \mathcal{D}_{alg}(E,\Pi)$ decomposing as
\begin{align*}
\beta=\sum_{k \ge 0} \beta_{k}
\end{align*}
where $\beta_{k} \in \Gamma(\bigwedge^{(k+1)}\mathcal{E}\otimes \bigwedge^{k} \mathcal{E}^{*})$ such that
$\beta_{0}=p^{!}(\mu)$ holds. This is a generalization of the construction of $\Omega$ given in the proof of Proposition \ref{prop:charge} in subsection \ref{s:BFV_charge}.

First consider $\Omega_{0}+p^{!}(\mu) \in \Gamma(\mathcal{E})$. Obviously 
\begin{align*}
[\Omega_{0}+p^{!}(\mu),\Omega_{0}+p^{!}(\mu)]_{G}=0
\end{align*}
holds. Recall that $G$ is the Poisson bivector field on $\mathcal{E}^{*}[1]\oplus \mathcal{E}[-1]$ that
corresponds to the fibre pairing between $\mathcal{E}$ and $\mathcal{E}^{*}$.
Consequently we obtain a differential
\begin{align*}
\delta[\mu](-):=[\Omega_{0}+p^{!}(\mu),-]_{G}=: \delta(-) + \partial_{\mu}(-).
\end{align*}
In the proof of Proposition \ref{prop:charge} in subsection \ref{s:BFV_charge} a homotopy $h$ for $\delta$ was defined satisfying
$h \circ \delta + \delta \circ h = id - p^{*}\circ i^{*}$ where
$p^{*}: \Gamma(\bigwedge E) \to BFV(E,\Pi)$ is essentially
given by the pull back $p^{!}$ and $i^{*}: BFV(E,\Pi) \to \Gamma(\bigwedge E)$ is given by natural restriction and projection maps. The concrete formula of $h$ implies that
\begin{align*}
h \circ \partial_{\mu} + \partial_{\mu} \circ h = 0
\end{align*}
since $p^{!}(\mu)$ is a section of $\Gamma(\mathcal{E}) \subset BFV(E,\Pi)$ that is constant along the fibres of $E \to S$ and consequently
\begin{align*}
h \circ \delta[\mu] + \delta[\mu] \circ h = id - p^{*}\circ i^{*}.
\end{align*}
Observe that the maps $p^{*}$ and $i^{*}$ are no longer morphisms of complexes with respect to $\delta_{\mu}$.

Consider the diffeomorphism $q_{\mu}: S_{\mu}:=\{(x,-\mu(x))|x\in S\} \to S$ and the pull back vector bundle
$q_{\mu}^{!}(E) \to S_{\mu}$.

\underline{claim A.1:} $H^{\bullet}(BFV(E,\Pi),\delta[\mu]) \cong \Gamma(\bigwedge^{\bullet} (q_{\mu}^{!}E))$
\newline 
Since $S_{\mu}$ and $S$ are diffeomorphic there is a vector bundle isomorphism between
$q_{\mu}^{!}(E)$ and $E$ which induces and isomorphism $\vartheta$ of graded algebras between $\Gamma(q_{\mu}^{!}(E))$ and $\Gamma(E)$. It is straightforward to check that 
\begin{align*}
p_{\mu}^{*}: \Gamma(q_{\mu}^{!}(E)) \xrightarrow{\vartheta} \Gamma(E) \xrightarrow{p^{*}} BFV(E,\Pi)
\end{align*}
and
\begin{align*}
i_{\mu}^{*}: BFV(E,\Pi) \xrightarrow{i^{*}} \Gamma(E) \xrightarrow{\vartheta^{-1}} \Gamma(q_{\mu}^{!}(E))
\end{align*}
are chain maps between $(BFV(E,\Pi),\delta_{\mu})$ and $(\Gamma(p_{\mu}^{!}E),0)$. In fact, $p_{\mu}^{*}$ is given by the unique extension of a section of $q	_{\mu}^{!}(E)$ to a section of $\bigwedge \mathcal{E} \otimes \bigwedge \mathcal{E}^{*}$ that is constant along the fibres of $\mathcal{E} \to E$. Futhermore $i_{\mu}^{*}$ is given by the composition of $BFV(E,\Pi) \to \Gamma(\bigwedge \mathcal{E})$ with
the evaluation at $S_{\mu}$. Obviously $i_{\mu}^{*} \circ p_{\mu}^{*} = id$ and
\begin{align*}
h \circ \delta[\mu] + \delta[\mu] \circ h = id - p_{\mu}^{*}\circ i_{\mu}^{*}
\end{align*}
hold. This implies the claim A.1.

Having established claim A.1 the constructions of elements $\gamma_{1}, \gamma_{2}, \dots$ with
$\gamma_{k} \in \Gamma(\bigwedge^{(k+1)}\mathcal{E}\otimes \bigwedge^{k}\mathcal{E}^{*})$ such that
$\Omega_{0} + p^{!}(\mu) + \gamma_{1} + \gamma_{2} + \cdots $ is a MC-elements goes through 
as in the proof of Proposition \ref{prop:charge} in subsection \ref{s:BFV_charge}: One tries to extend $\Omega_{0}+p^{!}(\mu)$ inductively and meets obstructions classes at each level. The first obstruction class vanishes if and only if
$S_{\mu}$ is a coisotropic submanifold of $E$: $2R_{0}:=[\Omega_{0}+p^{!}(\mu),\Omega_{0}+p^{!}(\mu)]_{\iota_{\nabla}(\Pi)}$ gives a cohomology class in $H(BFV(E,\Pi),\delta_{\mu})$, the evaluation of $2R_{0}$
at $S_{\mu}$ is $0$ if and only if the vanishing ideal of $S_{\mu}$ is a Lie subalgebra under the Poisson bracket $[-,-]_{\Pi}$. This is equivalent to $S_{\mu}$ being coisotropic. When the class $[R_{0}]$ is zero, we can find $\gamma_{1}$ with $R_{0}=-\delta_{\mu}(\gamma_{1})$ which will be our first correction term.
All higher obstruction classes vanish due to claim A.1.
Then setting $\beta_{0}:=p^{!}(\mu)$ and $\beta_{m}:=\gamma_{m}-\Omega_{m}$ for $m>1$ yields a MC-element
\begin{align*}
\beta:=\sum_{k \ge 0}\beta_{k}
\end{align*}
of the desired form.

\underline{claim B:} Given two elements $\alpha$ and $\beta$ of $\mathcal{D}_{alg}(E,\Pi)$ with
$T(\alpha)=T(\beta)=p^{!}(\mu)$ for some $\mu \in \Gamma(E)$, there is an element of
$Inn^{\ge 2}(BFV(E,\Pi))$ mapping $\alpha$ to $\beta$.
\newline

Observe that inner derivations given by the adjoint action of $\widetilde{BFV}^{0}_{\ge 2}(E,\Pi)$
are nilpotent and therefore always integrate to an inner automorphism.
Assume that $\beta$ and $\alpha$ coincide up to order $k>0$, i.e. 
\begin{align*}
\beta - \alpha = 0 \mbox{ mod } BFV_{\ge k}(E,\Pi).
\end{align*}
The MC-equation for $\beta$ and $\alpha$ implies that
\begin{align*}
\delta[\mu](\beta_{k})=F(\beta_{0},\dots,\beta_{(k-1)})=F(\alpha_{0},\dots,\alpha_{(k-1)})=\delta[\mu](\alpha_{k}).
\end{align*}
Here $F$ is a function that can be constructed from the MC-equation: the equation $1/2[\Omega+\beta,\Omega+\beta]_{BFV}=0$
can be decomposed with respect to the ghost-momentum degree. For the ghost-momentum degree $k-1$ one obtains
$\delta[\mu](\beta_k)$ plus other terms depending on $(\beta_0,\dots,\beta_{(k-1)})$ only. We denote the sum of this
other terms by $-F$.

Consequently $\delta[\mu](\beta_{k}-\alpha_{k})=0$. By claim A.1 and the assumption $k>0$ there is an element
$\varepsilon_{k} \in BFV^{(k+1,k+1)}(E,\Pi)$ such that $\beta_{k}-\alpha_{k}=\delta[\mu](\epsilon_{k})$.
Then 
\begin{eqnarray*}
exp(-[\epsilon_{k},-]_{BFV})(\alpha)&=&\alpha -[\epsilon_{k},\alpha]_{BFV} \mbox{ mod } BFV_{\ge (k+1)}(E,\Pi)\\
&=& \alpha + [\alpha,\epsilon_{k}]_{BFV} \mbox{ mod } BFV_{\ge (k+1)}(E,\Pi)\\
&=& \alpha + \delta[\mu](\epsilon_{k}) \mbox{ mod } BFV_{\ge (k+1)}(E,\Pi)\\
&=& \beta \mbox{ mod } BFV_{\ge (k+1)}(E,\Pi)
\end{eqnarray*}
so $exp(-[\epsilon_{k},-])(\alpha)$ and $\beta$ coincide up to order $k+1$.

Inductively one findes $\varepsilon_{1},\epsilon_{2},\dots,\epsilon_{N}$ such that
\begin{align*}
exp(-[\epsilon_{N},-])\cdots exp(-[\epsilon_{2},-])exp(-[\epsilon_{1},-])(\alpha)=\beta.
\end{align*}
Then the BCH-formula yields an $\varepsilon\in BFV^{0}_{\ge 2}(E,\Pi)$ such that
the inner automorphism generated by its adjoint action on $BFV(E,\Pi)$ maps $\alpha$ to $\beta$.
\end{proof}

\subsection{An Example}\label{s:example_Marco}

We consider an example that was first presented in \cite{Zambon} and that was also investigated in \cite{OhPark}.
Zambon showed that the space of coisotropic deformations $\mathcal{C}(E,\Pi)$ ``near'' a fixed coisotropic submanifold $S$
cannot carry the structure
of a (Fr\'echet) manifold because there exist ``tangent vectors'' whose sum is not tangent to $\mathcal{C}(E,\Pi)$. Oh and
Park showed that this can be understood with the help of the strong homotopy Lie algebroid $(\Gamma(\bigwedge E),\partial_{s}=\mu^{1},\mu^{2},\dots)$ associated to $S$, see Definition \ref{def:shLA} in subsection \ref{s:shLA}. The extension of elements in the first Lie algebroid cohomology to MC-elements meets several obstructions,
and the first of them is given by a quadratic relation. Hence, the sum of solutions might fail to be a solution again which
explains Zambon's observation.

Consider the vector bundle $E=\mathbb{R}^{2} \times (S^{1})^{4} \to (S^{1})^{4}$ with coordinates
$(x^{1},x^{2},\theta^{1},\theta^{2},\theta^{3},\theta^{4})$ ($\theta$ denotes the angle-coordinate on $S^{1}$). We
equip $E$ with the symplectic form
\begin{align*}
\omega = d\theta^{1} \wedge dx^{1} + d\theta^{2} \wedge dx^{2} + d\theta^{3} \wedge d \theta^{4}
\end{align*}
and define $S:=(S^{1})^{4}$ which is a coisotropic submanifold of $E$. 

The BFV-complex $BFV(E,\omega^{-1})$ is given by the smooth functions on the smooth graded manifold $E \times (\mathbb{R}^{*}[1])^{2} \times (\mathbb{R}[-1])^{2} \to E$. We introduce fibre coordinates $(c_{1},c_{2})$ on $(\mathbb{R}^{*}[1])^{2}$ and $(b^{1},b^{2})$ on $(\mathbb{R}[-1])^{2}$. Since the bundle $E$ is flat we can just set
\begin{align*}
[-,-]_{BFV}=:[-,-]_{G}+[-,-]_{\omega}
\end{align*}
where $[-,-]_{G}$ denotes the graded Poisson bracket given by the pairing between $(\mathbb{R}^{*}[1])^{2}$
and $(\mathbb{R}[-1])^{2}$ and $[-,-]_{\omega}$ is the Poisson bracket associated to the symplectic form $\omega$.

The element $\Omega_{0}$ reads $c_{1}x^{1}+c_{2}x^{2}$ and it is closed with respect
to the graded Poisson bracket on the BFV-complex, so no extension is needed. The BFV-differential $D_{BFV}$ of the
BFV-complex is given by
\begin{align*}
D=
x^{1} \frac{\partial}{\partial b^{1}} + x^{2} \frac{\partial}{\partial b^{2}}+
c^{1} \frac{\partial}{\partial \theta^{1}} + c^{2} \frac{\partial}{\partial \theta^{2}} .
\end{align*}
It is straightforward to check that
the cohomology with respect to $D_{BFV}$ is given by periodic functions in the variables $\theta^{3}$ and $\theta^{4}$ tensored
by the Grassmann-algebra generated by $c^{1}$ and $c^{2}$.

The MC-equation reads
\begin{align*}
[\Omega_{0}+\beta,\Omega_{0}+\beta]_{BFV}=0
\end{align*}
and if we assume that $\beta$ is a $D_{BFV}$-cocycle it reduces to
\begin{align*}
[\beta,\beta]_{\omega}=0.
\end{align*}
If we impose that $\beta$ is a normalized MC-element (see subsection \ref{s:(n)MC}) it only depends on the variables
$\theta^{1},\theta^{2},\theta^{3}$ and $\theta^{4}$. In this case the MC-equation reduces further to
\begin{align}\label{condition}
\{\beta,\beta\}_{S}=0.
\end{align}
where
\begin{align*}
\{f,g\}_{S}:= \frac{\partial f}{\partial \theta^{4}}\frac{\partial g}{\partial \theta^{3}}
- \frac{\partial f}{\partial \theta^{3}}\frac{\partial g}{\partial \theta^{4}}.
\end{align*}
Condition \eqref{condition} was also found in \cite{OhPark}.

Consider an element $c_{1}f^{1}+c_{2}f^{2}$ where $f^{1}$ and $f^{2}$ depend on the angle-variables only. When does this section define
a coisotropic submanifold? In the proof of Proposition \ref{thm:BFV_MC} in subsection \ref{s:(n)MC} we showed that this is equivalent to
\begin{align}\label{h2}
[\Omega_{0}+c_{1}f^{1}+c_{2}f^{2},\Omega_{0}+c_{1}f^{1}+c_{2}f^{2}]_{\omega}
\end{align}
being exact with respect to $\delta[c_{1}f^{1}+c_{2}f^{2}]:=(x^{1}+f^{1})\frac{\partial}{\partial b^{1}} + (x^{2}+f^{2})\frac{\partial}{\partial b^{2}}$. 
Computing the bracket \eqref{h2} yields
\begin{align*}
2c_{1}c_{2}(\frac{\partial f^{1}}{\partial \theta^{2}} - \frac{\partial f^{2}}{\partial \theta^{1}}+
\frac{\partial f^{1}}{\partial \theta^{4}}\frac{\partial f^{2}}{\partial \theta^{3}}-
\frac{\partial f^{2}}{\partial \theta^{4}}\frac{\partial f^{1}}{\partial \theta^{3}}).
\end{align*}
We denote this expression by $H$. It only depends on the angle-variables. Exactness of $H$ translates
into the condition that there exists a pair of functions $g_{1}$ and $g_{2}$ that might depend on
all variables on $E$ such that
\begin{align*}
\delta[c_{1}f^{1}+c_{2}f^{2}](b^{1}g_{1}+b^{2}g_{2})=(x^{1}+f^{1})g_{1}+(x^{2}+f^{2})g_{2}=H.
\end{align*}
Since $H$ is constant in $x^{1}$ and $x^{2}$ the left hand side $(x^{1}+f^{1})g_{1}+(x^{2}+f^{2})g_{2}$ is too. Evaluating it at $x^{1}=-f^{1}$ and $x^{2}=-f^{2}$ implies that $H$ must vanish identically.
Hence a section of the bundle $\bigwedge(\mathbb{R}^{2})\times E \to E$ given by $c_{1}f^{1}+c_{2}f^{2}$ defines a coisotropic submanifold iff
\begin{align*}
\frac{\partial f^{1}}{\partial \theta^{2}} - \frac{\partial f^{2}}{\partial \theta^{1}}+
\frac{\partial f^{1}}{\partial \theta^{4}}\frac{\partial f^{2}}{\partial \theta^{3}}-
\frac{\partial f^{2}}{\partial \theta^{4}}\frac{\partial f^{1}}{\partial \theta^{3}}=0.
\end{align*}
Up to different sign convenctions this condition coincides with the one given in \cite{Zambon}, where it was derived in an analytical context.

\subsection{Formal Deformations of coisotropic Submanifolds}\label{s:fordef_coisotropicsubmanifolds}
Let $E\to S$ be a finite rank vector bundle over a smooth finite dimensional manifold $S$. Assume
$(E,\Pi)$ is a Poisson manifold such that $S$ is coisotropic.

We introduce a formal parameter $\varepsilon$ of degree $0$ and consider the graded commutative algebra $BFV(E,\Pi)[[\varepsilon]]$. It inherits the structure of a differential graded Poisson algebra $(BFV(E,\Pi)[[\varepsilon]],D_{BFV},[-,-]_{BFV})$ from $BFV(E,\Pi)$, see Definition \ref{def:BFV} in subsection \ref{s:BFV_charge}. We define the space of \textsf{formal MC-elements} by
\begin{align*}
\mathcal{D}_{for}(E,\Pi):=\{\beta \in \varepsilon BFV(E,\Pi)[[\varepsilon]]: [\Omega + \beta,\Omega + \beta]_{BFV}=0\}.
\end{align*}
Recall that $\Omega$ is a degree $1$ element of $BFV(E,\Pi)$ such that $[\Omega,\Omega]_{BFV}=0$ and if one decomposes
$\Omega$ with respect to the ghost-momentum degree, i.e.
\begin{align*}
\Omega=\sum_{k \ge 0}\Omega_{k}
\end{align*}
with $\Omega_{k} \in \Gamma(\bigwedge^{(k+1)}\mathcal{E}\otimes \bigwedge^{k}\mathcal{E}^{*})$, $\Omega_{0}$ is required to be the tautological section of $\mathcal{E} \to E$.

In subsection \ref{s:(n)MC} we introduced $\widetilde{BFV}^{0}(E,\Pi)$ and its action by derivations on $BFV(E,\Pi)$.
In the formal setting one considers $\varepsilon \widetilde{BFV}^{0}(E,\Pi)[[\varepsilon]]$ and its action on $BFV(E,\Pi)[[\varepsilon]]$. Since the action by such a derivation is pro-nilpotent, it always integrates to an automorphism of $(BFV(E,\Pi)[[\varepsilon]],[-,-]_{BFV})$. We denote the subgroup of these automorphisms by $Inn_{for}(BFV(E,\Pi))$. As explained in subsection \ref{s:(n)MC} this group naturally acts on $\mathcal{D}_{for}(E,\Pi)$ by
\begin{eqnarray*}
Inn_{for}(BFV(E,\Pi)) \times \mathcal{D}_{for}(E,\Pi) &\to& \mathcal{D}_{for}(E,\Pi)\\
(\Psi, \beta) &\mapsto& \Psi(\Omega+\beta) - \Omega.
\end{eqnarray*}

Throughout subsection \ref{s:(n)MC} we had to fix a normalization condition on the MC-elements of
$(BFV(E,\Pi),D_{BFV},[-,-]_{BFV})$ in order to make connection to the geometry
of coisotropic submanifolds of $(E,\Pi)$. We considered the truncation map $T:BFV^{1}(E,\Pi) \to \Gamma(\mathcal{E})$
and imposed that the image of a MC-element $\beta$ under $T$ has to lie in the image of the pull back map $\Gamma(E) \to \Gamma(\mathcal{E})$.
In the formal setting no normalization condition is needed due to
the following

\begin{Lemma}\label{l:noformalnormalization}
For every $\beta \in \mathcal{D}_{for}(E,\Pi)$ there is a $\Psi \in Inn_{for}(BFV(E,\Pi))$ such that
the image of $\Psi(\beta)$ under $T: BFV^{1}(S,\Pi)[[\varepsilon]] \to \Gamma(\mathcal{E})[[\varepsilon]]$
is given by a pull back from $\varepsilon\Gamma(E)[[\varepsilon]]$.
\end{Lemma}

\begin{proof}
The element $\beta \in \mathcal{D}_{for}(E,\Pi) \subset \varepsilon BFV^{1}(E,\Pi)[[\varepsilon]]$ decomposes
uniquely into
\begin{align*}
\beta = \sum_{k\ge 0}\beta_{k}
\end{align*}
with $\beta_{k} \in \varepsilon \Gamma(\bigwedge^{(k+1)}\mathcal{E}\otimes \bigwedge^{k}\mathcal{E}^{*})[[\varepsilon]]$.
In particular $\beta_{0} \in \varepsilon \Gamma(\mathcal{E})[[\varepsilon]]$ which we further decompose as
\begin{align*}
\beta_{0}=\sum_{l\ge 1}\beta_{0}(l)\varepsilon^{l}.
\end{align*}

Consider the cocycle $[\beta_{0}(1)] \in H(BFV(E,\Pi),\delta)$. Using the homotopy $h$ introduced in the proof
of Proposition \ref{prop:charge} in subsection \ref{s:BFV_charge} one finds a section $\tilde{\beta}_{0}(1) \in \varepsilon\Gamma(\mathcal{E})$ that is a pull back from a section of $\varepsilon \Gamma(E)$ such that
$[\beta_{0}(1)]=[\tilde{\beta}_{0}(1)]$. Hence there is $\gamma(1) \in \varepsilon \Gamma(\mathcal{E}\otimes \mathcal{E}^{*})$ satisfying $\beta_{0}(1) = \tilde{\beta}_{0}(1) + \delta(\gamma(1))$.
The automorphism $exp([\gamma(1),-]_{BFV})$ maps the MC-element $\Omega + \beta$ to another one whose image under the truncation map modulo $\varepsilon^{2}$ is given by
\begin{align*}
\Omega_{0}+\beta_{0}(1) - \delta(\gamma(1))=\Omega_{0}+\tilde{\beta}_{0}(1),
\end{align*}
i.e. the new MC-element has the desired property modulo $\varepsilon^{2}$.

Let us assume that we established $\beta_{0}=p^{!}(\mu)$ modulo $\varepsilon^{k}$ for some $\mu \in \varepsilon \Gamma(E)[[\varepsilon]]$. Consider the $\delta$-cocycle $\beta_{0}(k)$. As before there is $\gamma(k)\in \varepsilon^{k}\Gamma(\mathcal{E}\otimes \mathcal{E}^{*})$ and a pull back section $\tilde{\beta}_{0}(k) \in \varepsilon^{k}\Gamma(\mathcal{E})$ such that
\begin{align*}
\beta_{0}(k)=\tilde{\beta}_{0}(k) + \delta(\gamma(k))
\end{align*}
holds. We consider the inner automorphism $exp([\gamma(k),-]_{BFV})$ which maps the MC-element $\Omega+\beta$ to an otherone whose truncation modulo $\varepsilon^{(k+1)}$ is given by
\begin{align*}
\Omega_{0}+\sum_{1\le m\le k}\beta_{0}(m) - \delta(\gamma(k)) = \Omega_{0} + \sum_{1 \le m \le (k-1)} \beta_{0}(m) + \tilde{\beta}_{0}(k).
\end{align*}
Using induction with respect to the polynomial degree in $\varepsilon$, the fact that the formal variable ring is complete with respect to the $\varepsilon$-adic topology and the BCH-formalua one finds an appropriate formal inner automorphism $\Psi$.
\end{proof}

In subsection \ref{s:Poisson} we stated that one possible characterization of coisotropic submanifolds uses vanishing ideals: a submanifold $S$ of a Poisson manifold $(E,\Pi)$ is coisotropic if and only if its vanishing ideal $\mathcal{I}_{S}:=\{f \in \mathcal{C}^{\infty}(E): f|_{C}=0\}$ is a Lie subalgebra of the Poisson algebra of functions. A multiplicative ideal of a Poisson algebra that in addition is a Lie subalgebra is called a \textsf{coisotrope}, see \cite{Weinstein}.

\begin{Definition}\label{def:formaldeformations}
Let $E\to S$ be a finite rank vector bundle over a smooth finite dimensional manifold $S$. Assume
$(E,\Pi)$ is a Poisson manifold such that $S$ is coisotropic.

A 
\textsf{formal deformation} of $S$ is a coisotrope $\mathcal{I}$ of $(\mathcal{C}^{\infty}(E)[[\varepsilon]],[-,-]_{\Pi})$ such that $\mathcal{I} \mbox{ mod } \varepsilon = \mathcal{I}_{S}$.
We denote the set of formal deformations of $S$ by $\mathcal{C}_{for}(E,\Pi)$.
\end{Definition}

\begin{Lemma}\label{l:coisotropes}
Let $E\to S$ be a finite rank vector bundle over a smooth finite dimensional manifold $S$. Assume
$(E,\Pi)$ is a Poisson manifold such that $S$ is coisotropic.

There is an $Inn_{for}^{\ge 1}(BFV(E,\Pi))$-equivariant map from $\mathcal{D}_{for}(E,\Pi)$ to $\mathcal{C}_{for}(E,\Pi)$ (with trivial action on $\mathcal{C}_{for}(E,\Pi)$)
such that the image of $0$ is the $\mathbb{R}[[\varepsilon]]$-linear extension of the vanishing ideal $\mathcal{I}_{S}$ of $S$.
\end{Lemma}

\begin{proof}
Given $\beta \in \mathcal{D}_{for}(E,\Pi)$ we want to construct a coisotrope $\mathcal{I}(\beta)$ of $(\mathcal{C}^{\infty}(E)[[\varepsilon]],[-,-]_{\Pi})$ in a way that is invariant under the action of the group $Inn_{for}^{\ge 1}(BFV(E,\Pi))$ on $\mathcal{D}_{for}(E,\Pi)$.

Consider the truncation $\beta_{0} \in \varepsilon\Gamma(\mathcal{E})[[\varepsilon]]$ of $\beta$.
We choose a trivializing atlas $\{U_{\alpha}\}_{\alpha \in A}$ for the vector bundle $E \to S$ which yields a trivializing atlas for the vector bundle $\mathcal{E}\to E$. On each chart $U_{\alpha}$ we pick a local frame $(c^{\alpha}_{j})_{j=1,\dots,e}$ for the bundle $\mathcal{E}$ and obtain a unique decomposition
\begin{align*}
(\Omega_{0}+\beta_{0})|_{U_{\alpha}}=\sum_{j=1}^{e}h_{\alpha}^{j}c^{\alpha}_{j}
\end{align*}
with $h_{\alpha}^{j} \in \mathcal{C}^{\infty}(U_{\alpha}\times \mathbb{R}^{e})[[\varepsilon]]$ for $j=1,\dots, e$. Let
$J_{\alpha}$ be the multiplicative ideal of $\mathcal{C}^{\infty}(U_{\alpha}\times \mathbb{R}^{e})[[\varepsilon]]$ generated by $(h_{\alpha}^{j})_{j=1,\dots, e}$. It is straightforward to conclude from $\beta \in \mathcal{D}_{for}(E,\Pi)$ that $J_{\alpha}$ is a coisotrope
of the Poisson algebra $(\mathcal{C}^{\infty}(U_{\alpha}\times \mathbb{R}^{e})[[\varepsilon]],[-,-]_{\Pi}|_{U_{\alpha}\times \mathbb{R}^{e}})$.

Observe that the family of ideals $(J_{\alpha})_{\alpha \in A}$ can be glued together, i.e. given $U_{\alpha} \cap U_{\beta} =:U_{\alpha\beta} \neq \emptyset$ then $f\in \mathcal{C}^{\infty}(U_{\alpha\beta}\times \mathbb{R}^{e})$ lies in the restriction of $J_{\alpha}$ to $U_{\alpha\beta}\times \mathbb{R}^{e}$ if and only if it lies in the restriction of $J_{\beta}$ to $U_{\alpha\beta}\times \mathbb{R}^{e}$. This steams from the fact that the transition matrices $U_{\alpha\beta}\times \mathbb{R}^{e} \to GL(\mathbb{R}^{e})$ for the vector bundle $\mathcal{E}$ are invertible.

We define $\mathcal{I}(\beta)$ to be the set of elements of $\mathcal{C}^{\infty}(E)[[\varepsilon]]$ whose restriction
to every coordinate domain $U_{\alpha}\times \mathbb{R}^{e}$ lies in $J_{\alpha}$. An argument similar to the gluing statement above shows that $\mathcal{I}(\beta)$ is in fact independent of the choice of atlas and one easily checks that it is a coisotrope and that $\mathcal{I}(\beta)\mbox{ mod } \varepsilon = \mathcal{I}_{S}$ holds.

Furthermore $\mathcal{I}(\beta)$ is not affected if we let a bundle automorphism act on the section $\beta_{0}$. Notice that
the action of $\varepsilon BFV_{\ge 1}(E,\Pi)[[\varepsilon]]$ on $BFV(E,\Pi)[[\varepsilon]]$ induces an action on $BFV^{(1,0)}(E,\Pi)[[\varepsilon]]=\Gamma(\mathcal{E})[[\varepsilon]]$ which coincides with the action given by
\begin{eqnarray*}
BFV^{(1,1)}(E,\Pi) \to \Gamma(\mathcal{E}\otimes \mathcal{E}^{*}) \cong \Gamma(End(\mathcal{E})) \xrightarrow{exp} \Gamma(GL(\mathcal{E})).
\end{eqnarray*}
From this the $Inn_{for}^{\ge 1}(BFV(E,\Pi))$-equivariance of $\beta \mapsto \mathcal{I}(\beta)$ follows.
\end{proof}

%% file: Appendix.tex
\section{Details on the Homotopy Transfer}\label{AppendixA}

This appendix provides background information on the material presented in subsection \ref{s:homotopytransfer}. The aim is to prove Theorem \ref{thm:transfer} which is a central technical tool in sections \ref{s:BFV} and \ref{s:Liealgebroid}. We first relate the homotopy transfer to integration over isotropic subspace in the BV-Formalism. Then we check that the formulae given in \ref{s:homotopytransfer} actually work. All the results are well-known to the experts and we do not claim any originality related to this treatment.

\subsection{Connection to the BV-Formalism}\label{s:BV-formalism}

We present a heuristic derivation of the formulae for the homotopy transfer as presented in subsection \ref{s:homotopytransfer}. It makes use of the BV-Formalism which was introduced by Batalin and Vilkovisky. In the case of $A_{\infty}$-algebras a similar treament can be found in \cite{Kajiura}.

In the finite dimensional setting the BV-Formalism was made rigorous by Schwarz, see \cite{Schwarz}. Although not justified at mathematical level of rigour in the infinite-dimensional setting in general, the BV-Formalism serves as a way to obtain formulae for the homotopy transfer which can be checked to work using purly algebraic manipulations a posteriori. We remark that there are certain (infinite-dimensional) situations where a mathematical treatment can be provided -- see \cite{Costello} for instance.

Let $V$ be a graded vector space. In subsection \ref{s:infty} the one-to-one correspondence between $L_{\infty}$-algebra structures on $V$ and codifferentials of $S(V[1])$ was explained. If one assumes that $V$ is finite dimensional, the space of coderivations of the coalgebra $S(V[1])$ is in bijection to the space of derivations of the algebra $S(V^{*}[-1])=:\mathcal{C}^{\infty}(V[1])$, i.e. vector fields on $V[1]$. Under this bijection codifferentials correpsond to so called \textsf{cohomological vector fields}, i.e. derivations of degree $1$ that square to zero. Hence there is a one-to-one correspondence between $L_{\infty}$-algebra strutures on $V$ and homological vector fields on $V[1]$. Moreover flat $L_{\infty}$-algebras are encoded in homological vector fields that vanish at $0 \in V[1]$.

The space of multivector fields on $V[1]$ can be described as the space of functions on the smooth graded manifold $T^{*}[1](V[1])=V[1]\oplus V^{*}[0]$. Being a shifted cotangent bundle, this smooth graded manifold carries a graded symplectic structure. Equivalently the graded commutative algebra $\mathcal{C}^{\infty}(T^{*}[1](V[1]))$ carries the structure of a graded Poisson bracket $[-,-]_{BV}$ of degree $-1$ called the \textsf{BV-bracket}.  The space of vector fields forms a Poisson subalgebra and a vector field $X$ on $V[1]$ is cohomological if and only if $[X,X]_{BV}=0$. This equation is called the \textsf{classical master equation}.

There is a bijection between the space of homomorphisms $Hom(V[1],V)$ of $V[1]$ of degree $-1$ and the graded vector space $V^{*}[-1]\otimes V[0]$. Using a basis $\{\gamma_{i}\}$ of $V[0]$ and the dual basis $\{\gamma^{i}\}$ of $V^{*}[-1]$ the identity $id \in End(V)$ yields an element $\gamma^{i}\otimes \gamma_{i}$. One defines the following operator of degree $-1$ on $\mathcal{C}^{\infty}(V[1]\oplus V^{*}[0])$
\begin{align*}
\Delta:=\frac{\partial^{2}}{\partial \gamma^{i} \partial \gamma_{i}}
\end{align*}
which is called the \textsf{BV-operator}. It is straightforward to check that $\Delta \circ \Delta =0$. However $\Delta$ is not a derivation with respect to the graded commutative associative product of $\mathcal{C}^{\infty}(V[1]\oplus V^{*}[0])$. The deviation to being a derivation is measured by the BV-bracket $[-,-]_{BV}$, i.e.
\begin{align*}
\Delta(X\cdot Y) - \Delta(X)\cdot Y - (-1)^{|X|}X\cdot \Delta(Y) = (-1)^{|X|}[X,Y]_{BV}
\end{align*}
for homogeneous $X$ and arbitrary $Y$ in $\mathcal{C}^{\infty}(V[1]\oplus V^{*}[0])$. The quadruple
$(\mathcal{C}^{\infty}(V[1]\oplus V^{*}[0]),\cdot,\Delta,[-,-]_{BV})$ is an example of a \textsf{BV-algebra}. Given such an algebra one can write down the \textsf{quantum master equation}:
\begin{align*}
\Delta(X)+\frac{1}{2}[X,X]_{BV}=0.
\end{align*}
The importance of this equation is due to the identity
\begin{align*}
\Delta(e^{X})=(\Delta(X)+\frac{1}{2}[X,X]_{BV})e^{X}.
\end{align*}
Hence $e^{X}$ is $\Delta$-closed if and only if $X$ satisfies the quantum master equation.

Let $X$ be a cohomological vector field on a graded vector space $V[1]$ which vanishes at $0$, i.e. $V[1]$ is equipped with the structure of a flat $L_{\infty}[1]$-algebra. Denote the differential of this $L_{\infty}[1]$-structure by $\delta$ and the corresponding cohomology by $H[1]$. Suppose that there are chain maps
$i: H[1] \hookrightarrow V[1]$ (injective) and $p: V[1]\to H[1]$ (surjective) such that $p \circ i= id_{H[1]}$. Hence $V[1]$ splits as a graded vector space into $A[1]\oplus H[1]$. 
We assume existence of a homotopy $h: V[1] \to V$ such that
\begin{align*}
\delta \circ h + h \circ \delta = id_{V[1]} - i \circ p.
\end{align*}
The kernel of this map is a graded vector subspace of $V[1]$. We consider its intersection with $A[1]$ which we denote by $K[1]$.
The conormal bundle of $K[1]$ as a graded vector subspace of $A[1]$ is a Lagrangian vector subspace $L[1]$ of $T^{*}[1](A[1])$ and an isotropic subspace of $T^{*}[1](V[1])$.

Given a Lagrangian vector subspace $L[1]$ of $T^{*}[1](A[1])$ there is a well-defined notion of integration
\begin{align*}
\int_{L[1]}:\mathcal{C}^{\infty}(T^{*}[1](A[1]) \to \mathbb{R}
\end{align*}
under suitable convergence assumptions, see \cite{Schwarz}. The connection between the quantum master equation and the integration theory is

\begin{Theorem}\label{thm:Schwarz}
\begin{itemize}
\item Assume $S \in \mathcal{C}^{\infty}(T^{*}[1](A[1]))$ is $\Delta$-closed and let $L[1]$ and $L'[1]$ be two cobordant Lagrangian submanifolds of $T^{*}[1](A[1])$. Then $\int_{L[1]}S = \int_{L'[1]}S$.
\item Assume $S \in \mathcal{C}^{\infty}(T^{*}[1](A[1]))$ is $\Delta$-exact and let $L[1]$ be any Lagrangian submanifold of $T^{*}[1](A[1])$. Then $\int_{L[1]}S=0$.
\end{itemize}
\end{Theorem}

Using the splitting $V[1]=A[1]\oplus V[1]$ and the induced splitting of $T^{*}[1](V[1])$ one can extend $\int_{L[1]}$ to a map
\begin{align*}
\int_{L[1]}: \mathcal{C}^{\infty}(T^{*}[1](V[1])) \to \mathcal{C}^{\infty}(T^{*}[1](H[1])).
\end{align*}
Furthermore the BV-operator $\Delta$ also decomposes into $\Delta_{A} + \Delta_{H}$. Due to Theorem \ref{thm:Schwarz}, $\int_{L[1]}$ is a chain map between the complexes $(\mathcal{C}^{\infty}(T^{*}[1](V[1])),\Delta)$ and $(\mathcal{C}^{\infty}(T^{*}[1](H[1])),\Delta_{H[1]})$.

One can apply the BV-formalism as follows: interpret a vector field $Z$ on $V[1]$ as a function on $T^{*}[1](V[1])$ and assume that it satisfies the quantum master equation. Hence $e^{Z}$ is $\Delta$-closed. Apply the map $\int_{L[1]}$ to obtain a function $Y'$ on $T^{*}[1](H[1])$ that satisfies the quantum master equation with respect to $\Delta_{H}$. If one assumes that there is a function $Z'$ such that $e^{Z'}=Y'$ it follows that $Z'$ is a vector field that satisfies the quantum master equation.
This procedure has a physical intepretation in terms of integrating out ultraviolet degrees of freedom. Moreover there is a purely algebraic interpretation of the integration map $\int_{L[1]}$ in terms of certain graphs, known as Feynman diagrams.

It can be physically justified that in the ``classical limit'' the whole procedure reduces to the following: start with a cohomological vector field $X$ on $V[1]$, translate it to a function on $T^{*}[1](V[1])$. Using the tree-level part of the Feynman diagrams to ``integrate'' over
the isotropic subspace $L[1]$ one obtains a cohomological vector field on $H[1]$. If one reinterprets this in terms of $L_{\infty}[1]$-algebra structures one recovers the procedure for homological transfer along contractions as presented in subsection \ref{s:homotopytransfer}. 

Going beyond the tree-level in this integration procedure yields richer structures, see \cite{Costello} and \cite{Mnev} for instance.

\subsection{Transfer of differential complexes}\label{s:dcom}

\begin{Lemma}\label{l:dcom}
Let $(X,d,h,i,pr)$ be a graded vector space equipped with contraction data and a finite compatible filtration, i.e. a collection of graded subvector spaces
\begin{align*}
X=\mathcal{F}_{0}X \supseteq \mathcal{F}_{1}X \supseteq \cdots \supseteq \mathcal{F}_{n}X \supseteq \mathcal{F}_{(n+1)}X \supseteq \cdots
\end{align*}
such that $\mathcal{F}_{N}X = \{0\}$ for $N$ large enough, satisfying
\begin{itemize}
\item $d(\mathcal{F}_{k}X) \subset \mathcal{F}_{k}X$ for all $k \ge 0$ and
\item $h(\mathcal{F}_{k}X) \subset \mathcal{F}_{k}X$ for all $k \ge 0$.
\end{itemize}
Futhermore suppose $X$ is equipped with the structure of a differential complex $(X,D)$ such that
\begin{itemize}
\item $(D-d)(\mathcal{F}_{k}X) \subset \mathcal{F}_{(k+1)}X$.
\end{itemize}
Then the cohomology $H$ of $(X,d)$ is naturally equipped with the structure
of a differential complex and there is a well-defined chain map
$\tilde{i}:H \to X$.
\end{Lemma}

\begin{proof}
Set $D_{R}:=D-d$, it follows from $D^{2}=(d+D_{R})^{2}=0$ and $d^{2}=0$ that
\begin{align*}
D_{R}\circ d + d\circ D_{R} + D_{R}^{2}=0
\end{align*}
holds. In this special case the formulae for the induced structure given in subsection \ref{s:homotopytransfer} reduce to
\begin{eqnarray*}
\mathcal{D}&:=& p \circ \tilde{\mathcal{D}} \circ i \mbox{ where}\\
\tilde{\mathcal{D}}&:=& D_{R} \left(\sum_{k\ge 0}(-h D_{R})^{k}\right).
\end{eqnarray*}

\underline{claim 1:} $\mathcal{D} \circ \mathcal{D} =0$.
\newline 
We compute
\begin{eqnarray*}
-d(\tilde{\mathcal{D}})&=& D_{R} D_{R} (\sum_{m\ge 0}(-hD_{R})^{m}) + D_{R} d (\sum_{m\ge 0}(-hD_{R})^{m}) \\
&=& D_{R} D_{R} (\sum_{m\ge 0}(-hD_{R})^{m}) + D_{R} ip (D_{R} \circ \sum_{m\ge 0}(-hD_{R})^{m})\\
&& \quad -D_{R} D_{R} (\sum_{m\ge 0}(-hD_{R})^{m}) + D_{R} h d(D_{R} \circ \sum_{m\ge 0}(-hD_{R})^{m})\\
&=& \tilde{\mathcal{D}} ip \tilde{\mathcal{D}} + \tilde{\mathcal{D}} d
\end{eqnarray*}
and consequently
\begin{align*}
\mathcal{D}^{2}=p \tilde{\mathcal{D}} i \circ p \tilde{\mathcal{D}} i=0.
\end{align*}

The forumlae for the $L_{\infty}[1]$-morphism given in subsection \ref{s:homotopytransfer} reduce to
\begin{align*}
\tilde{i}:= (\sum_{k \ge 0} (-h D_{R})^{k}) i.
\end{align*}

\underline{claim 2:} $\tilde{i}$ is a chain map from $(H,\mathcal{D})$ to $(X,D)$.
\newline 
First we rewrite $\tilde{i}$ as
\begin{align*}
\tilde{i}=(id - h \tilde{\mathcal{D}})i
\end{align*}
and compute
\begin{eqnarray*}
D\circ \tilde{i} &=& (d + D_{R})(id - h \tilde{\mathcal{D}})i= d(-h \tilde{\mathcal{D}})i + \tilde{\mathcal{D}}i\\
&=& ip\tilde{\mathcal{D}}i + hd(\tilde{\mathcal{D}})i= (id - h \tilde{\mathcal{D}})i \circ p \tilde{\mathcal{D}}i
=\tilde{i} \circ \mathcal{D}.
\end{eqnarray*}
\end{proof}

\subsection{Transfer of differential graded Lie algebras}\label{s:dgLA}
We prove Theorem \ref{thm:transfer}, subsection \ref{s:homotopytransfer}:
We are given contraction data $(X,d,h,i,p)$ and the structure of a differential graded Lie algebra $(X,D,[-,-])$. In subsection \ref{s:dcom} we set $D_{R}:=D-d$ and defined $\tilde{\mathcal{D}}$ and $\mathcal{D}$ respectively.
We use the d\'ecalage-isomorphism to translate the graded Lie bracket into a graded symmetric operation which we denote by $\{-,-\}$ from now on.

The description of the induced structure maps can be rephrased as follows: consider an oriented trivalent tree $T$ with $n$ leaves whose edges are decorated by non-negative integers as introduced in subsection \ref{s:homotopytransfer}. One can associate a map
$\Phi(T):=\tilde{m}_{T}: (X[1])^{\otimes n} \to X[1]$ to $T$ by placing $\{-,-\}$ at its trivalent vertices, copies of $D_{R}$ at
all its edges, $-h$ between two consecutive such operations, $i$ at every leaf and $p$ at the root.
Let
\begin{align*}
\tilde{\nu}^{k}:=\sum_{\sigma \in \Sigma_{k}}\sum_{[T]\in[\mathbb{T}](k)}\frac{1}{|Aut(T)|}\tilde{m}_{T}
\end{align*}
and observe that $\nu^{k}$ from subsection \ref{s:homotopytransfer} coincides with $p \circ \tilde{\nu}^{k}\circ i^{\otimes k}$.

In subsection \ref{s:infty} we introduced the family of Jacobiators associated to a family of maps. By definition a family of maps constitutes an $L_{\infty}[1]$-algebra structure if the associated Jacobiators vanish.
Denote the family of Jacobiators associated to $(\nu^{k}:S^{k}(H[1]) \to H[2])$ by $(J^{n})$. We can write
$J^{n}=p\circ \tilde{J}^{n} \circ i^{\otimes n}$ with
\begin{multline*}
\tilde{J}^{n}(x_{1} \cdots x_{n}):=\\=
\sum_{r+s=n} \sum_{\sigma \in (r,s)-\text{shuffles}} \mspace{-36mu} \sign(\sigma)\, \tilde{\nu}^{s+1}(ip\tilde{\nu}^{r}(x_{\sigma(1)} \otimes \cdots \otimes x_{\sigma(r)})
\otimes x_{\sigma(r+1)} \otimes \cdots \otimes x_{\sigma(n)}). 
\end{multline*}

\underline{claim A:} $-d\left(\sum_{\sigma \in \Sigma_k}\sum_{[T] \in [[\mathbb{T}]](n)}\frac{1}{|Aut(T)|}\sigma^{*}\tilde{m}_{T}\right)i^{\otimes n}=\tilde{J}^{n}i^{\otimes n}$
\newline 

To prove this claim we introduce an extented graphical calculus: we allow to add one special edge in every tree which is marked either by a ``$\circ$'' or a ``$\times$'' and require that the special edge is decorated by two non-negative integers:

\begin{minipage}[c]{0.8 \textwidth}
\centering \includegraphics[width=2cm]{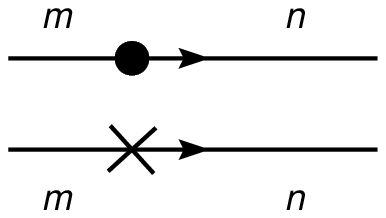}
\end{minipage}

We call oriented decorated trees with a special edge
of the first kind \textsf{pointed} and with a special edge of the second kind \textsf{truncated}.
Denote the space of pointed oriented decorated trees by $\mathbb{T}^{\circ}$ and the space of
truncated oriented decorated trees by $\mathbb{T}^{\times}$.
We extend $\Phi$ to trees with marked special edges: instead of composing two consecutive operations of degree $1$ by $\circ -h \circ$ we use
$\circ ip \circ$ at the pointed edge and ordinary composition
at the edge with a cross. Moreover one has to add the sign given by $(-1)$ to the powers of the sum of all
inputs left to the truncated or pointed edge.

One can easily check that
\begin{align*}
\sum_{\sigma \in \Sigma_{k}}\sum_{[T] \in [[\mathbb{T}]](n)}\frac{1}{|Aut(T)|}\sigma^{*}\Phi(P(T)) = \tilde{J}^{n}
\end{align*}
holds where $P(T)$ is the sum of all ways to change an ordinary edge of $T$ into a pointed one. Consequently claim A follows from

\underline{claim A.1:} 
\begin{align*}
-d\left(\sum_{\sigma \in \Sigma_{n}}\sum_{[T] \in [[\mathbb{T}]](n)}\frac{1}{|Aut(T)|}\sigma^{*}\tilde{m}_{T}\right)i^{\otimes n}= \left(\sum_{\sigma \in \Sigma_{k}}\sum_{[T] \in [[\mathbb{T}]](n)}\frac{1}{|Aut(T)|}\sigma^{*}\Phi(P(T))\right)i^{\otimes n}.
\end{align*}
\newline 

We prove claim A.1. by induction over the number of leaves $n$. For $n=1$ the claim is simply the equation
\begin{align*}
-d\tilde{\mathcal{D}}i = \tilde{\mathcal{D}}ip\tilde{\mathcal{D}}i
\end{align*}
which was established in subsection \ref{s:dcom}.
To inductive step uses the identities
\begin{eqnarray*}
-d \Phi(\includegraphics[width=1cm]{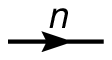}) =\mspace{-10mu}
\sum_{r+s=n+1}\mspace{-13mu} \Phi(\includegraphics[width=1cm]{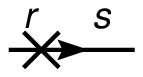})
-\sum_{r+s=n}\mspace{-13mu} \Phi(\includegraphics[width=1cm]{truncatedliners.eps}) +\sum_{r+s=n}\mspace{-13mu} \Phi(\includegraphics[width=1cm]{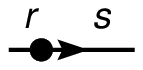}) + \Phi(\includegraphics[width=1cm]{simpleline.eps}) d
\end{eqnarray*}
and
\begin{eqnarray*}
-d\{X,Y\}&=&\{dX,Y\}+(-1)^{|X|}\{X,dY\}+\\
&+&\{D_{R}X,Y\}+(-1)^{|X|}\{X,D_{R}Y\}+D_{R}\{X,Y\}.
\end{eqnarray*}
Computing the left hand side of the equation in claim A.1. successivly leads to the right hand side plus
\begin{eqnarray*}
\sum_{\sigma \in \Sigma_{n}}\sum_{[T]\in [\mathbb{T}](n)}\frac{1}{|Aut(T)|}\sigma^{*}\Phi(X(T))
\end{eqnarray*}
where $X(T)$ is the sum of all ways to change an ordinary interior edge of $T$ into a truncated one which is decorated by $(0,0)$ . The evaluation of this sum at $x_{1}\otimes \cdots \otimes x_{n}$ contains terms of the form
\begin{eqnarray*}
&&\mspace{-20mu}\sum_{\sigma \in \Sigma_{n}} \text{sign}(\sigma) \sum_{r+s+t=n}
1/2 \sum_{[U] \in [\mathbb{T}](r), [V] \in [\mathbb{T}](s), [W] \in [\mathbb{T}](t)}
(\{\{-h \circ \Phi(U)(x_{\sigma(1)} \cdots x_{\sigma(r)}),\\
&&-h \circ \Phi(V)(x_{\sigma(r+1)} \cdots x_{\sigma(r+s)})\},-h \circ \Phi(W)(x_{\sigma(r+s+1)} \cdots x_{\sigma(n)})\}).
\end{eqnarray*}
Since the expression in the last two lines is of the form $\{\{a,b\},c\}$ and the sum runs
over all permutations with appropriate signs it vanishes due to the graded Jacobi identity.
Hence $J^{n}=p\tilde{J}^{n}i^{\otimes n}=p(d(\dots))=0$ and consequently the induced structure maps
$(\nu^{k}:S^{k}(H[1])\to H[2])$ define an $L_{\infty}[1]$-algebra structure on $H[1]$.

It remains to show that the maps $\lambda^{n}: S^{n}(H[1]) \to X[1]$ defined in subsection \ref{s:homotopytransfer} establish
an $L_{\infty}[1]$-morphism between $(H[1],\nu^{2},\nu^{2},\dots)$ and $(X[1],D,\{-,-\})$. We give explicite formulae for the identities that must be checked in order to prove that we obtain
an $L_{\infty}$-morphism:

\begin{eqnarray*}
&&-D(h \circ \tilde{\nu}^{n}(x_{1} \otimes \cdots \otimes x_{n}))\\
&&+1/2\sum_{r+s=n} \sum_{\sigma \in (r,s)-\text{shuffles}}\text{sign}(\sigma)\{h \circ \tilde{\nu}^{r}(x_{\sigma(1)}\otimes
\cdots \otimes x_{\sigma(r)}),\\
&& \hspace{5.5cm} h\circ \tilde{\nu}^{s}(x_{\sigma(r+1)}\otimes \cdots \otimes x_{\sigma(n)})\}\\
&&+\sum_{p+q=n}\sum_{\tau \in (q,p)-\text{shuffles}}\text{sign}(\tau) h \circ \tilde{\nu}^{p+1}(ip \circ
\tilde{\nu}^{q}(x_{\tau(1)}\otimes \cdots \otimes x_{\tau(q)})\otimes\\
&& \hspace{5.5cm} \otimes x_{\tau(q+1)} \otimes \cdots \otimes x_{\tau(n)})\\
&& -ip\tilde{\nu}^{n}(x_{1}\otimes \cdots \otimes x_{n})
\end{eqnarray*}
has to vanish identically for all $n \ge 2$ (the case $n=1$ was dealt with in subsection \ref{s:dcom}).

It is straightforward to check that the expression
\begin{itemize}
\item in the second and third line is equal to $B:=\tilde{\nu}^{n}+D_{R}h\tilde{\nu}^{n}$,
\item in the fourth line is equal to $C:=h\left(\sum_{\sigma \in \Sigma_{n}}\sum_{[T]\in [[\mathbb{T}]](n)}\frac{1}{|Aut(T)|}\sigma^{*}\Phi(P(T))\right)$,
\item in the first line is equal to $-\tilde{\nu}^{n}+ip\tilde{\nu}^{n} + hd\tilde{\nu}^{n}-D_{R}h\tilde{\nu}^{n}$.
\end{itemize}
The identity $-d\tilde{\nu}^{n}i^{\otimes i}=\left(\sum_{\sigma \in \Sigma_{n}}\sum_{[T]\in [[\mathbb{T}]](n)}\frac{1}{|Aut(T)|}\sigma^{*}\Phi(P(T))\right)i^{\otimes n}$ implies that everything cancels.